\magnification=1200

\loadmsam
\loadmsbm
\loadeufm
\loadeusm
\UseAMSsymbols

\font\BIGtitle=cmr10 scaled\magstep3
\font\bigtitle=cmr10 scaled\magstep1
\font\BIGtitle=cmr10 scaled\magstep3
\font\bigtitle=cmr10 scaled\magstep1
\font\boldsectionfont=cmb10 scaled\magstep1
\font\section=cmsy10 scaled\magstep1

\def\scr#1{{\fam\eusmfam\relax#1}}

\def\scrA{{\scr A}}
\def\scrB{{\scr B}}
\def\scrC{{\scr C}}
\def\scrD{{\scr D}}
\def\scrF{{\scr F}}

\def\scrH{{\scr H}}
\def\scrI{{\scr I}}
\def\scrL{{\scr L}}
\def\scrK{{\scr K}}
\def\scrJ{{\scr J}}
\def\scrM{{\scr M}}
\def\scrN{{\scr N}}

\def\scrP{{\scr P}}
\def\scrS{{\scr S}}

\def\scrT{{\scr T}}
\def\scrV{{\scr V}}

\def\scrW{{\scr W}}
\def\gr#1{{\fam\eufmfam\relax#1}}

\def\grC{{\gr C}}

	\def\grg{{\gr g}}

\def\grl{{\gr l}}

\def\db#1{{\fam\msbfam\relax#1}}

\def\dbA{{\db A}} 
\def\dbC{{\db C}} 
 \def\dbF{{\db F}}
\def\dbG{{\db G}} \def\dbH{{\db H}}

 \def\dbN{{\db N}}
 
\def\dbQ{{\db Q}} \def\dbR{{\db R}}

 \def\dbZ{{\db Z}}

\def\Ebar{\bar{E}}

\def\Gbar{\bar{G}}
\def\Hbar{\bar{H}}
\def\Kbar{\bar{K}}

\def\Xbar{\bar{X}}

\def\kbar{\bar{k}}

\def\vbar{\bar{v}}

\def\Gtil{\widetilde{G}}

\def\Rtil{\widetilde{R}}

\def\dbZhat{\widehat{\dbZ}}
\def\Ghat{\widehat{G}}

\def\Ker{\text{Ker}}
\def\der{\text{der}}
\def\Sh{\hbox{\rm Sh}}
\def\sh{\hbox{\rm sh}}
\def\Sp{\text{Sp}}

\def\sc{\text{sc}}
\def\Res{\text{Res}}
\def\ab{\text{ab}}
\def\ad{\text{ad}}

\def\Gal{\text{Gal}}

\def\End{\text{End}}
\def\Spec{\text{Spec}}
\def\Spf{\text{Spf}}

\def\Lie{\text{Lie}}

\def\leaderfill{\leaders\hbox to 1em
     {\hss.\hss}\hfill}
\def\nspace{\lineskip=1pt\baselineskip=12pt\lineskiplimit=0pt}

\def\Proclaim#1{\medbreak\noindent{\bf#1\enspace}\it\ignorespaces}
\def\finishproclaim{\par\rm
     \ifdim\lastskip<\medskipamount\removelastskip
     \penalty55\medskip\fi}
\def\proof{\par\noindent {\it Proof:}\enspace}
\def\references#1{\par
  \centerline{\boldsectionfont References}\bigskip
     \parindent=#1pt\nspace}
\def\Ref[#1]{\par\medskip\hang\indent\llap{\hbox to\parindent
     {[#1]\hfil\enspace}}\ignorespaces}
\def\Item#1{\par\smallskip\hang\indent\llap{\hbox to\parindent
     {#1\hfill$\,\,$}}\ignorespaces}
\def\ItemItem#1{\par\indent\hangindent2\parindent
     \hbox to \parindent{#1\hfill\enspace}\ignorespaces}

\def\Le{{\mathchoice{\,{\scriptstyle\le}\,}
  {\,{\scriptstyle\le}\,}
  {\,{\scriptscriptstyle\le}\,}{\,{\scriptscriptstyle\le}\,}}}
\def\Ge{{\mathchoice{\,{\scriptstyle\ge}\,}
  {\,{\scriptstyle\ge}\,}
  {\,{\scriptscriptstyle\ge}\,}{\,{\scriptscriptstyle\ge}\,}}}

\def\arrowsim{\,\smash{\mathop{\to}\limits^{\lower1.5pt
  \hbox{$\scriptstyle\sim$}}}\,}

\def\doublemaprights#1#2#3#4{\raise3pt\hbox{$\mathop{\,\,\hbox to
     #1pt{\rightarrowfill}\kern-30pt\lower3.95pt\hbox to
     #2pt{\rightarrowfill}\,\,}\limits_{#3}^{#4}$}}

\def\rightcapdownarrow{\raise9pt\hbox{$\ssize\cap$}\kern-7.75pt
     \Big\downarrow}

\def\rcapmapdown#1{\rightcapdownarrow\kern-1.0pt\vcenter{
     \hbox{$\scriptstyle#1$}}}

\def\rmapdown#1{\Big\downarrow\kern-1.0pt\vcenter{
     \hbox{$\scriptstyle#1$}}}
\def\rightsubsetarrow#1{{\ssize\subset}\kern-4.5pt\lower2.85pt
     \hbox to #1pt{\rightarrowfill}}
\def\longtwoheadedrightarrow#1{\raise2.2pt\hbox to #1pt{\hrulefill}
     \!\!\!\twoheadrightarrow}

\def\Gal{\operatorname{\hbox{Gal}}}

\NoBlackBoxes
\parindent=25pt
\document
\footline={\hfil}

\null
\vskip 0.6 cm
\centerline{\BIGtitle Shimura Varieties and the Mumford--Tate conjecture,}
\centerline{\BIGtitle part I}
\vskip 0.2in
\centerline{\bigtitle Adrian Vasiu}
\vskip 0.2 in
\centerline{Univ. of Arizona}
\centerline{vers. 9/14/01: few enlargements, improvements and corrections to the vers. dated 1/14/00}
\centerline{to be resubmitted for public. in Ann. of Math.} 
\footline={\hfill}
\null
\vskip 0.05 in
\centerline{ABSTRACT. We prove the Mumford--Tate conjecture for those abelian varieties}
\centerline{over number fields, whose simple factors of their adjoint Mumford--Tate groups have}
\centerline{over $\dbR$ certain (products of) non-compact factors. In particular, we prove this conjecture}  
\centerline{for the orthogonal case ($B_n$ and $D_n^{\dbR}$ Shimura types). We construct embeddings of Shimura}
\centerline{varieties of (whose adjoints are of prescribed abelian type) into unitary Shimura varieties. We}
\centerline{also prove two general theorems involving Frobenius tori of abelian varieties over number fields.} 
\smallskip
\centerline{Key words: abelian and Shimura varieties, reductive and p-divisible groups,}
\centerline{Galois representations and Hodge cycles.}
\smallskip
\centerline{MSC 2000: Primary 11G10, 11G18, 11R32, 14G35 and 14G40.}

\vskip 0.2 in
\centerline{\bigtitle Contents}

{\nspace{

\smallskip
\line{\item{\bf\S1.}{Introduction}\leaderfill 1}

\smallskip
\line{\item{\bf \S2.}{Some complements} \leaderfill 9}

\smallskip
\line{\item{\bf \S3.}{Good embeddings into unitary Shimura varieties} \leaderfill 19}

\smallskip
\line{\item{\bf \S4.}{The basic techniques} \leaderfill 35}

\smallskip
\line{\item{\bf \S5.}{Some conclusions} \leaderfill 58}

\smallskip
\line{\item{}{Erratum to [Va1]} \leaderfill 76}

\smallskip
\line{\item{}{References}\leaderfill 77}

}}

\footline={\hss\tenrm \folio\hss}
\pageno=1

\bigskip\smallskip
\noindent
{\boldsectionfont \S1. Introduction}
\bigskip

\Proclaim{1.1. An approach.} \rm
Let $A$ be an abelian variety over a number field $E$. Let $p$ be a rational prime. Let $H_A$ (resp. $G_{\dbQ_p}$) be the Mumford--Tate group of $A$ (resp. be the connected component of the origin of the algebraic envelope of the $p$-adic representation naturally attached to $A$), cf. 4.0. $H_{A\dbQ_p}$ is naturally identified with a subgroup of $GL(H^1_{\text{\'et}}(A_{\Ebar},\dbQ_p))$ (see 4.0). The Mumford--Tate conjecture (see 4.1) asserts: as subgroups of $GL(H^1_{\text{\'et}}(A_{\Ebar},\dbQ_p))$, we have $H_{A\dbQ_p}=G_{\dbQ_p}$. Our approach to the proof of this conjecture is via the following five main steps.

\medskip
{\bf S1.} A well known result (see 4.2.5) allows us to assume that $p$ is as big as desired. In particular, we can assume the existence of a prime $v$ of $E$ unramified over $p$ and such that $A$ has good reduction with respect to it.

\smallskip
{\bf S2.} Based on standard arguments on Frobenius tori (for their construction reviewed in 4.2.7.1, see [Chi] and [Pi]), we can assume $G_{\dbQ_p}$ is a split, reductive group (see 4.2.8).

\smallskip
{\bf S3.} Using the main result of [Wi] we define each simple factor $\scrF$ of $G_{\dbQ_p}^{\ad}$ or of $H_{A\dbQ_p}^{\ad}$ to be compact or non-compact, depending on the fact that a suitable injective cocharacter 
$$\mu\colon\dbG_m\hookrightarrow G_{\dbQ_p}$$ 
(see 4.3.3 and 4.4.4) has a trivial or a non-trivial image in $\scrF$ (for details see 5.1.3). The third step is an attempt to show that the product of the non-compact factors of $G^{\ad}_{\dbQ_p}$ is isomorphic to the product of the non-compact factors of $H^{\ad}_{A\dbQ_p}$.

\smallskip
{\bf S4.} Assuming that the attempt of S3 has been successful, we show that in fact $G^{\der}_{\dbQ_p}=H^{\der}_{A\dbQ_p}$.

\smallskip
{\bf S5.} Assuming $G^{\der}_{\dbQ_p}=H^{\der}_{A\dbQ_p}$, by standard arguments we show that $G_{\dbQ_p}=H_{A\dbQ_p}$.
\finishproclaim

Combining S1-5, we stated in [Va2, 1.15.4] the following expectation:

\medskip
\Proclaim{1.2. Expectation (the split criterion).}
There is $N(A)\in\dbN$ such that: 
if $p\Ge N(A)$ and if $G_{\dbQ_p}$ is a split group, then $G_{\dbQ_p}=H_{A\dbQ_p}$.
\finishproclaim

It is known (see 4.2.9) that this Expectation implies $G_{\dbQ_p}=H_{A\dbQ_p}$, for any prime $p$.

\medskip
\Proclaim{1.3. Conventions, notations and organization.} \rm
For types and generalities on (integral aspects of) Shimura varieties we refer to [Va1, 2.2-5 and \S 3]. For standard Hodge situations and standard PEL situations we refer to [Va2, 2.3.5]. The expression ``standard Hodge situation" is abbreviated as SHS, cf. loc. cit. Warning: we use freely [Va1, 4.3.4] and the notations and conventions of [Va2, 2.1].  
A simple, adjoint Shimura pair $(G,X)$ is said to be without involution (resp. with involution) if its reflex field $E(G,X)$ is (resp. is not) a totally real number field, cf. [De1, 3.8 (ii)]. For Shimura group pairs we refer to [Va2, 2.2.5 1), 2.2.6 5) and end of 4.1]. 

Always by $k$ we denote a perfect field of characteristic $p\Ge 2$. $W(k)$ is the Witt ring of $k$, $B(k):=W(k)\fracwithdelims[]1p$, while $\sigma$ is the Frobenius automorphism of any such $k$, $W(k)$ or $B(k)$. 
For the notions of (principally quasi-polarized) Shimura (filtered) $\sigma$-crystals or $p$-divisible groups we refer to [Va2, 2.2.8, 2.2.9 12) and 2.2.20]; for these notions in the context of a SHS, we refer to [Va2, 2.3.9-11]. For $p$-divisible objects (with tensors) of Fontaine's categories we refer to [Va2, 2.2.1 c)]. The notion of a Shimura $p$-divisible group (over a regular, formally smooth $W(k)$-scheme) being a versal deformation (at some maximal point which is algebraic over $k$), is understood in the sense of [Va2, 3.6.19 B)]. We refer to loc. cit. (resp. to [Va2, 2.2.1.3]) for uni plus quasi-versal (resp. for pull backs of) ($p$-divisible) objects of Fontaine's categories.

For $n\in\dbN$, the expression ``$n$ is big enough" is used in the sense $n>>0$ but with effective lower bounds. In similar senses (involving degrees or suitable Haar measures), we speak about large enough number fields and small enough compact, open subgroups of locally compact groups.

We denote by $E_v$ the completion of a number field $E$ w.r.t. a finite prime of it $v$. If $G$ is a group acting on a set $S$, we denote by $S^G$ the subset of $S$ formed by elements fixed by $G$. If $M$ is a $*$-module, with $*$ a group scheme or a Lie algebra over some ring, we speak as well about the representation of $*$ on $M$. If $*$ is a reductive group over an integral scheme, then $Z^0(*)$ denotes the maximal torus of the center $Z(*)$ of $*$.

If $X_S$ or $X$ is a scheme over a scheme (resp. over a $\dbZ$-algebra) $S$, then for any other $S$-scheme $T$ (resp. for any homomorphism $i_S:S\to T$ of $\dbZ$-algebras), we often denote, without any extra comments, by $X_T$ the $S$-product (resp. the $\Spec(S)$-product) of $X_S$ or $X$ and $T$ (resp. and $\Spec(T)$). The same applies for morphisms. Warning: if there are a couple of ways to take $i_S$, we always mention it a priori but we never add it as part of the notation $X_T$. We use freely the standard notations (see [Bou2, p. 188]) for the fundamental weights of a split, simple Lie algebra of $A_n$ Lie type. We use a lower left index $B$ to denote the Betti cohomology and homology. 
\finishproclaim                 

\Proclaim{1.3.1. On techniques.} \rm
In [Va1-2] we developed some techniques pertaining to Shimura varieties of preabelian type and to the deformation theory of abelian varieties endowed with de Rham (crystalline) tensors. Here (to be continued in [Va5] and [Va7]) we start exploiting them for achieving progress in the proof of the Mumford--Tate conjecture. In this paper, we concentrate mostly in using techniques from reductive groups (see [Va1, 4.3] and 4.3.3-4), from the local deformation theory (see [Va1, 5.2-4]), and from Shimura varieties (see [Va1, 6.5-6] and \S3). In [Va5] (resp. in [Va7]) we will enlarge these techniques by the ones of [LP1] and by extra techniques pertaining to Shimura varieties (resp. by the ones involving -- see [Va2, 3.6 and 3.15.1]-- global deformations and the ones involved --see the announcements of [Va1, 1.7.1] and [Va7]-- in our work on the Langlands--Rapoport conjecture). We do hope that all these techniques and ideas will match to a complete proof of the conjecture.
 
For samples of what the new techniques presented here can achieve see 1.4. These techniques are elaborated in \S3 (the ones involving Shimura varieties) and in \S4 (the ones involving reductive groups and the local deformation theory). For details of how the paper is structured see 1.5-8.
\finishproclaim

\Proclaim{1.3.2. The context.} \rm
We consider the situation of 1.2: $G_{\dbQ_p}$ is split and $p\Ge N(A)$. Let $(v_{\alpha})_{\alpha\in\scrJ}$ be a family of homogeneous tensors of (see [Va2, 2.1]) the essential tensor algebra 
$\scrT(H^1_{\text{\'et}}(A_{\Ebar},\dbQ_p))$ of $H_{\text{\'et}}^1(A_{\Ebar},\dbQ_p)\oplus H_{\text{\'et}}^1(A_{\Ebar},\dbQ_p)^*$ such that $G_{\dbQ_p}$ is the subgroup of $GL(H_{\text{\'et}}^1(A_{\Ebar},\dbQ_p))$ fixing $v_{\alpha}$, $\forall\alpha\in\scrJ$. Such a family exists (cf. [De3, 3.1 c)]). For our proposed number $N(A)\in\dbN$ see 4.3.0; the main new ideas behind its definition are based on the philosophy of 4.3.4.1. Here we just mention that $N(A)$ depends only on $A$ and is effectively (and easily) computable. Below we take $p$ big enough and odd and such that we are as well in the context of S1. For instance, $p\Ge N(A)$ works for all that follows; warning: each isolated part of 1.4-8 can be worked out under better bounds. Below it is irrelevant (cf. 4.0) which complex embeddings we choose for different number fields and so we do not mention them.
\finishproclaim

\Proclaim{1.4. The main results.} \rm
To state some of the results of this paper, let $\Sh(H_A,X_A)$ be the Shimura variety attached to $A$. We write the adjoint Shimura pair
$$(H_A^{\ad},X_A^{\ad})=\prod_{j\in\scrI} (H_j,X_j)\leqno (1)$$
as a finite product of simple, adjoint Shimura pairs. We have (see 5.1.2):
\finishproclaim

\smallskip
\Proclaim{Theorem 1.}
The split criterion holds if $\forall j\in\scrI$, the adjoint Shimura pair $(H_j,X_j)$ is one of the following types:

\medskip
{\bf a)} non-special $A_n$ type;

\smallskip
{\bf b)} $B_n$ type;

\smallskip
{\bf c)} $C_n$ type, with $n$ odd;

\smallskip
{\bf c')} $C_{2n}$ type, with $n$ odd and such that $4n$ is not of the form $C_{4m+2}^{2m+1}$, with $m\in\dbN$;

\smallskip
{\bf d)} $D_n^{\dbH}$ type, with $n$ odd and such that $2n$ is not of the form $C_{2^{m+1}}^{2^m}$, with $m\in\dbN$;

\smallskip
{\bf d')} non-inner $D_{2n}^{\dbH}$ type, with $n\in\{q,3q|q\,\text{a prime}\}\setminus\{9\}$; 

\smallskip
{\bf e)} $D_n^{\dbR}$ type.
\finishproclaim

\smallskip
Here we do not describe in detail what means non-special $A_n$ type (see 5.1.2.1 2)). We just mention four large families of cases of such types. 

\Proclaim{Samples (see 5.1.2.3).} 
$(H_j,X_j)$ is of non-special $A_n$ type if any one of the following four conditions is satisfied:

\medskip
{\bf i)} $H_{j\dbR}$ has two simple factors $SU(a_i,b_i)^{\ad}_{\dbR}$, $i=\overline{1,2}$, such that $1\Le a_1< a_2\le b_2< b_1$ and either $(a_1,b_1)=1$ or $(a_2,b_2)=1$;

\smallskip
{\bf ii)} $H_{j\dbR}$ has a simple, non-compact factor $SU(a,b)^{\ad}_{\dbR}$, with $g.c.d.(a,b)=1$ with the pair $(a,b)$ not of the form $(C_r^s,C_r^{s-1})$, for $s,r\in\dbN$, $2\le s\le r-1$;

\smallskip
{\bf iii)} $H_{j\dbR}$ has as a factor a product $SU(a,a)^{\ad}_{\dbR}\times SU(b,2a-b)_{\dbR}^{\ad}$ with $b$ different from $a$ and odd;

\smallskip
{\bf iv)} $H_j$ is of $A_{q-1}$ Lie type, with $q$ equal to $4$ or to a prime.
\finishproclaim

See 2.5 for the notion of $(H_j,X_j)$ being of non-inner type. Sample (cf. end of 2.5.1): $(H_j,X_j)$ is of non-inner $D_{2n}^{\dbH}$ if it is of $D_{2n}^{\dbH}$ with $\dbR$-involution in the sense of 2.5.1 b). Let $A_v$ be the reduction of $A$ w.r.t. $v$. Let $T_v$ be the Frobenius torus of $A_v$. Let $l$ be a rational prime different from $p$. Let $G_{\dbQ_l}$ have a meaning analogue to $G_{\dbQ_p}$ but for $l$. We have:

\smallskip
\Proclaim{Theorem 2 (see 3.3.2).}
Let $j\in\scrI$. We view $T_{v\dbQ_l}$ as a subgroup of $G_{\dbQ_l}$ (uniquely determined up to $G_{\dbQ_l}$-conjugacy). We have: all images of ${T_v}_{\overline{\dbQ_l}}$ in simple factors of $H_{j\overline{\dbQ_l}}$ have the same dimension.
\finishproclaim

\Proclaim{Theorem 3 (see 4.4.8 b) and c)).}
We assume $A_v$ is an ordinary abelian variety. Then $T_v$ is isomorphic to a torus of $H_A$ and $T_{v\dbQ_p}$ is isomorphic to a torus of $G_{\dbQ_p}$. Moreover, for any rational prime $q$, $Z^0(H_{A\dbQ_q})$ is a torus of $G_{\dbQ_q}$.
\finishproclaim

\Proclaim{1.5. Some complements to [Va2].} \rm
Our work on 1.2 grew out (see [Va2, 3.6.20 2)]) of Milne's conjecture ([Mi4, 0.1]; see also [Va1, 5.6.6]). A significant part of [Va7] is dedicated to this conjecture; so in 2.2 we first state it as a question in a more general context and then prove it in two generic situations (see 2.2.1-2) as well as in a special case (see 2.2.3), just to pay due homage to it and to be able to apply it in \S3. Extra complements to [Va1-2] are included in 2.3-5. We draw attention to 2.4-5. 2.4 complements the relative form [Va2, 3.6.18.7] of the inducing property of [Va2, 3.6.18.5], over a complete, ramified DVR $V$ of mixed characteristic and residue field $k$. 2.5 introduces some extra fields (accompanying reflex fields) and so implicitly types of simple, adjoint Shimura pairs.
\finishproclaim

\Proclaim{1.6. Embeddings into Siegel modular varieties and the shifting process.} \rm 
In 3.1-2 we prove a rational version of the result stated in [Va2, 1.15.8 1) to 4)]. It is just the $\dbQ$--version of the results (over $\dbR$) of [Sa], obtained following the pattern of [Va1, 6.5-6] and Erratum to [Va1]. In other words, starting from a Shimura pair $(G,X)$ of adjoint, abelian type, we construct (see 3.2 for the case when $(G,X)$ is simple and see 3.2.2 2) for the passage to the general case via the Hodge quasi products of 2.6) useful injective maps 
$$(G_0,X_0)\hookrightarrow (G_1,X_1)\operatornamewithlimits{\hookrightarrow}^{f_1} (GSp(W,\psi),S),$$
with $(G_0^{\ad},X_0^{\ad})=(G,X)$, with $f_1$ a PEL type embedding (see [Va2, 2.1] for the meaning of this; cf. also [Sh] and [Ko]), and with $G_1^{\der}$ a simply connected semisimple group whose adjoint has all simple factors of some $A_n$ Lie type ($n\in\dbN$ can depend on the factor). This rational version of [Va2, 1.15.8] is good enough to get as well an effective generic version of loc. cit. (see 3.2.1 and 3.2.2 3)). 3.2.1 has many applications (see 1.6.1 and 3.2.2 4) for a glimpse): as part of its proof is slightly sketched, they are not stated here. In some situations it is more practical to work with uniformized variants of 3.2.1 3); they are mentioned in 3.3.3.1. 

The proof of 3.2 relies on [De2, 2.3.10], being as well the (refined) unitary version of loc. cit., and on a natural twisting process (see 3.1 F)) of certain semisimple groups over number fields. The passage from [Sa] to [De2, 2.3.10] ``forgets" the fact that the embeddings of hermitian symmetric domains of classical Lie type into Siegel domains constructed in [Sa] are obtained by first achieving embeddings into hermitian symmetric domains associated to unitary groups. The mentioned twisting process is the extra ingredient needed to ``correct" this forgetfulness. 
\finishproclaim

\Proclaim{1.6.1. More on 3.3.} \rm
In 3.3.1-2 we prove Theorem 2 (under stronger conditions which just require $p$ big enough). The proof relies on 3.2.1 2), on standard arguments involving Frobenius tori and on the fact that the Langlands--Rapoport conjecture is known to be true (see [Mi4]) for Shimura varieties of $A_n$ type. In fact, for what we need from this last conjecture, the main result of [Zi] suffices. 

Theorem 2, though it is used just in 5.1.3-4 beginning with the Step 2 of 5.1.3, it is stated and proved in 3.3 for not interfering with the flow of ideas of \S4-5.
For it, we assume that the reader is familiar with some of the results reviewed in parts running from 4.0 to 4.2.7.1 and (very isolated) with 4.2.8, 4.2.11 and 4.2.14. So the reader less familiar with these reviewed results, should read 3.3 only after reading the mentioned parts of 4.0-2.

The shifting process is stated explicitly in 3.3.3. It can be summarized as: from the point of view of S3 or S4, we have total freedom to ``cook" an embedding into a Siegel modular variety of a Shimura variety whose adjoint variety is isomorphic to $\Sh(H_A^{\ad},X_A^{\ad})$. See 3.3.3 ADV for some advantages we get using a suitable such ``cooking"(which replaces $A$ by another abelian variety $A^{1}$ over a finite field extension $E^1$ of $E$,  unrelated to $A$ or $A_{E^1}$ in any usual way, like isogenies, Weil restrictions of scalars, etc.; however, it is appropriate to say that $A^1$ and $A_{E^1}$ are adjoint-isogeneous).

The ideas of 3.3 can be used in many other situations; in particular, they reobtain (see 3.3.2.1) the results of [Pi, 5.10-11] (recalled in 4.2.14), provided we know that the ordinary reduction conjecture (for instance, see  [Oo, p. 11]) is true for $A$. 
\finishproclaim

\Proclaim{1.7. On Theorem 1.} \rm
Most of the proof of Theorem 1 is standard (to be compared with [Pi, 5.14-15]). However, it involves three main new ideas. First, in the $A_n$ type case, besides Serre's result (see [Pi, 3.5]) and Pink's results (see [Pi, 5.10-11]), we rely on 3.2.1 2) and [Zi]. In other words, if all Shimura pairs $(H_j,X_j)$ ($j\in\scrI$) are of some $A_n$ type, then using 3.3.1 2) and the main result of [Zi] we check directly that we can assume (in connection to S3 or S4) that $T_v$ is isomorphic to a torus of an inner form of $H_A$ (we recall $p$ is big enough).

Second, we do not need assumptions on (the dimension or the structure of) $H_A^{\ab}$, as [Wi] and the simple Fact 4.2.7.2 make them unnecessary. They achieve S5, cf. Step 3 of 5.1.3. See also Exercise 5.1.7 for the stronger fact (generalizing [Bog]) that always $Z^0(H_A)_{\dbQ_p}$ is a torus of $G_{\dbQ_p}$.

Third, we use more intimately 3.3.3.1. More precisely: Theorem 2, the shifting process of 3.3.3 and 3.3.3.1 are the main new ingredients needed to accomplish S4 (via Faltings' well known result 4.2.3 and Zarhin's lemma stated in 4.2.5.1); see Step 2 of 5.1.3. We do not stop here to detail fully how these ingredients get combined in the general case. However, as an exemplification of their usage, we assume for the next 2 paragraphs that $\scrI$ has just one element and that $(H_A^{\ad},X_A^{\ad})$ is as in e) of Theorem 1. 

It is known (cf. [Pi, 4.3 and 5.10]) that the image $\scrF_1$ of $G_{\overline{\dbQ_p}}$ in a simple factor $\scrF_0$ of $H_{A\overline{\dbQ_p}}^{\ad}$ is either $\scrF_0$ itself or is a suitable product $PB$ of two simple, adjoint groups of some $B_m$ Lie type, $m\in\{0,...,n-1\}$; the $B_0$ Lie type corresponds to the trivial group. The shifting process and 3.3.3 1) allow us to assume that the representation of $\Lie(\scrF_0)$ on $H^1:=H^1_{\text{\'et}}(A_{\overline{\dbQ_p}},\dbQ_p)\otimes_{\dbQ_p} \overline{\dbQ_p}$ involves both the half spin representations (of the split $D_n$ Lie type) and that any non-trivial, irreducible subrepresentation of $\Lie(H_{A\overline{\dbQ_p}}^{\der})$ on $H^1$ factors through a simple factor of it (i.e. no tensor product representations of products of simple $\overline{\dbQ_p}$-Lie algebras are involved). The representation of $\Lie(PB)$ on a non-trivial, simple $\Lie(\scrF_0)$-submodule $H^1_\prime$ of $H^1$ which is as well a $H_{A\overline{\dbQ_p}}$-module, is (up to isomorphism) the same regardless of which half spin representation is involved in defining $H^1_\prime$. 

Moreover, based on 3.3.3 1) and 2) we can assume that, for any arbitrarily chosen such $\Lie(\scrF_0)$-submodule $H^1_\prime$, the representation on it of $Z^0(H_{A\overline{\dbQ_p}})$ is the same scalar representation as for a similar $\Lie(\scrF_0)$-submodule $H^1_{\prime\prime}$ of $H^1$ but involving the other half spin representation. So if $\scrF_1$ is not $\scrF_0$, we reach a contradiction with 4.2.3.
If no simple factor of $H_{A\dbR}^{\ad}$ is compact, an entirely similar argument (based on 3.3.3 3)) shows that $\Lie(\scrF_0)$ can not project onto two or more distinct simple factors of $\Lie(H^{\ad}_{A\overline{\dbQ_p}})$. 

Warning: in general, if $(H_j,X_j)$ is of non-special $A_n$ type with involution or if $H_{j\dbR}$ has compact factors, the above first and third ideas have to be intimately combined (via 3.3.3.1) in order to get that a simple factor of $\Lie(G^{\ad}_{\overline{\dbQ_p}})$ does not project onto two or more simple factors of $\Lie(H^{\ad}_{A\overline{\dbQ_p}})$ (see Step 2 of 5.1.3). 
\finishproclaim

\Proclaim{1.7.1. Comments.} \rm
To our knowledge, a) and e) of Theorem 1 are completely new general instances for the validity of the conjecture. For instance, e) was not known before even in the case when $\End(A_{\dbC})=\dbZ$. Very particular cases of a), b), c), c'), d) and d') of Theorem 1, but involving Shimura varieties of arbitrary dimension, were known previously; see 4.2.16 for references. From the point of view of possibilities of types of $(H_j^{\ad},X_j^{\ad})$, Theorem 1 handles about 80 percent; this percentage is computed as follows: 20 percent for each b) and e), 16 percent for a), 14 percent for c) and c'), and just 10 percent for d) and d'). See 5.1.2.4 for the simplest cases not settled by Theorem 1. 
\finishproclaim

\Proclaim{1.8. On Theorem 3.} \rm
To explain the ideas behind Theorem 3, we need some notations.
 Let $O_v$ be the ring of integers of $E_v$. So $W(\dbF)$ is the completion of the maximal unramified, finite extension of $O_v$. It is quite convenient (from notations' point of view) to work what follows over $W(\dbF)$ instead of over $O_v$ (however, see 4.4.7.2). 

After some nice arrangements (cf. most of 4.3) we can put everything in the context of a Shimura filtered $\sigma$-crystal
$$\scrC:=(M,F^1,\phi,G_{W(\dbF)},(t_{\alpha})_{\alpha\in\scrJ})$$ 
over $\dbF$. $M$, $F^1$ and $\phi$ are as explained below. For $\alpha\in\scrJ$, $t_{\alpha}$ is the homogeneous tensor of $\scrT(M\fracwithdelims[]1p)$ corresponding to $v_{\alpha}$ under Fontaine's comparison theory. $G_{W(\dbF)}$ is the Zariski closure in $GL(M)$ of the subgroup of $GL(M\fracwithdelims[]1p)$ fixing $t_{\alpha}$, $\forall\alpha\in\scrJ$. Here is the central place where we need $p$ to be big enough: for these arrangements we use (see 4.3.4) the theory of $\dbZ_p$-well positioned family of tensors of [Va1, 4.3]; we mainly need [Va1, 4.3.10.1 1)]. 

To reach to the context of $\scrC$ we need to pass to a $\dbZ\fracwithdelims[]1p$-isogeny $A^{\prime}$ of $A_{E^{\prime}}$ (with $E^{\prime}$ a finite field extension of $E$) and to a prime $v^{\prime}$ of $E^{\prime}$ dividing $v$, cf. 4.3.3 and 4.3.5-6. What we achieve by doing this: the Zariski closure $G_{\dbZ_p}$ of $G_{\dbQ_p}$ in $GL(H^1_{\text{\'et}}(A^\prime_{\overline{E^\prime}},\dbZ_p))$ is a (split) reductive group over $\dbZ_p$.  

If $v^{\prime}$ is unramified over $p$, then  $M:=H_{\text{crys}}^1(A_{\dbF}^{\prime}/W(\dbF)),$ with $A_{\dbF}^{\prime}$ as the special fibre of the (abelian scheme which is the) N\'eron model $A_{W(\dbF)}^{\prime}$ of $A^{\prime}_{B(\dbF)}$, $\phi$ is the resulting $\sigma$-linear endomorphism of $M$ and $F^1$ is the Hodge filtration of $M$ defined by $A_{W(\dbF)}^{\prime}$. 

If $v^{\prime}$ is ramified over $p$, then the filtered $\sigma$-crystal $(M,F^1,\phi)$ is still the one of an abelian scheme $A^0$ over $W(\dbF)$; $A^0_{\dbF}$ is the extension to $\dbF$ of the reduction of $A^{\prime}$ w.r.t. $v^{\prime}$. To explain how we get $A^0$, let $L$ be the completion of the maximal unramified, algebraic field extension of $E_{v^\prime}^\prime$ and let $O_L$ be its ring of integers. Let $e_A:=[O_L:W(\dbF)]$. Let $\Rtil e_A$ be the piano ring of $\dbF$ of resonance $e_A$ (see 2.1).
Let $A^{\prime}_{O_L}$ be the abelian scheme over $O_L$ which is the N\'eron model of $A_L^{\prime}$.
 As $G_{\dbZ_p}$ is a split, reductive group and as $p$ is big enough, we can assume (see 4.3.4) that the family of tensors $(v_{\alpha})_{\alpha\in\scrJ}$ has a subfamily which is $\dbZ_p$-very well position w.r.t. $G_{\dbQ_p}$ and is enveloped by $H^1_{\text{\'et}}(A^{\prime}_{\overline{L}},\dbZ_p)$. So the machinery of [Va1, 5.2-3] applies to $A_{O_L}^{\prime}$. We get (see 4.4.5; 1.8.1 points out that we can assume that we are in a polarized context) an abelian scheme $A_{\Spec(\Rtil e_A)}$ over $\Spec(\Rtil e_A)$ such that: 

\medskip
{\bf a)} its pull back through an $O_L$-valued point of $\Spec(\Rtil e_A)$ is $A^\prime_{O_L}$;

\smallskip
{\bf b)} the family $(w_{\alpha}^\prime)_{\alpha\in\scrJ}$ of tensors corresponding to $(v_{\alpha})_{\alpha\in\scrJ}$ under Fontaine's comparison theory, belongs to the $F^1$-filtration of $\scrT(H^1_{dR}(A_{\Spec(\Rtil e_A)}/\Rtil e_A)[{1\over p}])$ defined by the Hodge filtration of $A_{\Spec(\Rtil e_A[{1\over p}])}$. 

\medskip
We take $(A^0,(t_{\alpha})_{\alpha\in\scrJ}):=z_0^*(A_{\Spec(\Rtil e_A)},(w_{\alpha}^\prime)_{\alpha\in\scrJ})$, where $z_0$ is the $W(\dbF)$-valued Teichm\"uller lift of $\Spec(\Rtil e_A)$; see 2.1 for the Frobenius lift of $\Spec(\Rtil e_A)$. 
\finishproclaim

\Proclaim{1.8.1. Principal polarizations.} \rm
Regardless of how $v^{\prime}$ is over $p$, we can assume (this is standard, see 4.3.2-3 and 4.3.5-6)  that all abelian schemes considered above (i.e. $A$, $A^\prime$, $A_{\Spec(\Rtil e_A)}$ and $A^0$) are principally polarized and have some level-$N$ symplectic similitude structure (for $N\in\dbN\setminus\{1,2\}$ and relatively prime to $p$; as $p\Ge 3$, we take $N=4$). We can also assume that the $\dbZ\fracwithdelims[]1p$-isogeny between $A_{E^\prime}$ and $A^\prime$ is compatible with the chosen polarizations.

If $v^{\prime}$ is unramified over $p$, we denote $A_{W(\dbF)}^{\prime}$ from now on  by $A^0$. So $A^0$ is well defined regardless of how $v^\prime$ is over $p$. 
Let $\psi_M\colon M\otimes_{W(\dbF)} M\to W(\dbF)(1)$ be the perfect alternating form defined by the principal polarization of $A^0$. We get a principally quasi-polarized Shimura filtered $\sigma$-crystal 
$$\scrC\scrP:=(M,F^1,\phi,G_{W(k)},(t_{\alpha})_{\alpha\in\scrJ},\psi_M).$$ 
\finishproclaim

\Proclaim{1.8.2. Versal deformation.} \rm
Next we just apply the machinery of [Va1, 5.4]. We get (see 4.4.7) a versal deformation of $\scrC\scrP$ over a regular ring $R$ of formal power series with coefficients in $W(\dbF)$; it corresponds to a principally polarized abelian scheme $(A_R^0,p_{A_R^0})$ over $\Spec(R)$, endowed with a family of de Rham tensors. 2.4.1 shows that $A^\prime_{O_L}$ and its principal polarization are the pull back of $(A_R^0,p_{A_R^0})$ through an $O_L$-valued point of $\Spec(R)$. If $A^0_{\dbF}$ has ordinary reduction, then [Va2, 2.3.17.1-2 and 3.1.0 b)] implies that its canonical lift is obtained from $A_R^0$ through a $W(\dbF)$-valued point of $\Spec(R)$. In this way we get (via the fact that all Hodge cycles of $A^\prime_{O_L}$ are ``preserved" under this deformation over $R$; see [Va1, 5.4.5 and 5.2.10]) the $\dbQ$--part of Theorem 3. As [Fa3, ch. 5] and [Fa2, 2.6] deal just with global contexts, for the $\dbQ_p$-part of Theorem 3 we use besides loc. cit. the global deformations of [Va2, 3.6 and 3.15.1].
\finishproclaim

\Proclaim{1.8.3. A complementary way to Pink's results.} \rm
The above ideas of 1.8 can be used (see 4.4.10) to reobtain [Pi, 5.10-11] (via [Va2] and Serre's and Katz--Messing's results used in [Pi, \S3]). We hope they will be useful in other situations as well: in the part of [Va7] related to the Mumford--Tate conjecture, we start from where the ideas of 4.4.7-9 and 3.3 stop; see 4.4.11 and 5.2 for samples.
\finishproclaim

\Proclaim{1.8.4. On the splitness assumption.} \rm
The assumption $G_{\dbQ_p}$ is split is not needed for Theorem 2 and so implicitly for Theorem 1. Moreover, for most of the above part of 1.8 it is enough to assume $G_{\dbQ_p}$ is unramified, cf. 4.4.12. However, as in this paper we do not perform any extrapolations (they are hinted at in 4.4.11) and as for our present main hope (see 5.2) to a complete proof of the Mumford--Tate conjecture we do need $G_{\dbQ_p}$ to be split, it has been felt practical to make the ideas more palatable by assuming most of the time $G_{\dbQ_p}$ is split.
\finishproclaim

\Proclaim{1.9. Some applications.} \rm
In \S5 we also present (besides Theorem 1, examples for a) of it and S4) briefly some applications (see the definition of Serre's volumes in 5.3, the complement to [No] in 5.1.6 and the remark 5.4 2) on the ordinary reduction conjecture) and some other conclusions (like 5.1.5, 5.4 8) and 5.5.1-2). 5.1.5 is an application to direct sums of abelian motives over number fields; we refer to it as the independence property (of suitable Galois representations). 5.5.1 extends slightly the main result (pertaining to the existence of some cocharacters) of [Wi] to certain $p$-divisible groups endowed with $\dbQ_p$-\'etale Tate tensors over a $V$ as in 1.5. 5.5.2 is a complement to [Va1, 4.3.10.1 1)].
\finishproclaim

\Proclaim{1.10. Acknowledgements.} \rm
We would like to thank UC at Berkeley (resp. Univ. of Utah and Univ. of Arizona) for providing us with good conditions for the writing of this manuscript (resp. for improving on it). We would like to thank Prof. J. S. Milne for providing us with the manuscript [Mi4] at the right time for being able to build up properly on [Mi4, 0.1]. We would like to thank Prof. S. G. Tankeev for making us aware of 4.2.5 and Prof. G. Faltings for pointing out an inexactitude in the first version of 2.2.2.1. 

Also, we are obliged to [Pi]: though this work was done independently of [Pi], as afterthoughts, we did not hesitate to simplify the things using [Pi] in 5.1; moreover, the definition 5.1.2.1 1) (based on [Pi, 5.10]) and the numerical analysis of 5.1.2.2 are slightly more general than the ones we initially thought of. In other words: the ideas of 3.1-2 and \S4 were clear to us in February 1996, but we learned Serre's theory of [Chi] later on (1998-9). This research was partially supported by the NSF grant DMS 97-05376.
\finishproclaim

\bigskip
\noindent
{\boldsectionfont \S2. Some complements}
\bigskip 

In this chapter we include some complements to [Va1-2].
Let $p$ and $k$ be as in 1.3. 

\smallskip
\Proclaim{2.1. Piano rings.} \rm
Let $e\in\dbN$. Let $T$ be an independent variable. We denote by $\Rtil e$ (resp. by $Re$) the subring of $B(k)[[T]]$ formed by formal power series $\sum_{n=0}^\infty a_nT^n$ for which $a_n{\fracwithdelims[]{n}{e}}!\in W(k)$ for all $n$ (resp. for which the sequence $b_n:=a_n{\fracwithdelims[]{n}{e}}!$ is formed by elements of $W(k)$ and converges to 0). Sometimes, when we want to emphasize $k$, we write $\Rtil e(k)$ (resp. $Re(k)$) instead of $\Rtil e$ (resp. of $Re$). Their Frobenius lifts (see [Fa3] or [Va1, 5.2]; for $p=2$ see also [Va2, 2.3.18 B8]) are defined by: $T$ goes to $T^p$. Obviously we have a natural $W(k)$-monomorphism $Re\hookrightarrow\Rtil e$. The rings $Re$ were used for the first time with significant impact in the original versions of [Fa3], while the rings $\Rtil e$ were used for the first time with significant impact in [Va1, 5.2] (resp. in [Va2, 2.3.18 B]) for $p\Ge 3$ (resp. for $p=2$). We refer to $\Rtil e$ (resp. to $Re$) as the piano ring of $k$ of resonance $e$ (resp. of convergent resonance $e$). This terminology is inspired from their construction and from the following important Fact (for $p\ge 3$, see [Fa3, ch. 2] and [Va1, 5.2.1]; the case $p=2$ is very much the same):

\medskip\noindent
{\bf Fact.} {\it We assume $p\ge 3$ (resp. $p=2$). For any DVR $V$ which is a finite, flat, totally ramified extension of $W(k)$ of degree $e$, there are $W(k)$-epimorphisms $Re\twoheadrightarrow V$ and $\Rtil e\twoheadrightarrow V$ (resp. there are $W(k)$-epimorphisms $Re\twoheadrightarrow V$); they are defined by: $T$ is mapped into an arbitrarily chosen uniformizer of $V$. Moreover, mapping $T$ to 0, we get $W(k)$-epimorphisms $Re\twoheadrightarrow W(k)$ and $\Rtil e\twoheadrightarrow W(k)$ respecting Frobenius lifts.}

\medskip
We view the $W(k)$-epimorphisms of the above Fact as the keys of a piano tuned on convergent or non-convergent (resp. convergent) resonances of length $e$.  
If $e^\prime\in\dbN$ is such that $e^\prime\Ge e$, then we have natural $W(k)$-monomorphisms (inclusions)
$$Re^\prime\hookrightarrow Re$$
and $\Rtil e^\prime\hookrightarrow\Rtil e$. So in the Fact, we can substitute ``of degree $e$" by: of degree at most $e$.
\finishproclaim

\Proclaim{2.2. Digression on a conjecture of Milne.} \rm
For the original conjecture of Milne we have in mind we refer to [Mi4, 0.1] and [Va1, 5.6.6]. Let $A$ be an abelian variety over a DVR $V$ of mixed characteristic $(0,p)$. Let $K:=V\fracwithdelims[]1p$. Let $(s_{\alpha})_{\alpha\in\scrJ}$ be a suitable family of Hodge cycles of $A_K$. In the context of abelian varieties, the conjecture asserts that under some mild conditions the below question has a positive answer. 

\medskip
{\bf Q.} {\it Is there a faithfully flat, $p$-adically complete $V$-algebra $R$ such that there is an $R$-isomorphism
$$\rho\colon H^1_{\text{\'et}}(A_{\Kbar},\dbZ_p)\otimes_{\dbZ_p} R\arrowsim H^1_{dR}(A/V)\otimes_V R$$
taking the $p$-component of the \'etale component of $s_{\alpha}$ into the de Rham component of $s_{\alpha}$, $\forall\alpha\in\scrJ$?}

\medskip
A similar question can be put for other classes of varieties or cycles (and so of objects). In particular, using Fontaine's comparison theory to define \'etale counterparts of de Rham tensors, we have an analogue of Q for Shimura $p$-divisible groups. For an intensive study of this question (in the broader context of Shimura $p$-divisible groups) see [Va7]; for an analogue study for K3 surfaces we refer to [Va4]. Here we ``touch" this question just as far as it brings more light to parts of \S3. 

Warning: Tate-twists of Hodge cycles are allowed; $\dbZ_p(m)$, with $m\in\dbZ$, is always identified, as a $\dbZ_p$-module, with $\dbZ_p$, and the same applies with $\dbZ_p$ being replaced by $\dbQ_p$, $\dbZ$ or $\dbQ$. For abelian varieties the situation is very much simplified: they are always polarized and so, as in [Va2, 2.3.1], in practice we just need to keep track of such polarizations and of Hodge cycles which involve no Tate-twists and of whose different realizations (like Betti, de Rham, etc.) are homogeneous tensors of corresponding essential tensor algebras. 

To study {\bf Q} we can assume $V$ is a complete DVR having an algebraically closed residue field. If this is so and if the subgroup $G$ of $GL:=GL(H^1_{dR}(A/V))$ defined as the Zariski closure of the subgroup of $GL_K$ fixing the de Rham component of $s_{\alpha}$, $\forall\alpha\in\scrJ$, is a reductive group over $V$, then (this is standard) the answer to the above question is yes iff such an isomorphism $\rho$ exists for $R=V$.
\finishproclaim

\Proclaim{2.2.0. Definition.} \rm
We say $A$ has the MC property, if the answer to the above question is true, provided we work with the family of all Hodge cycles of $A$ (here MC stands for Milne's conjecture).
\finishproclaim

\Proclaim{2.2.1. The generic situation over a number field.} \rm
Let $A$ be an abelian variety over a number field $E$. Let $(s_{\alpha})_{\alpha\in\scrJ}$ be the family of all Hodge cycles of $A$. We assume it is equal to the family of all Hodge cycles of $A_{\Ebar}$ (i.e. we assume that by enlarging $E$ we do not gain extra Hodge cycles, cf. [De3]). We choose an embedding $E\hookrightarrow\dbC$; in what follows $A_{\dbC}$ refers w.r.t. it. Let 
$$L:=H^1_B(A_{\dbC},\dbZ).$$ 
Let $W:=L\otimes_{\dbZ} \dbQ$. Let $s_{\alpha}^B$ be the element of the tensor algebra of $W\oplus W^*$ which is the Betti realization of $s_{\alpha}$, $\alpha\in\scrJ$. 

Let $O_E$ be the ring of integers of $E$ and let $N\in\dbN$ be such that $A$ extends to an abelian scheme $\scrA$ over $\Spec(O_E\fracwithdelims[]1N)$. Let $D:=H^1_{dR}(\scrA/O_E\fracwithdelims[]1N)$. We have:
\finishproclaim

\smallskip\noindent
{\bf Fact.}  {\it There is a finite field extension $E_1$ of $E$ such that we have an $E_1$-isomorphism
$$i_{E_1}:W\otimes_{\dbQ} E_1\arrowsim D\otimes_{O_E\fracwithdelims[]1N} E_1$$
whose extension to an $E_1$-isomorphism between tensor algebras takes $s_{\alpha}^B$ into the de Rham realization $s_{\alpha}^{dR}$ of $s_{\alpha}$, $\forall\alpha\in\scrJ$. Moreover, there is $N_1\in\dbN$ a multiple of $N$ such that $i_{E_1}(L\otimes_{\dbZ} O_{E_1}\fracwithdelims[]1{N_1})=D\otimes_{O_E\fracwithdelims[]1N} O_{E_1}\fracwithdelims[]1{N_1}$ (here $O_{E_1}$ is the ring of integers of $E_1$).} 

\medskip
This Fact is a (standard) consequence of the following fact: under the de Rham isomorphism
$$W\otimes_{\dbQ} \dbC\arrowsim D\otimes_E \dbC,$$
$s_{\alpha}^B$ is taken into $s_{\alpha}^{dR}$, $\forall\alpha\in\scrJ$ (cf. [De3, \S 1-2]). From this Fact and the standard connection between \'etale and Betti cohomologies (see [SGA4, Exp. XI]), we deduce:

\Proclaim{2.2.1.1. Corollary.} There is a multiple $M(A)\in\dbN$ of $N$ such that the pull back of $\scrA$ to any local ring of $\Spec(O_E\fracwithdelims[]1{M(A)})$ which is a DVR, has the MC property. 
\finishproclaim

In [Va7] we will see that $M(A)$ is effectively computable.

\Proclaim{2.2.2. The generic situation in the context of Shimura varieties of Hodge type.} \rm
Let $f\colon (G,X)\hookrightarrow (G\Sp(W,\psi),S)$ be an embedding of the Shimura pair $(G,X)$ of Hodge type into a Shimura pair defining a Siegel modular variety. Let $L$ be a $\dbZ$-lattice of $W$ such that we get a perfect alternating form $\psi\colon L\otimes_{\dbZ} L\to\dbZ$. Let 
$$K_4:=\{\text{g}\in G\Sp(L\otimes_{\dbZ}\dbZhat,\psi)(\dbZhat)|\text{g mod 4 is the identity element}\}.$$ 
Let $\scrM$ be the $\Spec(\dbZ\fracwithdelims[]12)$-scheme parameterizing isomorphism classes of principally polarized abelian schemes of relative dimension equal to $\dim_{\dbQ}(W)/2$ and having level-4 symplectic similitude structure; we have a natural identification
$$\scrM_{\dbQ}=\Sh_{K_4}(G\Sp(W,\psi),S),$$ 
cf. [Va1, 3.2.9 and 4.1.0].  Let $H$ be a compact, open subgroup of $G(\dbA_f)$ contained in $K_4$. We have a natural morphism $\Sh_H(G,X)\to\Sh_{K_4}(G\Sp(W,\psi),S)$, cf. [Va1, 2.4]. Let $\scrN$ be the normalization of $\scrM$ in the ring of fractions of $\Sh_H(G,X)$. $\scrM$ is quasi-projective over $\Spec(\dbZ\fracwithdelims[]12)$ (see [MFK, p. 139]) and so an excellent scheme. So $\scrN$ is of finite type over $\Spec(\dbZ\fracwithdelims[]12)$. Let $(\scrA,\scrP_{\scrA})$ be the pull back to $\scrN$ of the universal principally polarized abelian scheme over $\scrM$. For the sake of simplicity we assume that all Hodge cycles of $\scrA_{\Spec(\overline{\scrK})}$ are as well Hodge cycles of $\scrA_{\Spec(\scrK)}$, for any $\scrN$-scheme $\Spec(\scrK)$ which is the spectrum of the field of fractions of a connected component of $\scrN_{E(G,X)}$. In other words, we assume $H$ is small enough. 
We have:
\finishproclaim

\Proclaim{2.2.2.1. Lemma.}
There is a number $M(f,L)\in\dbN$ such that for every DVR $R$ of mixed characteristic $(0,p)$, with $p\Ge M(f,L)$, any abelian scheme $A_R$ obtained from $\scrA$ by pull back through a morphism $m_R:\Spec(R)\to\scrN$ and of whose Mumford--Tate group is $G$ itself, has the MC property.
\finishproclaim

\proof
As $\scrN$ is of finite type over $\Spec(\dbZ\fracwithdelims[]12)$, generically this is very much the same as 2.2.1. Instead of number fields and their rings of integers we have to work with fields of fractions of different connected components of $\scrN_{E(G,X)}$ and with rings of global functions of different affine, integral, open subschemes of $\scrN$. We deduce the existence of an open, dense subscheme
$U$ of $\scrN$ such that $m_{V_1}^*(\scrA)$ has the MC property, for any morphism $m_{V_1}:\Spec(V_1)\to\scrN$, with $V_1$ a DVR of mixed characteristic, factoring through $U$. Let $N_d\in\dbN$ be such that for any prime $q\Ge N_d$, $U_{\dbF_q}$ is dense in $\scrN_{\dbF_q}$. 

Let $N_r\in\dbN$ be such that the Zariski closure of $G$ in $GL(L[{1\over {(N_r-1)!}}])$ is a reductive group scheme. From [Va1, 5.8.6] and [Va2, 2.3.6] we deduce that for any prime $v$ of $E(G,X)$ dividing a prime $q\Ge\max\{5,N_r,{1\over 2}\dim_{\dbQ}(W)\}$, the triple $(f,L\otimes_{\dbZ} \dbZ_{(q)},v)$ is a SHS.   

Now we can take 
$$M(f,L):=\max\{5,N_d,N_r,{1\over 2}\dim_{\dbQ}(W)\}.$$ 
To see this, let $m_R$ and $A_R$ be as in the text of the Lemma. The $\dbZ_p$-lattice $H^1_{\text{\'et}}(A_{\overline{R[{1\over p}]}},\dbZ_p)$ together with the family of $p$-components of \'etale components of Hodge cycles of $A_{\overline{R[{1\over p}]}}$ does not depend on the choice of such a $\dbZ_{(p)}$-morphism $m_R$ (cf. the moduli interpretations of [Va1, 4.1]). From this, the definition of $M(f,L)$ and [Va2, 2.3.13.1] we get: to check that $A_R$ has MC, we can assume $m_R$ factors through $U$. Conclusion: $A_R$ has the MC property. This ends the proof.

\smallskip
In [Va7] we will see that $M(f,L)$ is effectively computable.

\Proclaim{2.2.3. The case of canonical lifts.} \rm
We now assume $k=\kbar$ and $V=W(k)$. We also assume $A$ is the canonical lift of an ordinary abelian variety over $k$. Let $(M,F^1,\phi)$ be the filtered $\sigma$-crystal over $k$ of $A$. The slopes of $(M,\phi)$ are just 0 and 1. Let $\mu\colon\dbG_m\to GL(M)$ be the canonical split cocharacter of $(M,F^1,\phi)$ (cf. [Wi] and the conventions of [Va2, 2.1]). We get a direct sum decomposition $M=F^1\oplus F^0$, with $\gamma\in\dbG_m(W(k))$ acting through $\mu$  on $F^i$ as the multiplication with $\gamma^{-i}$, $i=\overline{0,1}$. We have $\phi(F^i)\subset F^i$, $i=\overline{0,1}$. For $p\Ge 3$ (resp. for $p=2$), let $B^+(W(k))$ be the Fontaine's ring as defined in [Fa3] (see also [Va1, 5.2]) (resp. as used in [Va2, 2.3.18.1 E]). 

We consider the $B^+(W(k))$-monomorphism 
$$\rho\colon M\otimes_{W(k)} B^+(W(k))\hookrightarrow H^1_{\text{\'et}}(A_{\overline{B(k)}},\dbZ_p)\otimes_{\dbZ_p} B^+(W(k))$$
of the integral version of Fontaine's comparison theory (see [Fa3, th. 5 or 7]; for $p=2$ see also [Va2, 2.3.18.1 E]). For $p\ge 3$ (resp. for $p=2$), let $\beta\in B^+(W(k))$ be defined as in [Va1, 5.2.3] but working with $B^+(W(k))$ (resp. be defined as in [Va2, 2.3.18.1 E]). We have:
$$\rho\bigl(\mu(\beta)(M\otimes_{W(k)} B^+(W(k)))\bigr)=H^1_{\text{\'et}}(A_{\overline{B(k)}},\dbZ_p)\otimes_{\dbZ_p} B^+(W(k)).$$
This can be checked easily starting from $W(k)$-basis $\{f_i|i\in\{1,...,\dim(A)\}\}$ and $\{e_i|i\in\{1,...,\dim(A)\}\}$ of $F^1$ and respectively of $F^0$ such that for all $i\in\{1,...,\dim(A)\}$ we have $\phi(f_i)=pf_i$ and $\phi(e_i)=e_i$. For $p=2$ we need to add: we can assume $W(k)f_i$ (resp. $W(k)e_i$) is the direct summand of $M$ defined naturally by a $p$-divisible group which is a direct summand of the $p$-divisible group of $A$. 

Any Hodge cycle of $A$ is a de Rham cycle in the sense of [Bl] (cf. [Va2, th. of 2.3.9 A]). So (to be compared with [Va2, cor. of 2.3.10]) $\phi$ fixes the de Rham component $HC_{dR}$ of any Hodge cycle $HC$ of $A$. This implies: $\mu(\beta)$, viewed as a $B^+(W(k))\fracwithdelims[]1{\beta}$-linear automorphism of $M\otimes_{W(k)} B^+(W(k))\fracwithdelims[]1{\beta}$, fixes $HC_{dR}$. So $\rho\circ\mu(\beta)$ takes $HC_{dR}$ into the $p$-component of the \'etale component of $HC$. We get the following Fact (it solves [Va2, 4.6.5]):

\medskip\noindent
{\bf Fact.} {\it $A$ has the MC property}.
\finishproclaim

\Proclaim{2.3. Representation well positioned families of tensors.} \rm
Let $O$ be a DVR. Let $\pi$ be a uniformizer of it. Let $K:=O[{1\over {\pi}}]$. Let $M_O$ be an $O$-lattice of a finite dimensional $K$-vector space $W$. Let $G_K$ be a reductive subgroup of $GL(W)$. For the notion (resp. theory) of $O$-well positioned families of tensors w.r.t. $G_K$ we refer to [Va1, 4.3.4 and 4.3.7 1)] (resp. to [Va1, 4.3]). Here we introduce some extra definitions, in order to complement [Va1, 4.3.10.1 1)] later on (see 4.3.4.8 b) and 5.5.2).
\finishproclaim

\Proclaim{2.3.1. Definition.} \rm
A family $(v_{\alpha})_{\alpha\in\scrJ}$ of $G_K$-invariant tensors of $\scrT(M_O)$ is said to be $M_O$-representation very well positioned for $G_K$ if the following 3 things hold:

\medskip
\item{{\bf a)}}
it is $O$-very well positioned for $G_K$;

\smallskip
\item{{\bf b)}}
for every local, integral, faithfully flat $O$-algebra $R$, and for any free $R$-module $M$ satisfying $M\fracwithdelims[]1{\pi}=W\otimes_K R\fracwithdelims[]1{\pi}$,  enveloping it, and such that $\forall\alpha\in\scrJ$, $v_{\alpha}$ is a tensor of $\scrT(M_O)$ having a non-zero image  modulo $\pi$ iff it is a tensor of $\scrT(M)$ having a non-zero image modulo $m_R$ (with $m_R$ the maximal ideal of $R$), there is an isomorphism
$$\rho\colon M_O\otimes_O R^{\sh}\arrowsim M\otimes_R R^{\sh},$$
taking $v_{\alpha}$ into $v_{\alpha}$, $\forall\alpha\in\scrJ$;

\smallskip
\item{{\bf c)}}
$G_K$ is the subgroup of $GL(M\fracwithdelims[]1{\pi})$ fixing $v_{\alpha}$, $\forall\alpha\in\scrJ$ (so $\rho$ of b) is defined by an element of $G_K(W\otimes_K R^{\sh}\fracwithdelims[]1p)$).
\finishproclaim

\smallskip
As in [Va1, 4.3.4] we have variants: we speak about a weakly (resp. strongly) $M_O$-representation very well positioned family of tensors; for this we keep c), we request $R$ in b) to be normal (resp. to be just reduced) instead of integral, and in a) we request the family to be weakly (resp. strongly) $O$-very well positioned for $G_K$. 

\Proclaim{2.3.2. Definition.} \rm
A family $(v_{\alpha})_{\alpha\in\scrJ}$ of $G_K$-invariant tensors of $\scrT(M_O)$ is said to be $M_O$-representation well positioned for $G_K$ if 2.3.1 a) is true, while 2.3.1 b) is true only for DVR's which are faithfully flat $O$-algebras (warning: for the sake of generality and flexibility, 2.3.1 c) is ignored here). 
\finishproclaim
 
\Proclaim{2.4. An inducing property for p-divisible groups in the ramified case.}  \rm 
What follows is a complement to the relative form [Va2, 3.6.18.5.7] of the inducing property of [Va2, 3.6.18.5]. It is desirable to have a variant of this relative form over complete, regular $W(k)$-algebras which are not formally smooth. Below we deal with the simplest case of a totally ramified, finite, DVR extension $V$ of $W(k)$. Let $K:=V[{1\over p}]$ and $e:=[V:W(k)]$. 

Let $H_V$ be a $p$-divisible group over $V$. Let 
$$H^1_{\text{\'et}}(H_V,\dbZ_p):=T_p(H_V)^*$$ 
be the dual of the Tate module of $H_V$. Let $(w_{\alpha})_{\alpha\in\scrJ}$ be a family of homogeneous tensors of $\scrT(H_{\text{\'et}}^1(H_V,\dbZ_p))$ fixed by $\Gal(K)$. We assume the existence of a reductive subgroup $\tilde G$ of $GL(H_{\text{\'et}}^1(H_V,\dbZ_p))$ having the following two properties:

\medskip
{\bf i)} its generic fibre is the subgroup of $GL(H_{\text{\'et}}^1(H_V,\dbZ_p))_{\dbQ_p}$ fixing $w_{\alpha}$, $\forall\alpha\in\scrJ$;

\smallskip
{\bf ii)} there is a subset $\scrJ_0$ of $\scrJ$ such that the family of homogeneous tensors $(w_{\alpha})_{\alpha\in\scrJ_0}$ is $\dbZ_p$-well positioned for $\tilde G_{\dbQ_p}$. 

\medskip
By the partial degrees of $w_{\alpha}$ we mean the numbers $a_{\alpha}$, $b_{\alpha}\in\dbN\cup\{0\}$ such that 
$$w_{\alpha}\in (H^1_{\text{\'et}}(H_V,\dbZ_p))^{\otimes a_\alpha}\otimes_{\dbZ_p} (H^1_{\text{\'et}}(H_V,\dbZ_p)^*)^{\otimes b_\alpha}.$$
We also assume that the partial degrees of $w_{\alpha}$ are smaller or equal to $\max\{1,p-2\}$, $\forall\alpha\in\scrJ_0$. We express this by: the partial degrees of $(w_{\alpha})_{\alpha\in\scrJ_0}$ are not bigger than $\max\{1,p-2\}$.
So in particular, we have:
$$\deg(w_{\alpha}):=a_{\alpha}+b_{\alpha}\Le \max\{2,2(p-2)\}.$$
Let $\Rtil e$ and $Re$ be as in 2.1. Let 
$$(M,F^1(M),\Phi_{M},\nabla)$$ 
be the filtered $F$-crystal over $Re/pRe$ obtained by taking the dual of the Lie algebra of the universal vector extension of $H_V$ (cf. [Me], [BBM] and [Fa3, ch. 6]). The main property (with $p\ge 2$) of $Re$ needed to construct it is:

\medskip\noindent
{\bf Fact.} {\it The $p$-th power of any non-invertible element of $Re/pRe$ is $0$. So $\forall n\in\dbN$, $Re/p^nRe$ is the inductive limit of its local, artinian $W_n(k)$-subalgebras taken by the reduction mod $2^n$ of the Frobenius lift of $Re$ into themselves (and ordered under the relation of inclusion).}

\medskip
Let $(t_{\alpha})_{\alpha\in\scrJ}$ be the family of tensors of $\scrT(M\fracwithdelims[]1p)$ corresponding to $(w_{\alpha})_{\alpha\in\scrJ}$ under Fontaine's comparison theory (see [Fa3, ch. 6]). 

We recall (see 2.1) that we have natural monomorphisms 
$$R:=W(k)[[T]]\subset Re\subset\Rtil e\subset B(k)[[T]];$$ 
the Frobenius lifts of $R$, $Re$ and $\Rtil e$ take $T$ into $T^p$. Let $FL$ be the Frobenius lift of $\Rtil e$. We fix a uniformizer $\pi_V$ of $V$. 

From now on we assume $p\ge 3$. $\pi_V$ defines naturally $W(k)$-epimorphisms $R\twoheadrightarrow V$, $Re\twoheadrightarrow V$ and $\Rtil e\twoheadrightarrow V$ (cf. Fact of 2.1); moreover, we have $W(k)$-epimorphisms $R\twoheadrightarrow W(k)$, $Re\twoheadrightarrow W(k)$ and $\Rtil e\twoheadrightarrow W(k)$ defined by: $T$ goes to 0. We refer to these six $W(k)$-epimorphisms as the natural ones (defined by $\pi_V$); at the level of spectra we speak about the six natural embeddings (defined by $\pi_V$).

We assume the existence of a cocharacter 
$$\mu_K\colon\dbG_m\to GL(M\otimes_{Re} K)$$ 
fixing the image $t_{\alpha}^{V}$ of $t_{\alpha}$ in $\scrT(M\otimes_{Re} K)$, $\forall\alpha\in\scrJ$, and such that it produces a direct sum decomposition 
$$M\otimes_{Re} K=F^1(M)\otimes_{Re} {K}\oplus F^0_{K}$$ 
with the property that $\beta\in\dbG_m(K)$ acts through $\mu_{K}$ trivially on the summand $F^0_K$ and as the multiplication with $\beta^{-1}$ on the summand $F^1(M)\otimes_{Re} K$. 

Let $G_{\Rtil e}$ (resp. $H_{\Rtil e}$) be the reductive group (resp. be a $p$-divisible group) over $\Rtil e$, obtained as in [Va1, 5.2 and 5.3.1] (resp. as $\tilde A$ of [Va1, 5.3.3]) starting from the pair $(H_V,(w_{\alpha})_{\alpha\in\scrJ})$. So $G_{\Rtil e}$ is the extension to $\Spec(\Rtil e$) of a reductive subgroup $G_{Re}$ of $GL(M)$. [Va1, 5.2-3] is worked out in the context of abelian schemes; but the same construction (finalized in [Va1, 5.3.3]) applies to the context of $p$-divisible groups. We can apply [Va1, 5.2-3] as the
 family of tensors $(w_{\alpha})_{\alpha\in\scrJ_0}$ is $\dbZ_p$-well positioned for $\tilde G_{\dbQ_p}$, is enveloped by $H^1_{\text{\'et}}(H_V,\dbZ_p)$, and is of partial degrees not bigger than $p-2$.
  The only difference from [Va1, 5.2-3]: as we did not assume $k$ to be algebraically closed, it might happen that the groups $\tilde G_{\Rtil e}$ and $G_{\Rtil e}$ are not isomorphic. 

Let $(M_0,F^1_0,\phi,G_{W(k)},(t^0_{\alpha})_{\alpha\in\scrJ})$ be the Shimura filtered $\sigma$-crystal over $k$ associated to the pull back $(H_{W(k)},(t_{\alpha}^0)_{\alpha\in\scrJ})$ of the pair $(H_{\Rtil e},(t_{\alpha})_{\alpha\in\scrJ})$ through the natural embedding $\Spec(W(k))\hookrightarrow\Spec(\Rtil e)$ (to be compared with [Va1, (5.3.5-12)]). The pair $(H_{W(k)},(t_{\alpha}^0)_{\alpha\in\scrJ})$ is a Shimura $p$-divisible group over $W(k)$.  
Let $G^0_{W(k)}$ be a smooth subgroup of $G_{W(k)}$ such that the natural morphism $m^0:G^0_{W(k)}\to G_{W(k)}/P$, with $P$ as the parabolic subgroup of $G_{W(k)}$ normalizing $F_0^1$, is smooth. We consider a Shimura filtered $F$-crystal with an emphasized family of tensors (see [Va2, 2.2.10])
$$\grC=(M_0,F^1_0,\phi,G_{W(k)},G^0_{W(k)},\tilde f,(t_{\alpha}^0)_{\alpha\in\scrJ}).$$ 
Let $(H,(t_{\alpha}^0)_{\alpha\in\scrJ})$ be the Shimura $p$-divisible group over the completion $\Ghat^0$ of $G^0_{W(k)}$ in its origin whose associated $p$-divisible object with tensors of $\scrM\scrF_{[0,1]}^\nabla(\Ghat^0)$ is $\grC$, cf. [Va2, 2.2.21 UP and 3.4.18.5.1]. We have:
\finishproclaim

\Proclaim{Proposition.}
$H_V$ together with its de Rham tensors $(t_{\alpha}^{V})_{\alpha\in\scrJ}$ is obtained from $H$ and its crystalline tensors $(t_{\alpha}^0)_{\alpha\in\scrJ}$, by pull back through a $V$-valued point of $\Ghat^0$. 
\finishproclaim

\proof 
Based on [Va2, 2.2.21 UP] we can assume $m^0$ is \'etale. As in [Va1, 5.3.1] we get: we can assume $\mu_{K}$ extends to a cocharacter $\mu_{V}\colon\dbG_m\to GL(M\otimes_{Re} V)$ factoring through $G_{V}$. The construction of $H_{\Rtil e}$ corresponds (cf. [Va1, 5.3.3]) to a lift of $\mu_{V}$ to a cocharacter $\mu_{\Rtil e}$ of $G_{\Rtil e}$. From [Va2, 2.3.17.2 and 2.2.21 UP] we get: it is irrelevant which lift we choose (i.e. it is irrelevant which $F_0^1$ we choose to construct $\grC$). 

For $n\in\dbN\cup\{0\}$, let $I_n$ (resp. $J_n$) be the ideal of $\Rtil e$ (resp. of $R$) formed by formal power series as in 2.1 but with $a_0=a_1=...=a_{n-1}=0$ (resp. generated by $T^n$). For $n\Le ep-1$ we have $\Rtil e/I_n=R/J_n$; but this does not hold for $n\Ge ep$. So to just ``copy" [Fa3, proof of th. 10 and rm. iii) after it] in our present context of (filtered $F$-crystals with tensors associated to) $(H,(t_{\alpha}^0)_{\alpha\in\scrJ})$ and $(H_{\Rtil e},(t_{\alpha})_{\alpha\in\scrJ})$, we need to point out that the triple 
$$(\Rtil e,(I_n)_{n\in\dbN\cup\{0\}},FL)$$ 
has all three properties required in order to be allowed to ``appeal" to loc. cit. They are:

\medskip
{\bf a)} $I_{n}/I_{n+1}$ is a free $W(k)$-module of finite rank (it is one), $\forall n\in\dbN\cup\{0\}$;

\smallskip
{\bf b)} $FL(I_{n})\subset I_{n+1}$, $\forall n\in\dbN\cup\{0\}$;

\smallskip
{\bf c)} $\Rtil e$ is complete w.r.t. the topology defined by $(I_n)_{n\in\dbN}$.

\medskip
So we can apply loc. cit. entirely to get that the pair $(H_{\Rtil e},(t_{\alpha})_{\alpha\in\scrJ})$ is obtained from the pair $(H,(t_{\alpha}^0)_{\alpha\in\scrJ})$ by pull back through a morphism $m_{\Rtil e}:\Spec(\Rtil e)\to\Ghat^0$ (to be compared with [Va2, 2.3.18 B9]; loc. cit. is the first place to mention that [Fa3, proof of th. 10 and rm. iii) after it] apply in larger contexts not necessarily involving just regular $W(k)$-algebras of formal power series). We just need to add: as $V/pV$ is a complete intersection, the fully faithfulness result of [BM, 4.3.2] and Grothendieck--Messing deformation theory imply that any $p$-divisible group over $\Rtil e$ is uniquely determined by its associated filtered $F$-crystal over $\Rtil e/p\Rtil e$. The Proposition follows from the existence of $m_{\Rtil e}$.

\Proclaim{2.4.1. Variant.} \rm
A result similar to Proposition of 2.4 can be stated if we add a principal quasi-polarization $p_{H_V}$ of $H_V$ with the property that the perfect alternating form on $H^1_{\text{\'et}}(H_V,\dbZ_p)$ associated to it is normalized by $\tilde G$. The same proof applies (as [Fa3, proof of th. 10 and rm. iii) after it] ``allows" Tate-twists as needed for this form). 
\finishproclaim

\Proclaim{2.4.2. Remark.} \rm
We consider a quintuple as above
$$\scrC_V:=(M,\Phi_M,\nabla,G_{Re},(t_{\alpha})_{\alpha\in\scrJ})$$ 
over $\Spec(Re/pRe)$, such that $G_{Re}$ is the reductive subgroup of $GL(M)$ whose fibre over $\Spec(Re[{1\over p}])$ is the subgroup of $GL(M[{1\over p}])$ fixing $t_{\alpha}$, $\forall\alpha\in\scrJ$, and $\mu_K$ exists. We claim it is obtained from a similar quintuple over $R/pR$ by extension of scalars, without assuming that 2.4 ii) holds or that the partial degrees of $(w_{\alpha})_{\alpha\in\scrJ_0}$ are not bigger than $\max\{1,p-2\}$. This is so due to the following two reasons.

\medskip
{\bf i)} Let $\mu_{\Rtil e}:\dbG_m\to G_{\Rtil e}$ be as in the proof of 2.4; let $F^1(M\otimes_{Re} \Rtil e)$ be the maximal direct summand of $M\otimes_{Re} \Rtil e$ on which $\dbG_m$ acts via $\mu_{\Rtil e}$ as the inverse of the identical cocharacter of $\dbG_m$. The resulting sextuple $(M\otimes_{Re} \Rtil e,F^1(M\otimes_{Re} \Rtil e),\Phi_M\otimes 1,G_{\Rtil e},(t_{\alpha})_{\alpha\in\scrJ})$ over $\Rtil e/p\Rtil e$ is obtained by pulling back a similar sextuple over $\Ghat^0_k$ through a $W(k)$-morphism $m_{\Rtil e}:\Spec(\Rtil e)\to\Ghat^0$ (cf. the paragraph before 2.4.1). 

\smallskip
{\bf ii)} The composite of the natural embedding $\Spec(V)\hookrightarrow\Spec{(\Rtil e)}$ with $m_{\Rtil e}$, factors through $\Spec(R)$. If $m_R$ is the $W(k)$-morphism $\Spec(R)\to\Ghat^0$ involved in this factorization, then (as we are dealing with $F$-crystals and crystalline tensors) $\scrC_V$ is uniquely determined by $m_R^*(\grC)$ (viewed without filtration) and so it is obtained by extension of scalars from a similarly defined quintuple over $R/pR$. 

\medskip
The same remains true for quintuples over $\Rtil e/p\Rtil e$ or when we are (as in 2.4.1) in a principally quasi-polarized context. Moreover, i) holds as well for $p=2$.
\finishproclaim

\Proclaim{2.5. Fields accompanying reflex fields.} \rm
We now introduce some notations and notions pertaining to simple, adjoint Shimura varieties, in order to facilitate future references and to have (later on) a better understanding of 1.5 d'). Let $n\in\dbN$.

Let $(G,X)$ be a simple, adjoint Shimura variety. It is well known (see [De1, 3.8]) that there is a (unique) totally real number field $E_1(G,X)$ such that $E(G,X)$ is either it or a totally imaginary quadratic extension of it, depending on the fact that $(G,X)$ is without or with involution. If $(G,X)$ is not of some $A_{2n-1}$ type, then (cf. [De1, 3.8 (ii)] and [Bou1, planche I to VI]) it is with involution iff it is of one of the following types $A_{2n}$, $D_{2n+3}^{\dbH}$, $D_{2n+3}^{\text{mixed}}$ and $E_6$. If $(G,X)$ is of $A_{2n-1}$ type then it can be with or without involution. For $K$ an arbitrary field of characteristic 0, we denote by $(G_K,[\mu_K])$ the extension to $K$ of the Shimura group pair $(G,[\mu])$ associated naturally (see [Va2, 2.2.6 5)]) to $(G,X)$. If $K$ is a $p$-adic field, we denote $[\mu_K]$ by $[\mu_p]$.

Let $F(G,X)$ be the totally real number field such that $G$ is the Weil restriction from $F(G,X)$ to $\dbQ$ of an absolutely simple, adjoint $F(G,X)$-group $G(F(G,X))$ (cf. [De2, 2.3.4]). Let $G^{\text{split}}$ be the split, simple, adjoint $\dbQ$--group of the same Lie type as $G$. It is well known that $Out(G):=Aut(G^{\text{split}})/G^{\text{split}}$ is either trivial, or $\dbZ/2\dbZ$, or the symmetric group $S_3$. So let $I(G,X)$ be the smallest field extension of $F(G,X)$ such that $G(F(G,X))_{I(G,X)}$ is an inner form of $G^{\text{split}}_{I(G,X)}$; its existence is a standard piece of Galois cohomology (to be compared with 3.1 F) below). We have $[I(G,X):F(G,X)]\in\{1,2,3,6\}$. 

\medskip\noindent
{\bf Definition.} $(G,X)$ is said to be of non-inner (resp. of inner) type if $I(G,X)\neq F(G,X)$ (resp. if $I(G,X)=F(G,X)$).

\medskip
Let $EF_1(G,X)$ (resp. $EF(G,X)$) be the composite field of $F(G,X)$ with $E_1(G,X)$ (resp. with $E(G,X)$). Let $G_1(G,X)$ (resp. $G(G,X)$) be the Galois extension of $\dbQ$ generated by $EF_1(G,X)$ (resp. by $EF(G,X)$). Let $G_2(G,X)$ be the Galois extension of $\dbQ$ generated by $F(G,X)$. $G_1(G,X)$ and $G_2(G,X)$ are totally real number fields.

We now assume $(G,X)$ is with involution. So $G(G,X)$ is a Galois extension of $\dbQ$ which is a totally imaginary quadratic extensions of $G_1(G,X)$; this follows from the fact (see [De1, p. 139]) that $E(G,X)$ is a CM field. Using standard facts on number fields (for instance, see [La, p. 168]), we get:

\medskip\noindent
{\bf Fact.} {\it The set of primes $p$ which split in $G_1(G,X)$ but do not split in $G(G,X)$ and are such that $G$ is unramified over $\dbQ_p$, is of positive Dirichlet density.} 

\medskip
If $p$ is such a prime, then $G_{\dbQ_p}$ is a product of absolutely simple, unramified, adjoint groups and any prime of $E_1(G,X)$ (resp. any prime $v$ of $E(G,X)$) dividing $p$ has residue field $\dbF_p$ (resp. $\dbF_{p^2}$). So at least one simple factor $(G_1,[\mu_1])$ of the Shimura group pair $(G_{\dbQ_p},[\mu_p])$ is such that the following two conditions hold with $s=2$:

\medskip
{\bf a)} $G_1$ is absolutely simple, non-split, unramified over $\dbQ_p$ and splits over $B(\dbF_{p^s})$;

\smallskip
{\bf b)} $\mu_1$ is a cocharacter of $G_{1B(\dbF_{p^s})}$ and the $G_1(B(\dbF_{p^2}))$-conjugacy classes $[\mu_1]$ and $[\gamma(\mu_1)]$ are distinct; here $\gamma\in\Gal(B(\dbF_{p^s})/\dbQ_p)$ is an arbitrary non-identity element and it acts naturally on the cocharacters of $G_{1B(\dbF_{p^s})}$.

\medskip
We also get that $(G,X)$ is of non-inner type.
\finishproclaim

\Proclaim{2.5.1. Special involutions.} \rm
We now assume $(G,X)$ is of $D_{2n+2}^{\dbH}$ or $D_{2n+2}^{\text{mixed}}$ type. So $E(G,X)=E_1(G,X)$ and $I(G,X)$ is a totally real number field. Let $l$ be a rational prime.

\medskip\noindent
{\bf Definition.} We say $(G,X)$ is: 

\smallskip
{\bf a)} without (resp. with) apparent $\dbR$-involution if $EF(G,X)=F(G,X)$ (resp. if $EF(G,X)\neq F(G,X)$);

\smallskip
{\bf b)} with $\dbR$-involution if $E(G,X)$ is not a subfield of $G_2(G,X)$;

\smallskip
{\bf c)} with $\dbQ_l$-involution if there is a finite field extension $F^1_l$ of $\dbQ_l$ with the property that the Shimura group pair $(G_{F^1_l},[\mu_l])$ has a simple factor $(G_1,[\mu_1])$ such that $G_1$ is an absolutely simple, non-split, unramified $F^1_l$-group and $\mu_1$ is a cocharacter of $G_{1F^2_l}$ whose $G_1(F^2_l)$-conjugacy class is not fixed by $\Gal(F^2_l/F^1_l)$; here $F^2_l$ is an unramified extension of $F^1_l$ of degree 2 or (in case $(G,X)$ is of $D_4^{\text{mixed}}$ type) 3.

\medskip
Warning: here $\dbR$ refers to the fact that we are dealing just with totally real number fields while the word involution points out that suitable Shimura group pairs (not necessarily over such fields; often we have to be in a $\dbQ_l$-context similar to the one of c)) are (in the logical sense) with involutions. Warning: the notions of apparent $\dbR$-involution and of non-inner type are unrelated (i.e. no one implies the other one).

As above, if $(G,X)$ is with $\dbR$-involution, then the set or primes $p$ which split in $F(G,X)$ but do not split in $EF(G,X)$ and are such that $G$ is unramified over $\dbQ_p$, is of positive Dirichlet density. Moreover, for such a prime $p$, if $(G,X)$ is not (resp. is) of $D_4^{\text{mixed}}$ type we similarly get the existence of a simple factor $(G_1,[\mu_1])$ of $(G_{\dbQ_p},[\mu_p])$ for which 2.5 a) and b) hold for $s=2$ (resp. for some $s\in\{2,3\}$). So $(G,X)$ is with $\dbQ_p$-involution. It is also with apparent $\dbR$-involution as well as of non-inner type.
\finishproclaim

\Proclaim{2.6. Hodge quasi products.} \rm
We consider a finite family of injective maps $f_i:(G_i,X_i)\hookrightarrow (GSp(W_i,\psi_i),S_i)$, $i\in I:=\{1,...,n\}$, into Shimura pairs defining Siegel modular varieties. For what follows we refer to [Va1, Example 5 of 2.5] and to [Va2, 4.9.2] (they deal with the case $n=2$). Let 
$$(W,\psi):=\oplus_{i=1}^n (W_i,\psi_i).$$ 
Let $G_0:=\prod_{i=1}^n G_i$. It is a subgroup of $\prod_{i=1}^n GSp(W_i,\psi_i)$ and so of $GL(W)$. Let $G$ be the connected component of the origin of the intersection $G_0\cap GSp(W,\psi)$. The monomorphism $G\hookrightarrow G_0$ gives birth to an isomorphism $G^{\der}\tilde\to G_0^{\der}$. Any monomorphism $h_x:\Res_{\dbC/\dbR} \dbG_m\hookrightarrow G_{0\dbR}$ defining an element $x\in\prod_{i=1}^n X_i$ factors through $G_{\dbR}$. Denoting by $X$ (resp. by $S$) the $G(\dbR)$-conjugacy (resp. $GSp(W,\psi)(\dbR)$-conjugacy) class of $h_x$, we get an injective map
$$f:(G,X)\hookrightarrow (GSp(W,\psi),S),$$
which we call a Hodge quasi product of $f_i$'s; warning: it does depend on the choice of $x$ (cf. [Va1, 2.5.1]). We denote any such $f$ by $\times^{\scrH}_{i\in I} f_i$. $(G,X)$ itself is referred as a Hodge quasi product of $(G_i,X_i)$'s and is denote by ${\times_{\scrH}}_{i\in I} (G_i,X_i)$. If $n=2$, then we also denote $f$ by $f_1\times^{\scrH} f_2$ and $(G,X)$ by $(G_1,X_1)\times_{\scrH} (G_2,X_2)$.

Let $L_i$ be a $\dbZ$-lattice of $W_i$ such that we get a perfect form $\psi_i:L_i\otimes_{\dbZ} L_i\to\dbZ$. We say the $\dbZ$-lattice $L:=\oplus_{i=1}^n L_i$ of $W$ is well adapted for Hodge quasi products of $f_i$'s. 

\medskip\noindent
{\bf Fact. 1)} {\it $Z^0(G)=\dbG_m$ iff $Z^0(G_i)=\dbG_m$, $\forall i\in I$.}

\smallskip
{\bf 2)} {\it If $f_i$ are PEL type embeddings, then $f$ is as well a PEL type embedding.}

\medskip
\proof 
1) follows from the existence of a short exact sequence
$$0\to Z^0(G)\hookrightarrow Z^0(G_0)\twoheadrightarrow \dbG_m^{n-1}\to 0.$$
\indent
For 2) we just need to mention: if $\scrB_i$ is the $\dbQ$--subalgebra of $\End(W_i)$ formed by elements fixed by $G_i$, then $G$ is the connected component of the origin of the subgroup of $GSp(W,\psi)$ fixing each element of the $\dbQ$--subalgebra $\prod_{i=1}^n \scrB_i$ of $\End(W)$.

\bigskip
\noindent
{\boldsectionfont \S3. Good embeddings into unitary Shimura varieties}
\bigskip

This chapter has two distinctive parts. First, in 3.1-2 (resp. in 3.2.1-2) we prove a rational version (resp. an effective generic version) of the result stated in [Va2, 1.15.8]. Second, in 3.3 we apply 3.2 to the understanding of Frobenius tori of principally polarized abelian varieties over number fields. As 3.1 and 3.3 contain many important features, its subsections are carefully itemized using letters or numbers. As part of the proof of 3.2.1 is slightly sketched, 3.2.1 is starred (cf. also 3.2.2 4)).

\Proclaim{3.1. A setting.} \rm
{\bf A)} What follows is very close in spirit to [Va1, 6.5-6] and Erratum. We start with a simple Shimura variety $\Sh(G,X)$ of adjoint, abelian type and of positive dimension. So $G$ is a non-trivial, simple, adjoint $\dbQ$--group. Let $F:=F(G,X)$; so (see 2.5) we have
$$G=\Res_{F/\dbQ} G(F).$$ 
For any totally real number field $F_1$ containing $F$, let
$$G^{F_1}:=\Res_{F_1/\dbQ} G(F)_{F_1}.$$ 
As in [Va1, 6.6.1] we get naturally a Shimura pair $(G^{F_1},X^{F_1})$ and an injective map 
$$(G,X)\hookrightarrow (G^{F_1},X^{F_1}),$$
with $X^{F_1}$ as in loc. cit. Let $p\Ge 3$ be a prime such that $G$ is unramified over $\dbQ_p$; this forces $F$ to be unramified above $p$ (this is the same as [Va1, 6.5.1]). In what follows we always choose $F_1$ to be unramified over $p$; this implies $G^{F_1}$ is as well unramified over $\dbQ_p$.

\smallskip
{\bf B)} In practice, it is convenient (though not really needed) to choose whenever possible $F_1$ to be large enough so that $G(F)_{F_1}$ splits over the completion of $F_1$ w.r.t. any finite prime of it. If $G$ is of $B_n$ or $C_n$ Lie type, then such a field $F_1$ always exists due to the fact that $G(F)$ itself splits over all but finite completions of $F$ w.r.t. finite primes of it; argument: an absolutely simple, adjoint group of $B_n$ or $C_n$ Lie type over the ring of integers of a local field, is a split group.

However, for the sake of generality, below even when possible, we do not assume $F_1$ is such that $G(F)_{F_1}$ splits over the completion of $F_1$ w.r.t. any finite prime of it. This generality is hinted at in E) below but in fact is not really used here.

\smallskip
{\bf C)} Till end of Proposition 3.2 we assume $G$ is not of $A_n$ Lie type, $n\in\dbN$. We now introduce the ``simplest" Shimura pair $(\tilde G,\tilde X)$ with the property that $\tilde G_{\dbR}$ is isomorphic to a simple, non-compact factor of $G_{\dbR}$. So $\tilde G$ is an absolutely simple $\dbQ$--group of adjoint type and of the ``simplest" possible nature. The choice of such $\tilde G$ is supported by the fact (see [De2]) that Shimura varieties of $D_n^{\text{mixed}}$ type are not of preabelian type. Depending on the type of $(G,X)$ we take $\tilde G$ as follows.

\medskip
$\bullet$ If $(G,X)$ is of $B_n$ type ($n\in\dbN\setminus\{1\})$, then $\tilde G:=SO(2,2n-1)$ is the group of $2n+1$ by $2n+1$ matrices of determinant 1 leaving invariant the quadratic form $-x_1^2-x_2^2+x_3^2+...+x_{2n+1}^2$ on $\dbQ^{2n+1}$.

\smallskip
$\bullet$ If $(G,X)$ is of $C_n$ type ($n\in\dbN\setminus\{1,2\})$, then we take $\tilde G$ to be the split, simple $\dbQ$--group of $C_n$ Lie type. 

\smallskip
$\bullet$ If $(G,X)$ is of $D_n^{\dbR}$ type ($n\in\dbN$, $n\Ge 4$), then $\tilde G:=SO(2,2n-2)^{\ad}$ ($SO(2,2n-2)$ is defined in the same manner as $SO(2,2n-1)$).

\smallskip
$\bullet$ If $(G,X)$ is of $D_n^{\dbH}$ type ($n\in\dbN$, $n\Ge 4$), then $\tilde G:=SO^*(2n)^{\ad}$ is the adjoint group of the $\dbQ$--version $SO^*(2n)$ of the $\dbR$-version described in [He, p. 445] (and denoted in the same way); so $SO^*(2n)(\dbQ)$ is the subgroup of $SO(2n)(\dbQ(i))$ formed by matrices leaving invariant the skew hermitian form $-z_1\overline{z}_{n+1}+z_{n+1}\overline{z}_1-....-z_n\overline{z}_{2n}+z_{2n}\overline{z}_n$ over $\dbQ(i)$ (here $\overline{z}_i$ means the complex conjugate of the complex number $z_i$, while $SO(m)$, $m\in\dbN$, is the semisimple $\dbQ$--group leaving invariant the quadratic form $x_1^2+...+x_m^2$ on $\dbQ^m$; in our case $m=2n$ and the connection between $z_i$'s and $x_i$'s is achieved --over $\dbQ(i)$-- via $z_i=x_i$).

\medskip
The values of $n\in\dbN$ are chosen so that we do not get repetitions among the types (so for the $D_4^{\dbH}$ and $D_4^{\dbR}$ types, the convention of [De2, end of 2.3.8] applies). If $\tilde G$ is not $SO^*(2n)^{\ad}$ for some $n\in\dbN$, $n\Ge 4$, then there is only one $\tilde G(\dbR)$-conjugacy class $\tilde X$ of homomorphisms $SO(2)_{\dbR}\to\tilde G_{\dbR}$ given birth to a Shimura pair $(\tilde G,\tilde X)$; in the excluded cases we have two choices for such a conjugacy class, cf. [De2, 1.2.7-8]. In what follows we are interested just in $\tilde G$: the role of $\tilde X$ is irrelevant; $\tilde X$ is mentioned just to emphasize that $\tilde G$ can be used as the first entry of a simple, adjoint Shimura pair.

\smallskip
{\bf D)} Similarly, in what follows $SU(m,m)$ is the $\dbQ$--version of the $\dbR$--version described in [He, p. 444] (and denoted in the same way). Till end of 3.2.1 we have the same restrictions on $n$ as in C). We introduce a simply connected semisimple $\dbQ$--group $G(A)$ by the rules:

\medskip
-- if $(G,X)$ is of $B_n$ or $D_n^{\dbR}$ type, then $G(A):=SU(2^{n-1},2^{n-1})$;

\smallskip
-- if $(G,X)$ is of $C_n$ or $D_n^{\dbH}$ type, then $G(A):=SU(n,n)$.

\medskip
The above choices for $G(A)$ conform with [Sa, 3.3-7 or p. 461]. So we naturally have a homomorphism 
$$h_{\dbR}\colon\tilde G^{\sc}_{\dbR}\to G(A)_{\dbR}$$
of finite kernel which at the level of Lie algebras satisfies the condition $(H_2)$ of [Sa, p. 427]. Here $\tilde G^{\sc}$ is the simply connected semisimple group cover of $\tilde G$. $\text{Ker}(h_{\dbR})$ is trivial if $(G,X)$ is of $B_n$, $C_n$ or $D_n^{\dbR}$ type. 

The only difference from [Sa, 3.3-7 or p. 461] is: if $(G,X)$ is of $D_n^{\dbR}$ type, we consider the spin representation and not just a half spin representation; so we have $G(A)=SU(2^{n-1},2^{n-1})$ and not $SU(2^{n-2},2^{n-2})$, i.e. we implicitly work with a product of two copies of $SU(2^{n-2},2^{n-2})$ embedded naturally in $SU(2^{n-1},2^{n-1})$. Of course, if $(G,X)$ is of $D_n^{\dbR}$ type we can often work as well with just one half spin representation (i.e. we can often work just with $SU(2^{n-2},2^{n-2})$). But then, for $n$ even, $\text{Ker}(h_{\dbR})$ is non-trivial and we can not get 3.2 f) below for the $D_n^{\dbR}$ type case. 

We can assume $h_{\dbR}$ is the extension to $\dbR$ of a homomorphism
$$h\colon\tilde G^{\sc}\to G(A)$$
(so the notation $h_{\dbR}$ is justified). This is a consequence of constructions in [Sa, 3.3-7]: instead of $\dbR$, $\dbC$ and $\dbH$, in loc. cit. we can work as well with $\dbQ$, $\dbQ(i)$ and respectively the standard quaternion algebra $\dbQ(i,j)$ over $\dbQ$. Let $\tilde G^d$ be the image of $h$ in $G(A)$.

\smallskip
{\bf E)} Of course, instead of $\dbQ(i)$, in C) and D) we can work as well with a totally imaginary quadratic extension $K$ of $\dbQ$ unramified over $p$, and instead of $\dbQ(i,j)$ of D) we can work as well with a quaternion $\dbQ$--algebra containing $K$ and ramified at the infinite place. What we mean by this: we have variants, where instead of $\tilde G$ or $G(A)$, we work with some forms of them, which split over $K$. However, as $\dbQ(i)$ is unramified above any odd prime, C) and D) as they stand are all we need in what follows; so we do not detail below the case of an arbitrary such $K$. Just one example: in some cases we can work with the SO-group of the quadratic form $-x_1^2-x_2^2+lx_3^2...+lx_5^2$ on $\dbQ^5$ instead of $SO(2,3)$, and with $K=\dbQ(\sqrt{-l})$ instead of $\dbQ(i)$, etc. (here $l\in\dbN$ is such that $K$ is unramified above $p$).

\smallskip
{\bf F)} After stating a Fact, we divide the rest of F) in three subparts.

\medskip\noindent
{\bf Fact.} {\it For $F_1$ large enough (subject to requirements of A)), $G(F)_{F_1}$ is an inner form of $\tilde G_{F_1}$.}

\medskip
{\bf i)} We recall the arguments for this well known Fact. If $G$ is of $B_n$ or $C_n$ Lie type, then we can take $F_1=F$. Argument: $G(F)$ has no outer automorphisms. 

If $G$ is of $D_n$ Lie type, $n\Ge 4$, then the Fact is a consequence of the following remarks: 

\medskip
-- $G(F)$ is a form (not a priori inner) of $\tilde G_F$;

\smallskip
-- the group scheme of automorphisms modulo inner automorphisms of $\tilde G_F$ is $\dbZ/2\dbZ$ if $n\Ge 5$, and is the symmetric group $S_3$ if $n=4$;

\smallskip
-- we have a short exact sequence $0\to\dbZ/3\dbZ\hookrightarrow S_3\twoheadrightarrow\dbZ/2\dbZ\to 0$;

\smallskip
-- in the \'etale topology of $\Spec(F)$, the non-trivial torsors of $\dbZ/2\dbZ$ (resp. of $\dbZ/3\dbZ$) correspond to quadratic extensions of $F$ (resp. to Galois extensions of $F$ of degree 3);

\smallskip
-- for any embedding $i_F\colon F\hookrightarrow\dbR$, $G(F)_{\dbR}$ is an inner form of $\tilde G_{\dbR}$ ([De2, 1.2.6]).

\medskip
So if $n\Ge 5$ (resp. if $n=4$), then there is a totally real field extension $F_1$ of $F$ of degree at most 2 (resp. of degree at most 6) such that $G(F)_{F_1}$ is an inner form of $\tilde G_{F_1}$. As $G$ and $\tilde G$ are unramified over $\dbQ_p$, we can assume $F_1$ is unramified over $p$. 

\smallskip
{\bf ii)} If $(G,X)$ is of $D_n^{\dbR}$ type, it is worth making slightly more explicit the situation. We consider an inner form $G_F^\prime$ of $\tilde G_F$ with the property that w.r.t. any embedding $i_F\colon F\hookrightarrow\dbR$, ${G^\prime}_{\dbR}$ is isomorphic to $G(F)_{\dbR}$. For instance, we can take $G_F^\prime$ to be the adjoint group of the SO-group of the quadratic form $a_1x_1^2+a_2x_2^2+x_3^2+...+x_{2n}^2$ on $F^{2n}$, where $a_1,a_2\in F$ are such that (cf. approximation theory) we have: 

\medskip
-- under each embedding $i_F$ as above, they are both negative or positive, depending on the fact that the resulting group $G(F)_{\dbR}$ is or is not compact;

\smallskip
-- $G_F^\prime$ is unramified over $\dbQ_p$.
 
\medskip
Enlarging $F$ to $F_1$ in the way allowed by A), we can assume (see end of i)) that $G(F)_{F_1}$ is an inner form of $G_{F_1}^{\prime}$ and so of $\tilde G_{F_1}$. 

\smallskip
{\bf iii)} Let $\tilde G^e$ be the image of $\tilde G^d$ in $G(A)^{\ad}$ and let $h_m:\tilde G^e\hookrightarrow G(A)^{\ad}$ be the resulting monomorphism. We come back to the general situation of the Fact. We use the inner form it mentions to twist the homomorphism $h_{F_1}$ (more precisely $h_{mF_1}$); as it is mentioned below, if $(G,X)$ is of $D_n^{\dbR}$, with $n$ even, we might have to first replace $F_1$ by a totally real quadratic extension of it unramified above $p$. We get a monomorphism
$$h^t_m\colon G(F)^d_{F_1}\hookrightarrow G(A)^t,$$
where $G(A)^t$ is the resulting inner form of $G(A)_{F_1}$ and $G(F)^d_{F_1}$ is the central isogeny cover of $G(F)_{F_1}$ corresponding to the group cover $\tilde G^d$ of $\tilde G$. After potentially passing to a totally real quadratic extension of $F_1$, this twist is possible due to the following three facts.

\medskip
{\bf F1.} If $(G,X)$ is not of $D_n^{\dbR}$ type, then the center of $\tilde G^d$ is contained in the center of $G(A)$ and so $\tilde G^e$ is an adjoint group.

\smallskip
{\bf F2.} If $(G,X)$ is of $D_n^{\dbR}$ type, then (as we are working with the spin representation, cf. D)) $\tilde G^e$ is nothing else but $SO(2,2n-2)$. To argue this we use [Bou1, (VIII) of planche IV] and distinguish two cases.

{\bf n is odd.} In this case the half spin representations are dual to each other (see [Bou2, p. 210]). Moreover, $Z(\tilde G^d)$ is $\mu_4$ and so has only one subgroup $Z_1(\tilde G^d)$ of order 2. $\tilde G^e$ is the quotient of $\tilde G^d$ by $Z_1(\tilde G^d)$ and so it is $SO(2,2n-2)$.

{\bf n is even.} Referring to loc. cit. with $l:=n\in 2\dbN$, the kernels of the two half spin representations are naturally identified with $P(R)/Q(R)+\dbZ\overline{w}_{l-1}$ and respectively $P(R)/Q(R)+\dbZ\overline{w}_l$. So the kernel of the quotient morphism $\tilde G^{\sc}\twoheadrightarrow\tilde G^e$ is naturally identified with $P(R)/Q(R)+\dbZ\overline{w}_1$ and so with $SO(2,2n-2)$. 

\smallskip
{\bf F3.} The compact inner form $SO(2n)_{\dbR}$ of $SO(2,2n-2)_{\dbR}$ is obtained by twisting through an element of the set $H^1(\dbR,SO(2,2n-2)_{\dbR})$ (and not only of $H^1(\dbR,SO(2,2n-2)^{\ad}_{\dbR})$). So the obstruction of lifting the class in $H^1(F_1,\tilde G_{F_1})$ defining the inner twist $G(F)_{F_1}$ of $\tilde G_{F_1}$ to a class in $H^1(F_1,\tilde G^e_{F_1})$ is measured by a class $\gamma\in H^2(F_1,\mu_2)$ which over $\dbR$, under any embedding $F_1\hookrightarrow\dbR$, becomes trivial (i.e. becomes the class defined by matrix $\dbR$-algebras). Based on the standard Galois interpretation of $H^2(F,\mu_2)$ and on [Ha, 5.5.1], we conclude: by passing to a totally real quadratic extension of $F_1$ unramified above $p$, we can assume $\gamma$ is the trivial class. 

\medskip
So in fact we first twist $h_m{F_1}$, and then we ``lift" the things to the simply connected semisimple group cover $G(A)^t$ of the resulting twist of $G(A)^{\ad}_{F_1}$. However, in what follows, by abuse of language, we say we twist $h_{F_1}$; similarly, in other situations (like over $\dbR$) we use the same language.

\smallskip
{\bf G)} We have:

\medskip\noindent
{\bf Fact.} {\it If by extension of scalars under an embedding $i_{F_1}\colon F_1\hookrightarrow\dbR$, $G(F)_{F_1}$ becomes a compact group (resp. it becomes $\tilde G_{\dbR}$) then, under the same extension of scalars, $G(A)^t$ becomes a compact group (resp. it becomes $G(A)_{\dbR}$).}

\medskip
\proof
What we need to prove is the following thing. If we twist $h_{\dbR}$ through the compact inner form of $\tilde G^{\sc}_{\dbR}$, then the inner form $IF$ of $G(A)_{\dbR}$ we get is the compact form. As $h_{\dbR}$ satisfies the condition $(H_2)$ of [Sa], this is a particular case of a general fact pertaining to symmetric domains and which can be checked using Cartan decompositions (referred to in [Sa, p. 427]; see also [He, \S7 of ch. 3]). Here we would like to include a proof of this Fact, in a form needed for future references.

\medskip
{\bf Case 1: $(G,X)$ is of $C_n$ or $D_n^{\dbH}$ type.} This case can be proved in many ways. We choose two such ways. For both we use:

\medskip
-- the classical fact of the inner conjugacy of maximal compact subgroups of a semisimple $\dbR$-group (for instance, see [He, 2.2 (ii) of p. 256]), and

\smallskip
-- the classification of forms of simple, adjoint, split $\dbR$-groups (see [He, p. 517-518]). 

\medskip
First, looking at the inner forms of $SU(n,n)_{\dbR}$ (they are $SU(n_1,2n-n_1)_{\dbR}$, with $n_1\in \{0,...,n\}$), the only one which can contain a compact subgroup of $C_n$ or $D_n$ Lie type, is the compact inner form (corresponding to $n_1=0$). Second, Zarhin's lemma (recalled in 4.2.5.1) applies (to the extension to $\dbC$ of the inclusion of a maximal compact subgroup of $IF$ into $IF$ and to the standard faithful representation of $IF_{\dbC}$ of dimension $2n$); the fact that $IF$ has a compact subgroup of the same rank as itself (condition needed to be checked in order to apply this lemma) is an easy consequence of [Mi1, B.19]. 

\medskip
{\bf Case 2: $(G,X)$ is of $B_n$ or $D_n^{\dbR}$ type.} Let us first remark that the standard infinite sequence of monomorphisms 
$$SO(2,1)\hookrightarrow SO(2,2)\hookrightarrow ... \hookrightarrow SO(2,n)\hookrightarrow ...$$ extends to an infinite sequence of injective maps
$$(SO(2,1),X(1))\hookrightarrow (SO(2,2),X(2))\hookrightarrow ...\hookrightarrow (SO(2,n),X(n))\hookrightarrow ...\leqno (2)$$  
between Shimura pairs.
Here $X(n)$, $n\in\dbN$, are the logical Hermitian symmetric domains (they are uniquely determined, cf. [De2, 1.2.7-8]). Using the fact that $\tilde G^{\sc}_{\dbR}$ is the simply connected semisimple group cover of $SO(2,m)_{\dbR}$ for some $m\in\dbN$, $m\Ge 3$, we just need to show: twisting the resulting monomorphism (obtained by composing the inclusions of (2) with $h_{m\dbR}$) 
$$SO(2,1)_{\dbR}\hookrightarrow SU(2^{n-1},2^{n-1})_{\dbR}^{\ad}=G(A)_{\dbR}^{\ad},$$ 
through the compact inner form of $SO(2,1)_{\dbR}$, we get the compact inner form of $G(A)_{\dbR}^{\ad}$. 

So we need to show: under the natural diagonal embedding of the simply connected semisimple group cover $SU(1,1)_{\dbR}$ of $SO(2,1)_{\dbR}$ in $SU(q,q)_{\dbR}$ (i.e. under the composite of the diagonal monomorphism $SU(1,1)_{\dbR}\hookrightarrow SU(1,1)^q_{\dbR}$, with the natural monomorphism $SU(1,1)_{\dbR}^q\hookrightarrow SU(q,q)_{\dbR}$; here $q\in\dbN$), through its twist corresponding to the compact inner form of $SU(1,1)_{\dbR}$ we get the compact inner form of $SU(q,q)_{\dbR}$. This is an easy exercise which can be solved directly (by direct calculation, or by reduction of the situation to the $C_q$ type of Case 1, or by using induction and the same lemma of Zarhin, etc.). 
This ends the proof.   

\smallskip
{\bf H)} Let $G_1^{\sc}:=\Res_{F_1/\dbQ} G(A)^t$. From  G) we deduce the existence of an adjoint Shimura pair of the form $(G_1^{\ad},X_1^{\ad})$, with $G_1^{\ad}$ as the adjoint group of $G_1^{\sc}$; [De2, 1.2.8 (ii)] guarantees that such a Shimura pair is uniquely determined.

\smallskip
{\bf I)} We consider (cf. [Va2, 2.3.5.1]) an embedding 
$$f_1\colon (G_1,X_1)\hookrightarrow (G\Sp(W,\psi),S)$$ 
 such that the following things hold: 

\medskip
-- $G_1^{\der}=G_1^{\sc}$ (so the notations are justified);

\smallskip
-- there is a $\dbZ_{(p)}$-lattice $L_{(p)}$ of $W$ for which we get a perfect form $\psi\colon L_{(p)}\otimes_{\dbZ_{(p)}} L_{(p)}\to\dbZ_{(p)}$ and, denoting by $\scrB_{(p)}$ the $\dbZ_{(p)}$-subalgebra of elements of $\End(L_{(p)})$ fixed by $G_1$, a standard PEL situation $(f,L_{(p)},v,\scrB_{(p)})$, for any prime $v$ of $E(G_1,X_1)$ dividing $p$.  

\medskip
{\bf I')} [Va2, 2.3.5.1] is just a $p\Ge 3$ version of parts of [Va1, 6.5.1.1 and 6.6.5] and Erratum where we had $p\Ge 5$; moreover, these loc. cit. are rooted in [De2, 2.3.9-10]. For the convenience of the reader, we recall the parts of [Va1, 6.5.1.1 and 6.6.5] and Erratum which are relevant for what follows, by working rationally. Warning: in i) to iii) below we do not use the particular way $G_1^{\ad}$ was obtained; we just use the fact that $(G_1^{\ad},X_1^{\ad})$ is of some $A_m$ Lie type.

\smallskip
{\bf i)} 
Let $T(A)^t$ be a maximal torus of $G(A)^t$. Let $E_1$ be a finite field extension of $F_1$, Galois over $\dbQ$ and such that $T(A)^t_{E_1}$ is split. Let $T_1:=\Res_{F_1/\dbQ} T(A)^t$. It is a maximal torus of $G_1^{\sc}$. It splits over $E_1$. Let $W_{E_1}$ be a $E_1$-vector space of dimension $2d_1$ equal with $1$ plus the rank of $G(A)$ (so depending on the two Cases of G), $d_1$ is $n$ or $2^{n-1}$); when we view it as a $\dbQ$--vector space, we denote it by $\tilde W_{E_1}$. Identifying $G(A)^t_{E_1}$ with $SL(W_{E_1})$, we consider a monomorphism
$$m_1:G^{\sc}_1\hookrightarrow SL(\tilde W_{E_1})$$
which is the composite of the natural embeddings
$$G^{\sc}_1\hookrightarrow \Res_{E_1/\dbQ} G(A)^t_{E_1}$$ and 
$\Res_{E_1/\dbQ} G(A)^t_{E_1}\hookrightarrow SL(\tilde W_{E_1})$ (the second one is defined via the mentioned identification). Let $W:=\tilde W_{E_1}\oplus \tilde W_{E_1}$. We consider the ``sum" of $m_1$ and of the monomorphism $m_1^*:G_1^{\sc}\hookrightarrow SL(\tilde W_{E_1})$ obtained by composing $m_1$ with the automorphism $SL(\tilde W_{E_1})\tilde\to SL(\tilde W_{E_1})$ which at the level of $2d_1\times 2d_1$ matrices takes $M$ into $(M^t)^{-1}$. We get a monomorphism
$$m_2:G^{\sc}_1\hookrightarrow SL(W).$$
\indent
{\bf ii)} 
Let (as in the proof of [De3, 2.3.10]) $\scrS_1$ be the set of extremal points of the Dynkin diagram $\scrD_1$ of $\Lie(G_{1\dbC}^{\sc})$ w.r.t. $\Lie(T_{1\dbC})$ (and some fixed Borel subalgebra of $\Lie(G_{1\dbC}^{\sc})$ normalized by $T_{1\dbC}$). $\Gal(\dbQ)$ acts naturally on $\scrS_1$ and so we can identify $\scrS_1$ with the $\Gal(\dbQ)$-set of $\dbQ$--homomorphisms from $K_{\scrS_1}$ to $\dbC$, for $K_{\scrS_1}$ a totally imaginary quadratic extension of $F_1$. $\Res_{K_{\scrS_1}/\dbQ} \dbG_m$ acts naturally on $W$ as follows. We can work over $\dbC$. If $W_0$ is a simple $G^{\sc}_{1\dbC}$-submodule of $W\otimes_{\dbQ} \dbC$ corresponding to the fundamental weight associated to some node $s_1\in\scrS_1$, then the extension of $\Res_{K_{\scrS_1}/\dbQ} \dbG_m$ to $\dbC$ acts on it via its character corresponding naturally to $s_1$. 
Let $\tilde G_1$ be the reductive subgroup of $GL(W)$ generated by $G^{\sc}_1$ and by $\Res_{K_{\scrS_1}/\dbQ} \dbG_m$. 

\smallskip
{\bf iii)} [Va1, 6.5.1.1 and 6.6.5] and Erratum explain how to construct a monomorphism $h_x:\Res_{\dbC/\dbR} \dbG_m\hookrightarrow \tilde G_1$ such that the Hodge $\dbQ$--structure on $W$ we get has type $\{(-1,0),(0,-1)\}$.  
[De2, 2.3.3] shows the existence of a perfect alternating form $\psi$ on $W$ such that $h_x$ factors through $G_{1\dbR}$, where $G_1$ is a reductive subgroup of $\tilde G_1$ normalizing $\psi$, containing the maximal subgroup of $\Res_{K_{\scrS_1}/\dbQ} \dbG_m$ which over $\dbR$ is compact and having $G_1^{\sc}$ as its derived subgroup. [Va1, 6.6.5.1] and Erratum point out that we can assume as well that the resulting injective map $f_1$ is a PEL type embedding and that a lattice $L_{(p)}$ as mentioned in I) does exist.

\smallskip
{\bf iv)} All factors of $G_{1\dbR}^{\ad}$ are either compact or are $SU(d_1,d_1)_{\dbR}$ group. So $(G_1^{\ad},X_1^{\ad})$ is without involution. If all simple factors of $G_{\dbR}$ are non-compact, then all simple factors of $G_{1\dbR}^{\ad}$ are non-compact; so (cf. [De2, 2.3.13] and Erratum) $h_x$ factors through the subgroup of $G_1$ generated by $G^{\sc}_1$ and $Z(GL(W))$.

\smallskip
{\bf J)} We come back to I). There is $M(G,X)\in\dbN$ effectively computable and such that we can assume
$$\dim_{\dbQ}(W)\Le M(G,X),$$ 
regardless of the prime $p\Ge 3$ we keep track of. Its existence is a consequence of:

\medskip
$\bullet$ the review of 3.1 I') ($[E_1:\dbQ]$ is bounded in terms of $\dim_{F_1}(T(A)^t)$), and of 

\smallskip
$\bullet$ the fact that we can choose (cf. i) and iii) of F)) $F_1$ so that we have $[F_1:F]\Le 12$.  
\finishproclaim

\Proclaim{3.2. Proposition.} There are injective maps
$$(G_2,X_2)\hookrightarrow (G_3,X_3)\hookrightarrow (G_1,X_1)$$
such that:

\medskip
{\bf a)} $(G_2^{\ad},X_2^{\ad})=(G,X)$ and $(G_3^{\ad},X_3^{\ad})=(G^{F_1},X^{F_1})$;

\smallskip
{\bf b)} we get a monomorphism $G_2^{\ad}\hookrightarrow G_3^{\ad}$ which (under the identifications of a)) is nothing else but the natural monomorphism $G\hookrightarrow G^{F_1}$;

\smallskip
{\bf c)} $G_3^{\der}=\Res_{F_1/\dbQ} G^d_{F_1}$, and under this identification and of the one of 3.1 H), the monomorphism $G_3^{\der}\hookrightarrow G_1^{\der}$ becomes the Weil restriction from $F_1$ to $\dbQ$ of $h^t_m$;

\smallskip
{\bf d)} if $(G,X)$ is not (resp. is) of $D_n^{\dbH}$ type, then we have $E(G_2,X_2)=E(G_3,X_3)=E(G_1,X_1)$ (resp. $E(G_2,X_2)=E(G_3,X_3)$);

\smallskip
{\bf d')}  if $(G,X)$ is of $D_n^{\dbH}$ type with involution, then $E(G_2,X_2)=E(G_3,X_3)$ is an extension of $E(G_1,X_1)$ of degree at most $2$;

\smallskip
{\bf e)} $Z^0(G_2)$ and $Z^0(G_3)$ are subtori of $Z^0(\tilde G_1)$;

\smallskip
{\bf f)} if $(G,X)$ is not of $D_n^{\dbH}$ type, then $G_2^{\der}$ and $G_3^{\der}$ are simply connected. 
\finishproclaim

\smallskip
\proof
We take $G_3$ to be generated by $Z^0(\tilde G_1)$ and by the subgroup $\Res_{F_1/\dbQ} G(F)_{F_1}^d$ of $G_1^{\der}$. This takes care of c). We take $G_2$ to be the subgroup of $G_3$ having $Z^0(G_1)$ as a subgroup and such that b) holds. So e) holds as well.

$G_1(\dbQ)$ permutes the connected components of $X_1$, cf. [Ko, 4.3]; this implies $X_1=X_1^{\ad}$. So the fact that $h_{\dbR}$ at the level of Lie algebras satisfies the condition $(H_2)$ of [Sa, p. 427] gets translated (cf. 3.1 G) and H)) in: there is $x_1\in X_1$ such that its attached monomorphism $\Res_{\dbC/\dbR} \dbG_m\hookrightarrow G_{1\dbR}$ factors through $G_{2\dbR}$ in such a way that the resulting homomorphism $\Res_{\dbC/\dbR} \dbG_m\to G_{\dbR}$ is an element of $X$. This takes care of a). 

d) is a easy consequence of [De1, 3.8] and [De2, 2.3.6]. f) is a direct consequence of the construction of $h^t_m$. This ends the proof.

\Proclaim{3.2.0. The $\dbZ_{(p)}$-version.} \rm
It is now easy to see (as in [Va1, 6.5-6]) that $h^t_m$ has a version over the normalization ${F_1}_{(p)}$ of $\dbZ_{(p)}$ in $F_1$. First, referring to 3.1 E), instead of $\dbQ$, a totally imaginary quadratic extension $K$ of $\dbQ$ unramified over $p$, and a quaternion $\dbQ$--algebra $H_{\dbQ}$ containing $K$ and ramified at the infinite place, we need to work respectively with: $\dbZ_{(p)}$, the normalization $K_{(p)}$ of $\dbZ_{(p)}$ in $K$, and with a maximal order of $H_{\dbQ}$ containing $K_{(p)}$. Typical situation: we work with $\dbZ_{(p)}$, with $\dbZ_{(p)}(i)$ and with $\dbZ_{(p)}(i,j)$. Second, the fact that ${F_1}_{(p)}$ is \'etale over $\dbZ_{(p)}$ implies that by the twisting process of iii) of 3.1 F), the extension of $h^t_{m}$ to $F_1\otimes_{\dbQ} \dbQ_p$ extends to a homomorphism between reductive groups over ${F_1}_{(p)}\otimes_{\dbZ_{(p)}} \dbZ_p$; this homomorphism is a monomorphism, cf. [Va1, 3.1.2.1 c)]. So $h^t_m$ extends to a monomorphism of reductive groups over ${F_1}_{(p)}$, cf. [Va1, 3.1.3.1]. 

Based on [Va2, 2.3.5.1] (cf. also 3.1 I)) we get: we can assume that the Zariski closure of $G_2$ and of $G_3$ in $GL(L_{(p)})$ are reductive groups over $\dbZ_{(p)}$.  So $H_i:=G_i(\dbQ_p)\cap GL(L_{(p)})(\dbZ_{p})$ is a hyperspecial subgroup of $G_i(\dbQ_p)$, $i=\overline{1,3}$.
\finishproclaim

\Proclaim{3.2.1*. Corollary.}
The result stated in [Va2, 1.15.8] is true for $(G,X)$ provided 
$$p\Ge\max\{5,M(G,X)/2\}.$$ 
\finishproclaim

\proof
Let $i\in\{2,3\}$. We denote by $f_i$ the inclusion $(G_i,X_i)\hookrightarrow (GSp(W,\psi),S)$. The inequality $p\Ge\max\{5,M(G,X)/2)\}$ conforms with [Va1, 5.8.6]. The fact that each triple $(f_i,L_{(p)},v_i)$, with $v_i$ a prime of $E(G_i,X_i)$ dividing $p$, is a SHS is a consequence of 3.1 J) and of loc. cit. (cf. also [Va2, 2.3.6]). So 3.2 and 3.2.0 take care of [Va2, 1.15.8 1) to 4)]. 

The rest of the proof is just sketched. [Va2, 1.15.8 5)] (i.e. the part involving ``$G$-ordinary" things) is a consequence of [Va2, 4.11.2-3]: the equivalent conditions of [Va2, 4.11.2] can be easily checked; see [Va7] for details and for specifications on the part ``other (simple to check) conditions" of [Va2, 1.15.8 5)]. This ends the proof.

\medskip
Warning: if we assume [Va2, 2.3.8 2)] (see [Va7]), then ``$p\Ge\max\{5,M(G,X)\}$" gets replaced by ``$p\Ge 3$".

\Proclaim{3.2.2. Remarks.} \rm
{\bf 1)} If $(G,X)$ is of $A_n$ type ($n\in\dbN$), then 3.2.1 gets substituted by [Va2, 2.3.5.1]. We repeat: i) to iii) of 3.1 I') review as well the rational variant of loc. cit.

{\bf 2)} Using 2.6 2) we get that 3.2 has a version for the case when $(G,X)$ is an arbitrary Shimura pair of adjoint, abelian type. To state it, let $(G,X)=\prod_{i\in I} (G_i,X_i)$ be written as a product of such simple Shimura pairs. Let $F^i:=F(G_i,X_i)$ and let $\Gtil_i$ be the adjoint $F^1$-group such that $G_i=\Res_{F^i/\dbQ} \Gtil_i$ (cf. 2.5).

Then there is a product $(G^1,X^1)=\prod_{i\in I} (G_i^1,X_i^1)$ of simple, adjoint Shimura pairs of some $A_n$ type ($n\in\dbN$ depends on $i\in I$), and there are injective maps
$$(G_2,X_2)\hookrightarrow (G_3,X_3)\hookrightarrow (G_1,X_1)\operatornamewithlimits{\hookrightarrow}\limits^{f_1} (GSp(W_1,\psi_1),S_1)$$
such that:

\medskip
{\bf a)} $(G_2^{\ad},X_2^{\ad})=(G,X)$ and (cf. 3.1 A)) $(G_3^{\ad},X_3^{\ad})=\prod_{i\in I} (G_i^{F^i_1},X_i^{F^i_1})$, with $F_1^i$ a totally real number field containing $F^i$;

\smallskip
{\bf b)} we get a monomorphism $G_2^{\ad}\hookrightarrow G_3^{\ad}$ which is a product of the natural monomorphisms $G_i\hookrightarrow G_i^{F_1^i}$;

\smallskip
{\bf c)} the monomorphism $G_3^{\der}\to G_1^{\der}$ is a product indexed by $i\in I$ of monomorphisms which are either isomorphisms or are obtained as in 3.2 c), depending on the fact that $G_i$ is or is not of $A_{n_i}$ Lie type for some $n_i\in\dbN$;

\smallskip
{\bf d)} $(G_1^{\ad},X_1^{\ad})=(G^1,X^1)$, $G_1^{\der}$ is a simply connected semisimple group, and $f_1$ is a PEL type embedding which (cf. 3.2 and 2.6 2)) is a Hodge quasi product indexed naturally by elements of $I$ of PEL type embeddings;

\smallskip
{\bf e)} if $\forall i\in I$ (resp. if there is $i\in I$ such that) $(G_i,X_i)$ is of $D_n^{\dbH}$ type, then $E(G_2,X_2)=E(G_3,X_3)=E(G_1,X_1)$ (resp. $E(G_2,X_2)=E(G_3,X_3)$);

\smallskip
{\bf f)} $Z^0(G_2)$ and $Z^0(G_3)$ are subtori of $Z^0(G_1)$;

\smallskip
{\bf g)} if none of the Shimura pairs $(G_i,X_i)$ is of some $D_n^{\dbH}$ type, then $G_2^{\der}$ and $G_3^{\der}$ are simply connected semisimple groups.

\medskip
{\bf 3)} Using 2) we get that the result stated in [Va2, 1.15.8] remains true for any Shimura quadruple $(\Gbar,\Xbar,\Hbar,\vbar)$ of adjoint, abelian type, with $\vbar$ dividing a prime $p\Ge\max\{5,M(\Gbar,\Xbar)\}$; here $M(\Gbar,\Xbar)\in\dbN$ is effectively computable (for instance, we can take it --cf. [Va1, 6.5.1]--  to be the product of the numbers $M(*)$ of 3.1 J), with $*$ running through the simple factors of $(\Gbar,\Xbar)$). 

{\bf 4)} 3.2.1 has many consequences (cf. [Va1, 6.8], [Va2, 1.15.8, 1.15.8.1 and 4.4.6]). We will come back to 3.2.1 and its consequences in [Va7]. In what follows we do not need 3.2.1: we just need just 3.2 and its variant of 2). 

{\bf 5)} We view 3.2 and 3.2.2 2) as the (refined) unitary version of [De2, 2.3.10]. Moreover, combining 3.2 and 3.2.2 2) with the PEL type embeddings of Shimura (see [Sh]), one can reobtain the essence of [De2, 2.3.10] in many situations (in particular, one can reobtain in many situations the fact that any Shimura variety of $A_n$, $B_n$, $C_n$, $D_n^{\dbH}$ or $D_n^{\dbR}$ is of preabelian type). However, as pointed out in [Va2, 4.6.11 C], we do not favor the approach to [De2, 2.3.10] via [Sh].

{\bf 6)} In 3.2 we can replace $F_1$ by any other totally real number field containing it (and, in case we want 3.2.1, unramified above $p$).
\finishproclaim

\Proclaim{3.2.3. Definition.} \rm
We refer to (2) of 3.1 G) as the standard SO-sequence of injective maps between Shimura pairs (to be abbreviated shortly as SSO).
\finishproclaim

\Proclaim{3.2.4. Point of view.} \rm
It is a widely spread opinion that the closest generalization of the elliptic modular curve is provided by Siegel or Hilbert--Blumenthal  modular varieties. In our opinion this is just partially correct: the Shimura varieties defined by Shimura pairs of SSO are equally well a very close generalization of the elliptic modular curve. We refer to these Shimura varieties as classical spin modular varieties. 
\finishproclaim

\Proclaim{3.3. Applications of 3.2 to Frobenius tori.} \rm
We start by introducing a new setting.

\smallskip
{\bf A)} We use the notations of 2.2.2. We still denote by $(\scrA,\scrP_{\scrA})$ its restriction to $\Sh_H(G,X)$. Let 
$$m\colon\Spec(E)\to\scrN$$ 
be a morphism. Let $A:=m^*(\scrA)$. We assume that the Mumford--Tate group of $A$ is $G$.

\smallskip
{\bf B)} Let 
$$(G_2,X_2)\hookrightarrow (G_3,X_3)\hookrightarrow (G_1,X_1)\operatornamewithlimits{\hookrightarrow}\limits^{f_1} (GSp(W_1,\psi_1),S_1)$$
be injective maps constructed following the pattern of 3.2 but starting from the adjoint Shimura pair $(G_2^{\ad},X_2^{\ad}):=(G^{\ad},X^{\ad})$ (cf. 3.2.2 2)). So $G_1^{\ad}$ has all its simple factors of $A_n$ Lie type (for potentially different values of $n$) and the embedding $(G_1,X_1)\hookrightarrow (GSp(W_1,\psi_1),S_1)$ is a PEL type embedding. Let $\scrB$ be the semisimple $\dbQ$--subalgebra of $\End(W_1)$ formed by elements fixed by $G_1$. Let 
$$f_2:(G_2,X_2)\hookrightarrow (GSp(W_1,\psi_1),S_1)$$
be the composite map.

We can assume that the intersection $X_{20}:=X_2\cap X$ inside $X_2^{\ad}=X^{\ad}$ is non-empty, cf. [Va1, 2.4.0 and rm. 3) of 3.2.7] and [Va2, 4.9.1]. As $G(\dbQ)$ is dense in $G(\dbR)$, we can assume that $m$ can be lifted to a point $\Spec(\dbC)\to\Sh(G,X)$ defined by an equivalence class $[x_m,g_m]$ (see [Va1, 4.1]), with $x_m\in X_{20}^0$; here $X_{20}^0$ is an arbitrary (but fixed) connected component of $X_{20}$ and $g_m\in G(\dbA_f)$. Let $L_1$ be a $\dbZ$-lattice of $W_1$ such that we get a perfect form $\psi_1\colon L_1\otimes_{\dbZ} L_1\to\dbZ$. Let $H_2:=\{\text{g}\in G_2(L_1\otimes_{\dbZ} \dbZhat)|\text{g mod 4 is the identity}\}$.
As in 2.2.2, we get a principally polarized abelian scheme $(\scrA_2,\scrP_{\scrA_2})$ over $\Sh_{H_2}(G_2,X_2)$.

\smallskip
{\bf C)} Let $(W_0,\psi_0):=(W\oplus W_1,\psi\oplus\psi_1)$. We consider as in  [Va2, 4.9.2] a Hodge twist $(G_4,X_4)$ of $(G,X)$ and $(G_2,X_2)$. It is equipped with a (diagonal) injective map 
$$f_3\colon (G_4,X_4)\hookrightarrow (GSp(W_0,\psi_0),S_0).$$ 
$f_3$ is the composite of an injective map 
$$f_4\colon (G_4,X_4)\hookrightarrow (G,X)\times_{\scrH} (G_2,X_2)$$ 
with a Hodge quasi product $f\times^{\scrH} f_2$. 
$f_3$ gives birth to identifications (to be denoted by $ID$) 
$$(G_4^{\ad},X_4^{\ad})=(G^{\ad},X^{\ad})=(G_2^{\ad},X_2^{\ad}).$$ 
These identifications $ID$ are defined via the composite of $f_4$ with the natural injective map $i_2:(G,X)\times_{\scrH} (G_2,X_2)\hookrightarrow (G,X)\times (G_2,X_2)$ and with the
projections of $(G,X)\times (G_2,X_2)$ onto its factors; they conform with the identification of B). Warning: we take $X_4$ such that $X_{20}^0$ is a connected component of it.

Let $H_4:=G_4(\dbQ_p)\cap H\times H_2$. Let $L_0:=L\oplus L_1$; it is a $\dbZ$-lattice of $W_0$ and we have a perfect form $\psi_0\colon L_0\otimes_{\dbZ} L_0\to\dbZ$. We get naturally (starting from $f_3$ and $L$), as in 2.2.2, a universal principally polarized abelian scheme $(\scrA_4,\scrP_{\scrA_4})$ over $\Sh_{H_4}(G_4,X_4)$. 

\smallskip
{\bf D)} We assume in what follows that the reader is familiar with the part of \S4 running from 4.0 to 4.2.7.1, with 4.2.8, 4.2.11 and 4.2.14. We also assume that $E$ contains $E(G_4,X_4)$ and that it is large enough so that we have a morphism
$$m_4\colon\Spec(E)\to\Sh_{H_4}(G_4,X_4),$$
whose projection on $\Sh_H(G,X)$ (via $i_2\circ f_4$) is $m$. The assumptions $x_m\in X_{20}^0\subset X_4$ guarantee the existence of such a morphism $m_4$, with $E$ large enough, provided we pass to isogenies and change a little bit the level structures (to be compared with the definition of $\grC_0$ in [Va2, 4.9.1.1]). Though we do not need this, we would like to point out that [Va1, 3.3.1] guarantees that this passage and these changes of level structures can be performed so that we still keep track of any a priori given odd prime $q$ for which the Zariski closure of $G$ (resp. of $G_2$) in $GL(L\otimes_{\dbZ} \dbZ_{(q)})$ (resp. in $GL(L_1\otimes_{\dbZ} \dbZ_{(q)})$) is a reductive group over $\dbZ_{(q)}$. What we do need here: these two operations are ``irrelevant" in connection to Mumford--Tate conjecture (see 4.2.11); so, not to complicate the notations, we assume that the performance of these two operations is not needed. Let 
$${}_2A$$ 
be the abelian variety over $E$ obtained from $\scrA_2$ by pull back via $m_4$, $i_2\circ f_4$ and $f_2$ (we use lower left index $2$ not to create confusion with the $A_2$ Lie type). So we have
$$m_4^*(\scrA_4)=A\times_E {}_2A.$$ 
We also assume $E$ is large enough so that all Hodge cycles of $A_{\Ebar}$, $A_{\Ebar}^2$ and $A_{\Ebar}\times A_{\Ebar}^2$ are defined over $E$.

\smallskip
{\bf E)} Let $v=v(A,E)$ be a prime of $E$ such that: 

\medskip
{\bf i)} $A$ and ${}_2A$ have good reductions $A_v$ and respectively ${}_2A_v$ w.r.t. $v$;

\smallskip
{\bf ii)} it is unramified over an odd rational prime $p(A)\Ge M(f_1,L_1)$ (see 2.2.2.1);

\smallskip
{\bf iii)} for any prime $l\Ge p(A)$ we get standard Hodge situations $(f_3,L_0\otimes_{\dbZ} \dbZ_{(l)},v_4(l))$, $(f_2,L_1\otimes_{\dbZ} \dbZ_{(l)},v_2(l))$ and $(f_1,L_1\otimes_{\dbZ} \dbZ_{(l)},v_1(l))$ (for $i\in\{1,2,4\}$, $v_i(l)$ is an arbitrary prime of $E(G_i,X_i)$ dividing $l$);

\smallskip
{\bf iv)} for any prime $l\Ge p(A)$, $\scrB_{(l)}:=\scrB\cap\End(L_1\otimes_{\dbZ} \dbZ_{(l)})$ is a semisimple $\dbZ_{(l)}$-algebra which is self dual w.r.t. $\psi_1$ (so the quadruple $(f_1,L_1\otimes_{\dbZ} \dbZ_{(l)},v_1(l),\scrB_{(l)})$ is a standard PEL situation).

\medskip
The existence of $v$ is implied by [Va1, 5.8.7] and [Va2, 2.3.6]: they take care of iii), while i), ii) and iv) are obviously true for $p$ big enough (for ii) cf. also [Va2, 4.6.13]). The main conditions are i), iii) and iv); strictly speaking, we need iii) and iv) just for $l=p(A)$. Warning: the use of ii) can be entirely avoided. 

\smallskip
{\bf F)} In all that follows we denote $p(A)$ just by $p$; also, let $v_i(p)$ be the prime of $E(G_i,X_i)$ divided by $v$, $i\in\{1,2,4\}$. Let $l$ be a rational prime different from $p$. In what follows we work with sufficiently high, positive, integral powers of the Frobenius automorphism of $\overline{k(v)}$ fixing $k(v)$, so that all Frobenius elements to be introduced below define valued points of the Frobenius tori involved (i.e. we prefer to work with tori instead of groups of multiplicative type). 

Let $Fr_i\in T_v^i(\dbQ_l)$, $i\in\{2,4\}$, be a Frobenius element obtained as in 4.2.7.1. Similarly, we get a Frobenius element $Fr\in T_v(\dbQ_l)$. The choice of a Frobenius automorphism of a prime of $\bar E$ dividing $v$, results in natural monomorphisms $T_{v\dbQ_l}\hookrightarrow G_{\dbQ_l}$ and $T_{v\dbQ_l}^i\hookrightarrow G_{i\dbQ_l}$. We can assume that under the identifications $ID$, the images of $Fr$ and $Fr_2$ in $G^{\ad}_{\dbQ_l}(\dbQ_l)$ are the same.
\finishproclaim

\Proclaim{3.3.1. CM lifts.} \rm
The key point is: as the embedding $f_1$ is a PEL type embedding, with $G_1^{\der}$ a simply connected semisimple group which is a product of semisimple groups of $A_n$ Lie type, the Langlands--Rapoport conjecture is true for the SHS $(f_1,{L_0}\otimes_{\dbZ} \dbZ_{(p)},v_1(p))$, (cf. [Mi4, 6.12] and 2.2.2.1; the reference to 2.2.2.1 can be substituted by [Mi4, (0.2) c)], cf. also [Va2, 4.6.13]). What we need from this key point in our context is: 

\medskip\noindent
{\bf Fact.} {\it There is a torus $T_1$ of $G_1$ and an element $Fr_{\dbQ}\in T_1(\dbQ)$ such that there is $g\in G_1(\dbQ_l)$ taking (under inner conjugation) $T_{1\dbQ_l}$ into $T_{v\dbQ_l}^2$ and $Fr_{\dbQ}$ into $Fr_2$.}

\medskip
As the manuscript [Mi4] is not yet published, this Fact is not used in what follows. However, for our present needs, instead of [Mi4] we can use as well the main result of [Zi] as follows. A suitable abelian variety ${}_2A_{\vbar}^{\prime}$ which is $\dbZ[{1\over p}]$-isogeneous (in the sense of [Va2, 2.1]) to ${}_2A_{v\overline{k(v)}}$, has a lift ${}_2A^{\prime}$ to a complete DVR of mixed characteristic in such a way that (see [Zi, 4.4]):

\medskip\noindent
{\bf CM.} {\it ${}_2A^{\prime}$ has complex multiplication.}

\smallskip\noindent
{\bf ISOG.} {\it The $\dbQ$--endomorphisms of ${}_2A_{\vbar}^{\prime}$ corresponding to elements of $\scrB$ and the $\dbQ$--polarization of ${}_2A_{\vbar}^{\prime}$ (obtained naturally from $\scrP_{\scrA_2}$) lift to $\dbQ$--endomorphisms and respectively to a $\dbQ$--polarization $p_{{}_2A^{\prime}}$ of ${}_2A^{\prime}$. Moreover, the elements of $\scrB_{(p)}$ correspond to $\dbZ_{(p)}$-endomorphisms of ${}_2A^{\prime}$ and $p_{{}_2A^{\prime}}$ is in fact a principal polarization.} 

\medskip
Using the realization of the triple $({}_2A^{\prime},p_{{}_2A^{\prime}},\scrB)$ in the Betti homology with coefficients in $\dbQ$ (of any pull back of ${}_2A^{\prime}$ to an abelian variety over $\dbC$), we get a symplectic space $(W_1^\prime,\psi_1^\prime)$ over $\dbQ$ and a $\dbQ$--monomorphism $\scrB\hookrightarrow\End(W_1^\prime)$. Let $G_1^\prime$ be the subgroup of $GSp(W_1^\prime,\psi_1^\prime)$ fixing all elements of $\scrB$; it is an inner form of $G_1$. More precisely, the ``difference" between the isomorphism classes of the following two triples $(W_1^\prime,\psi_1^\prime,\scrB)$ and $(W_1,\psi_1,\scrB)$ is ``measured" by a class $\gamma\in H^1(\dbQ,G_1)$ (for instance, this can be seen working in the \'etale cohomology with $\dbQ_l$-coefficients). We say that this class is trivial if the two triples are in fact isomorphic; the same applies if we replace $\dbQ$ by (i.e. if we tensor $W_1^\prime$ and $W_1$ with) any other field containing $\dbQ$.
We have:

\medskip\noindent
{\bf HW.} {\it Let $q$ be a rational prime. The images of $\gamma$ in $H^1(\dbQ_q,G_{1\dbQ_q})$ and $H^1(\dbR,G_{1\dbR})$ are trivial.}

\medskip
This is a consequence of the first four paragraphs of [Ko, \S8]; it is iv) of 3.3 E) which allows us to refer to [Ko]. The fact that $\scrB$ is not necessarily a $\dbQ$--simple algebra does not represent an impediment as, by using $\dbQ$--simple factors of $\scrB$, the situation gets reduced to the one used in [Ko, \S5 and \S8] (cf. constructions in ii) of 3.1 I') and d) of 3.2.2 2)). It is $ISOG$ which takes care of the case $q=p$ (cf. loc. cit.; to be compared as well with the part of the proof of [Va2, 4.12.12] referring to a $\dbQ_p$-context).

\smallskip
In particular, $HW$ implies: we have an identification (unique up to inner conjugation)
$$G_{1\dbQ_q}=G_{1\dbQ_q}^{\prime}$$ 
as well as a canonical identification
$$Z(G_1)=Z(G_1^\prime).$$
In fact we can assume $G_1^\prime=G_1$ (cf. [Mi4, 6.12]; see also the first five paragraphs of [Ko, \S7] for the cases when $G^{\ad}$ has no simple factors of $A_{2s+1}$ Lie type, with $s\in\dbN$); in what follows we will not use this identification at all, due to two main reasons. The first reason is: the ideas to be introduced below can be used as well in many other cases involving PEL type embeddings (which are not necessarily pertaining to $A_s$ Lie types); we recall that [Va1, 1, 7.1] does not ``handle" fully the D case (see [Ko]) of the PEL type embeddings and these points out that these ideas might be used later on. 
The second reason is: for applications to \S5 we need just the part of HW involving ($\dbQ_p$ and) $\dbR$. 

To prepare the background for these applications, let $(G_1^\prime,X_1^\prime)$ be the Shimura pair attached to ${}_2A^{\prime}$. From $HW$ and [Ko, 4.2] we deduce the existence of an isomorphism
$$W_1\otimes_{\dbQ} \dbR\arrowsim W_1^\prime\otimes_{\dbQ} \dbR$$ 
such that:

\medskip
{\bf a)} it takes $\psi$ into $\psi^\prime$ and $b\in\scrB$ into $b$, $\forall b\in\scrB$ (so it takes the subgroup $G_{1\dbR}$ of $GL(W_1\otimes_{\dbQ} \dbR)$ into the subgroup $G_{1\dbR}^\prime$ of $GL(W_1^\prime\otimes_{\dbQ} \dbR)$);

\smallskip
{\bf b)} it takes $X_1$ into $X_1^\prime$. 

\medskip
In particular, we get that the 0 dimensional Shimura varieties $\Sh(G_1^{\ab},X_1^{\ab})$ and $\Sh(G_1^{\prime\ab},X_1^{\prime\ab})$ are isomorphic. Even more particular, we get:

\medskip\noindent
{\bf AB.} {\it The smallest torus $Z_{\text{min}}$ of $Z(G_1)$ such that all monomorphisms $\Res_{\dbC/\dbR} \dbG_m\hookrightarrow G_{1\dbR}$ defined by elements of $X_1$, factor through the extension to $\dbR$ of the subgroup of $G_1$ generated by $G_1^{\der}$ and by $Z_{\text{min}}$, is the same as the smallest torus $Z_{\text{min}}^\prime$ of $Z(G_1^\prime)$ with the property that all monomorphisms $\Res_{\dbC/\dbR} \dbG_m\hookrightarrow G_{1\dbR}^\prime$ defined by elements of $X_1^\prime$, factor through the extension to $\dbR$ of the subgroup of $G_1^\prime$ generated by $G_1^{\prime\der}$ and by $Z_{\text{min}}^\prime$.}

\medskip
Above, HW points out to a Hasse--Witt type property, while AB stands for an abelian counterpart picture.

Coming back to our Frobenius tori, from $CM$ and $ISOG$ we deduce (cf. the third paragraph of 4.2.7.1; above, we can work as well with a finite field extension of $k(v)$ instead of $\overline{k(v)}$) the existence of a torus $T_1^\prime$ of $G_1^\prime$ and of an element $Fr_{\dbQ}^\prime\in T_1^\prime(\dbQ)$ such that:

\medskip\noindent
{\bf Fact'.} {\it There is $g^\prime\in G_1(\dbQ_l)$ taking (under conjugation) $T_{1\dbQ_l}^\prime$ into $T_{v\dbQ_l}^2$ and $Fr^\prime_{\dbQ}$ into $Fr_2$.}
\finishproclaim

\Proclaim{3.3.1.1. Remark.} \rm
Directly from the form of the monomorphisms $G_{2\dbR}^{\der}\hookrightarrow G_{3\dbR}^{\der}\hookrightarrow G_{1\dbR}^{\der}$ (see c) of 3.2.2 2)) we get: each simple factor of $\Lie(G_{3\dbR}^{\der})$ is absolutely simple and contained in a uniquely determined absolutely simple factor of $\Lie(G_{1\dbR}^{\der})$; moreover, each simple factor $\Lie(\scrF)$ of $\Lie(G_{2\dbR}^{\der})$ is embedded (diagonally) in a product of simple factors of $\Lie(G_{3\dbR}^{\der})$ which are all isomorphic to $\Lie(\scrF)$. The same can be stated over $\overline{\dbQ_l}$ instead of over $\dbR$.
\finishproclaim

\Proclaim{3.3.1.2. Two important things. A.} \rm
Two simple factors of $G_{1\overline{\dbQ_l}}^{\ad}$ are connected over $\dbQ$ through $G_1^{\ad}$ (i.e. are factors of the extension to $\overline{\dbQ_l}$ of the same simple factor of $G_1^{\ad}$) iff they are connected over $\dbQ$ through $G_1^{\prime{\ad}}$; here we use the identification $G_{1\overline{\dbQ_l}}^{\ad}=G_{1\overline{\dbQ_l}}^{\prime{\ad}}$. This is a consequence of the fact that $G_1^\prime $ is an inner form of $G_1$.

\smallskip
{\bf B.} In other order of ideas, from the construction of Frobenius tori, we get: the Zariski closure in $G_1^{\prime\ad}$ of the images of the integral powers of $Fr_{\dbQ}^\prime$ in $G_1^{\prime\ad}(\dbQ)$ is the image $T_1^{\prime\ad}$ of $T_1$ in $G_1^{\prime\ad}$. 
So these images, when viewed as $\overline{\dbQ_l}$-valued points of $T_{1\overline{\dbQ_l}}^{\prime\ad}$, are Zariski dense in $T_{1\overline{\dbQ_l}}^{\prime\ad}$. From 3.3.1 Fact' we deduce: the image $T_{v\overline{\dbQ_l}}^{\ad}$ of $T_{v\overline{\dbQ_l}}$ in $G^{\ad}_{\overline{\dbQ_l}}$ is equal to the conjugate of $T_{1\overline{\dbQ_l}}^{\prime\ad}$ through $g^\prime\in G_1(\overline{\dbQ_l})$. 
\finishproclaim 

Based on 3.3.1.2 B, 3.3.1.1 and 3.3.1 Fact' we conclude:

\Proclaim{3.3.2. Theorem.}
If two simple factors of $G^{\ad}_{\overline{\dbQ_l}}$ are connected over $\dbQ$ (in the sense of 3.3.1.2 A), then the images of $T_{v\overline{\dbQ_l}}^{\ad}$ in these two factors have the same dimension.
\finishproclaim

\Proclaim{3.3.2.1. Remark.} \rm
A natural question arises. Can we reobtain the results [Pi, 5.10-11] (see 4.2.14 below), using the above ideas of 3.3? We do not know its answer, as there is one obstruction. It consists in the fact that we do not know if (or when) we can assume that the Mumford--Tate group $T({}_2A^{\prime})$ of ${}_2A^{\prime}$ is (isomorphic to) $T_v^2$. We always have a natural monomorphism 
$$m_T:T_v^2\hookrightarrow T({}_2A^{\prime}).$$ 
If $m_T$ is an isomorphism and if the rank of $T_v^2$ is the same as the rank of the algebraic envelope of the $l$-adic representation attached to $A$, then we easily get (as $T({}_2A^{\prime})$ is a Mumford--Tate group) that [Pi, 5.10] holds for ${}_2A$. 

If ${}_2A_v$ is ordinary, then $m_T$ is an isomorphism (easy consequence --it is detailed in the proof of 4.4.8 b) below-- of the fact that we can assume ${}_2A^{\prime}$ is a canonical lift, cf. iii) of 3.3 E) and [Va2, a) of 4.4.1 2) and 4.6 P2]). So (based as well on Serre's results of [Chi, 3.8] and [Pi, 3.8]) we can reobtain [Pi, 5.10-11] for ${}_2A$ provided the ordinary reduction conjecture holds for ${}_2A$ (or equivalently, cf. [Va2, 4.9.23] for $A$). [Va2, a) of 4.4.1 2) and 4.6 P2] and iii) of 3.3 E) imply that the same can be said about $A$ itself.
\finishproclaim

\Proclaim{3.3.3. The shifting process.} \rm
We assume the reader is familiar with 4.0-1. 
Let $G_p$ (resp. $G_p(2)$) be the connected component of the origin of the $p$-adic representation (see 4.0) naturally attached to $A$ (resp. to ${}_2A$). The images of $G_p^{\der}$ and $G_p(2)^{\der}$ in $G^{\ad}_{\dbQ_p}=G_{2\dbQ_p}^{\ad}$ (cf. the identifications $ID$) are the same. This is a consequence of the fact that the Mumford--Tate group of $A\times_E {}_2A$ is a subgroup of $G_4$, cf. the choice of $m_4$. So to show that $G_p^{\der}=G^{\der}_{\dbQ_p}$ is the same as showing that $G_p(2)^{\der}=G_{2\dbQ_p}^{\der}$. In other words:

\medskip\noindent 
{\bf Fact.} {\it To prove the part of the Mumford--Tate conjecture involving derived groups for the pair $(A,p)$, is the same as proving it for the pair $({}_2A,p)$.}
 
\medskip\noindent
{\bf Definition.} An irreducible representation of a split, simple Lie algebra of classical Lie type $\scrT$ over a field $K$ of characteristic 0 is said to be an SD-standard representation, if it is associated to a minimal weight and if, in case $\scrT=A_n$, $n\in\dbN\setminus\{1,2\}$, it is of dimension $n+1$ (i.e. it is associated to one of the minimal weights $\overline{w}_1$ and $\overline{w}_n$). SD stands for Shimura--Deligne (cf. [De2, 2.3.10] and [Va1, 6.5-6]). Similarly, we speak about SD-standard representations of split, semisimple groups of classical Lie type over $K$, whose adjoints are simple. 

\medskip
We now list some of advantages (besides getting 3.3.2) we get by replacing $A$ by ${}_2A$.

\medskip\noindent
{\bf ADV.} {\bf 1)} {\it The irreducible subrepresentations of the representation of $\Lie(G_{2\dbC}^{\der})$ on $W_1\otimes_{\dbQ} \dbC$ factor through a simple factor of $\Lie(G_{2\dbC}^{\der})$ (i.e. we do not have to deal with tensor products of irreducible representations of simple factors of $\Lie(G_{2\dbC}^{\der})$) and are SD-standard representations of these simple factors. 

\smallskip
{\bf 2)} If $(G^{\ad},X^{\ad})$ has a simple factor $(G_0,X_0)$ of $D_n^{\dbR}$ type (resp. of $A_n$ type without involution), then we can assume that for each simple factor of $\Lie(G_{0\dbC}^{\der})$ both the half spin representations (resp. both the representations associated to the weights $\overline{w}_1$ and $\overline{w}_{n}$) of it are among the irreducible subrepresentations of $\Lie(G_{2\dbC}^{\ad})$ on $\scrW$, where $\scrW$ is a suitable irreducible subrepresentation of the representation of $Z^0(G_{2\dbC})$ on $\tilde W_1\otimes_{\dbQ} \dbC$.

\smallskip
{\bf 3)} We assume that all simple factors of $G^{\ad}_{\dbR}$ are non-compact and that none of the simple factors of $(G,X)$ are of some $A_n$ type with involution. Then we can assume $Z^0(G_2)=\dbG_m$.}

\medskip
{\bf Proof:} 1) is a consequence of 3.3.1.1 and of 3.1 D) and I'). 2) for the $D_n^{\dbR}$ type (resp. $A_n$ type without involution) is a consequence of 3.1 D) and I') (resp. of [De2, 2.3.10 and 2.3.13]; see [Va2, 6.5.1.1 and Case 1 of 6.6.5.1] for the $\dbZ_{(p)}$-context). 3) is a consequence of iv) of 3.1 I').

\medskip
The main disadvantages we get by replacing $A$ by ${}_2A$ are:

\medskip\noindent
{\bf DADV.} {\it We lose to a great extend the control on $Z(G)$. Moreover, often the dimension of ${}_2A$ is much bigger than the dimension of $A$.}  
\finishproclaim

\Proclaim{3.3.3.1. Remark.} \rm
Warning: $A$ and ${}_2A$ are by now means, in the general case, related through some isogeny. However, we can call them {\it adjoint-isogeneous}. Of course, instead of ${}_2A$ we get different variants of ADV by using any other abelian variety over $E$ which is adjoint-isogeneous to $A$ and which is obtained similarly to ${}_2A$ but via another injective map 
$$\tilde f_2\colon (\tilde G_2,\tilde X_2)\hookrightarrow (GSp(\tilde W_1,\tilde\psi_1),\tilde S_1),$$ 
with $(\tilde G_2^{\ad},\tilde X_2^{\ad})=(G^{\ad},X^{\ad})$ and such that some other desired properties hold. So for getting 3.3.3 ADV, we could have just quoted the proof of [De2, 2.3.10] and [De2, 2.3.13]; we worked 3.3 the way we did, to get (besides 3.3.3 ADV) 3.3.2 as well.

In practice, $\tilde f_2$ is a Hodge quasi product of injective maps indexed by simple factors of $(G^{\ad},X^{\ad})$. It is convenient to uniformize $\tilde f_2$ as follows. We write (as in 3.2.2 2)) $(G^{\ad},X^{\ad})=\prod_{i\in I} (G_i,X_i)$ as a product of simple, adjoint Shimura pairs. Let $F^i$ be as in 3.2.2 2). Let $F^{\text{com}}$ be the totally real number field generated by $F^i$'s, $i\in I$ (so $F^{\text{com}}$ is unramified above any prime over which all $F^i$'s are unramified). Let $F$ be a totally real number field containing $F^{\text{com}}$ and all $F_1$'s fields obtained as in 3.2 but for the different simple factors of $(G^{\ad},X^{\ad})$ which are not of $A_n$ type. Let 
$$(G_5,X_5):=\prod_{i\in I} (G_i^{F},X_i^{F}).$$
For each $i\in I$ we get injective maps 
$$(G_{i2},X_{i2})\hookrightarrow (G_{i3},X_{i3})\hookrightarrow (G_{i1},X_{i1})\operatornamewithlimits{\hookrightarrow}\limits^{\tilde f^i} (GSp(W_{i1},\psi_{i1}),S_{i1}),$$ 
with $(G_{i2}^{\ad},X_{i2}^{\ad})=(G_i,X_i)$ and $(G_{i3}^{\ad},X_{i3}^{\ad})=(G_i^{F},X_i^{F})$, cf. 3.2.2 6). Using Hodge quasi product, we ``put" them together. We get injective maps 
$$(\tilde G_2,\tilde X_2)\hookrightarrow (\tilde G_3,\tilde X_3)\hookrightarrow (\tilde G_1,\tilde X_1)\operatornamewithlimits{\hookrightarrow}\limits^{\tilde f_1} (GSp(\tilde W_1,\tilde\psi_1),\tilde S_1),$$
whose composite is (a uniformized) $\tilde f_2$. We have the following properties:

\medskip
{\bf 1)} $(\tilde G_3^{\ad},\tilde X_3^{\ad})=(G_5,X_5)$ and $(\tilde G_2^{\ad},\tilde X_2^{\ad})=(G^{\ad},X^{\ad})$;

\smallskip
{\bf 2)} the resulting injective map $(G^{\ad},X^{\ad})\hookrightarrow (G_5,X_5)$ is a product (indexed by elements of $I$) of the standard ones $(G_i,X_i)\hookrightarrow (G_i^F,X_i^F)$; 

\smallskip
{\bf 3)} $\tilde f_1$ is a PEL type embedding which is the Hodge quasi product of $\tilde f^i$'s and for which the analogue of c) of 3.2.2 2) holds;

\smallskip
{\bf 4)} the analogues of ADV 1) and 2) hold. 
\finishproclaim

3.3.2 (resp. 3.3.3) is used below just in Step 2 of 5.1.3 and so implicitly in 5.1.4 (resp. in 5.1.3-4, beginning with Step 2 of 5.1.3).

\Proclaim{3.3.4. Remark.} \rm
If 3.3.3 DADV looks unpleasant or if one would like to ``stick around" the injective map $f$ (and so $A$), then we can consider the PEL-envelope $(\tilde G,\tilde X)$ of $f$ as defined in [Va1, 4.3.12]. So we have natural injective maps
$$(G,X)\hookrightarrow (\tilde G,\tilde X)\operatornamewithlimits{\hookrightarrow}\limits^{\tilde f} (GSp(W,\psi),S)$$
whose composite is $f$, with $\tilde G$ as the connected component of the origin of the intersection of $GSp(W,\psi)$ with the centralizer in $GL(W)$ of the centralizer of $G$ in $GL(W)$. So $\tilde f$ is a PEL type embedding and we get a variant of 3.3.1 Fact' in its context: with $p$ big enough (so that the analogues of i) to iv) of 3.3 E) hold), we can view $Fr$ as a $\dbQ$--valued point of a torus $\tilde T^\prime$ of a suitable form $\tilde G^\prime$ of $\tilde G$; this form is automatically inner if $\tilde G^{\ad}$ has no simple factors of $D_s$ Lie type, with $s\in\dbN$, $s\Ge 4$ (cf. the connectedness aspects of [Ko, \S 7]). From the point of view of understanding the intersection of $G_p$ with $Z(G_{\dbQ_p})$, it is more convenient to work with $\tilde T^\prime$ instead of $T_1^\prime$. Warning: we do not know how one could get 3.3.2 via $\tilde T^\prime$ (i.e. by using $\tilde f$ instead of $f_1$). 

We now consider the centralizer $\tilde C$ of $Z^0(G)$ in $GSp(W,\psi)$. Let $\tilde G_1^{\ad}$ be the maximal factor of $\tilde C^{\ad}$ with the property that each simple factor of it is non-compact over $\dbR$. Let $\tilde G_1$ be the maximal normal, reductive subgroup of $\tilde C$ whose adjoint is $\tilde G_1^{\ad}$. We have the following variant of the above paragraph: instead of $\tilde G$ we work with $\tilde G_1$. The inclusion $\tilde G_1\hookrightarrow GSp(W,\psi)$ extends uniquely to an injective map
$$\tilde f_1:(\tilde G_1,\tilde X_1)\hookrightarrow (GSp(W,\psi),S)$$ 
through which $f$ factors. $\tilde f_1$ is a PEL type embedding and so, for $p$ big enough, we can view $Fr$ as a $\dbQ$--valued point of a torus $\tilde T_1^\prime$ of a suitable form of $\tilde G_1$. What we gain in this way: $\tilde C^{\ad}$ (and so also $\tilde G_1^{\ad}$) has no factors of $D_s$ Lie type and so this form is inner.  
\finishproclaim

\Proclaim{3.3.5. Remark.} \rm
We come back to 3.3 B). $f_1$ is a Hodge quasi product indexed naturally by the set $I$ of 3.3.3.1, cf. d) of 3.2.2 2). So if we choose $L_1$ to be well adapted to it, then ${}_2A$ is a product (indexed by $I$) of abelian varieties whose Mumford--Tate groups have simple adjoints.
\finishproclaim

\bigskip
\noindent
{\boldsectionfont \S4. The basic techniques}
\bigskip

\Proclaim{4.0.} \rm
Let $E$ be a number field. Let $A$ be an abelian variety over $E$. Let $p$ be a rational prime. Let 
$$H^1:=H_{\text{\'et}}^1(A_{\Ebar},\dbZ_p).$$ 
Let $G_{\dbQ_p}$ be the algebraic group over $\dbQ_p$ defined (see [Se]) as the connected component of the origin of the algebraic envelope of the natural $p$-adic representation
$$\rho\colon\Gal(E)\to GL(H_1[{1\over p}])(\dbQ_p).$$
So $G_{\dbQ_p}$ is a connected subgroup of $GL(H^1\fracwithdelims[]1p)$ and there is a finite field extension $E_1$ of $E$ such that $\rho(\Gal(E_1))$ is a compact, open subgroup of $G_{\dbQ_p}(\dbQ_p)$.

Let $H_A$ be the Mumford--Tate group of $A$. We recall that $H_A$ is a reductive group over $\dbQ$ and that any embedding $E\hookrightarrow\dbC$ allows the interpretation of $H_A$ as the smallest subgroup of $GL(H_B^1(A_{\dbC},\dbQ))$ with the property that the Hodge cocharacter 
$$\mu(\dbC)\colon\dbG_m\to GL(H_B^1(A_{\dbC},\dbQ)\otimes_{\dbQ} \dbC)$$ 
factors through $H_A$ (cf. [De3, 2.11 and 3.4]). We recall: if $H_B^1(A_{\dbC},\dbQ)\otimes_{\dbQ} \dbC=F^{1,0}\oplus F^{0,1}$ is the usual Hodge decomposition, then $\beta\in\dbG_m(\dbC)$ acts through $\mu(\dbC)$ as the multiplication with $\beta^{-1}$ on $F^{1,0}$ and as identity on $F^{0,1}$. From now on we fix (arbitrarily) an embedding $i_E\colon E\hookrightarrow\dbC$. $H^1_B:=H^1_B(A_{\dbC},\dbZ)$ is obtained via it. The smallest property of $H_A$ is specific to Mumford--Tate groups; it is used freely in what follows. Let $T^0:=Z^0(H_A)$. 

The connection between the Betti cohomology and the \'etale cohomology with coefficients in $\dbQ_p$ of $A_{\dbC}$ (see [SGA4, Exp. XI]), identifies canonically $H^1=H^1_B\otimes_{dbZ_{(p)}} \dbZ_p$ and so we view $H_{A\dbQ_p}$ as a subgroup of $GL(H^1\fracwithdelims[]1p)$ (it does not matter which embedding $i_E$ we fixed). We have the following remarkable prediction (cf. [Mu1] and [Se]):
\finishproclaim

\Proclaim{4.1. The Mumford--Tate conjecture for the pair $(A,p)$.}
As subgroups of $GL(H^1\fracwithdelims[]1p)$, $G_{\dbQ_p}=H_{A\dbQ_p}$.
\finishproclaim

\Proclaim{4.1.1. Goal.}
The main goal in \S4-5 is to approach (and prove in many cases) this conjecture, via the strategy of 1.2.
\finishproclaim

\Proclaim{4.1.2. A variant of 1.2.} \rm
There is a simpler variant of 1.2, where we assume that both $G_{\dbQ_p}$ and $H_{A\dbQ_p}$ are split, reductive groups. References to this variant are made below in: 4.2.8.1, 4.2.9-10, 4.3.3.1, 4.3.6.1.1 and Step 3 of 5.1.3.
\finishproclaim

\Proclaim{4.2. Previously used tools.} \rm
In 4.2 we first review some well known facts about 4.1 and then we show how it is implied by these facts and by (the assumed validity of) 1.2 or (of) its variant 4.1.2. $G_{\dbQ_p}$ commutes (cf. 4.2.3) with $Z(H_{A\dbQ_p})$. So let $\tilde G_{\dbQ_p}$ be the reductive subgroup of $H_{A\dbQ_p}$ generated by $T^0_{\dbQ_p}$ and $G_{\dbQ_p}$. We have $G_{\dbQ_p}^{\der}=\tilde G_{\dbQ_p}^{\der}$.
For any rational prime $q$, $G_{\dbQ_q}$ (resp. $\tilde G_{\dbQ_q}$) has the same significance as $G_{\dbQ_p}$ (resp. as $\tilde G_{\dbQ_p}$) but working in the $\dbQ_q$-\'etale context. We view $G_{\dbQ_q}$ and $\tilde G_{\dbQ_q}$ as subgroups of $H_{A\dbQ_q}$ and so of $GL(H_B^1\otimes_{\dbZ} \dbQ_q)$.
\finishproclaim

\Proclaim{4.2.1. Faltings' theorem.} 
$G_{\dbQ_p}$ is a reductive group (cf. [Fa1]).
\finishproclaim

\Proclaim{4.2.2. Pyatetskii-Shapiro--Deligne--Borovoi's theorem.} 
As subgroups of $GL(H^1\fracwithdelims[]1p)$, $G_{\dbQ_p}$ is a subgroup of $H_{A\dbQ_p}$ (for instance, cf. [De3, 2.9 and 2.11]).
\finishproclaim

\Proclaim{4.2.3. Faltings' theorem.}
We have
$$\End(A_{\Ebar})^{\text{op}}\otimes_{\dbZ_p} \dbQ_p\arrowsim\End(H^1\otimes_{\dbZ_p} \dbQ_p)^{G_{\dbQ_p}(\dbQ_p)}=\End(H^1\otimes_{\dbZ_p} \dbQ_p)^{H_{A\dbQ_p}(\dbQ_p)},$$
where the first isomorphism is defined by the action of endomorphisms of $A_{\Ebar}$ on $H^1$ (cf. [Fa1]; here $\text{op}$ stands for opposite).
\finishproclaim

\Proclaim{4.2.4. Serre's theorem.}
The rank $r_p$ of $G_{\dbQ_p}$ does not depend on $p$ (cf. [Chi, th. (3.10)]).
\finishproclaim

\Proclaim{4.2.5. Theorem.} 
If the Mumford--Tate conjecture for the abelian variety $A$ is true for a prime $p$, then it is true for all rational primes (cf. [LP2, 4.3] or [Ta1, 3.4]).
\finishproclaim

The proof of this result (cf. [Ta1, 3.4]) relies on 4.2.1-4 and on: 

\Proclaim{4.2.5.1. Zarhin's lemma.}
Let $\scrV$ be a finite dimensional vector space over a field of characteristic 0. Let $\grg_1\subset\grg_2\subset\End(\scrV)$ be Lie algebras of reductive subgroups of $GL(\scrV)$. We assume that the centralizer of $\grg_i$ in $\End(\scrV)$ is independent on $i\in\{1,2\}$ and that the rank of $\grg_1$ is equal to the rank of $\grg_2$. Then $\grg_1=\grg_2$ (cf. [Za, key lemma of \S5]). 
\finishproclaim

\Proclaim{4.2.6. CM type case.} \rm
The Mumford--Tate conjecture is true for abelian varieties (over number fields) having complex multiplication over $\overline{\dbQ}$ (this is well known; it is an easy exercise which can be solved using 4.2.7.1 below, 4.2.4 and the main result of [Wi]).
\finishproclaim

\Proclaim{4.2.7. Frobenius tori.} 
There is a number field $E^T$ and a torus $T$ of $H_{AE^T}$ (called a Frobenius torus of $A$) such that for any finite prime $v$ of $E^T$, $T_{E_v^T}$ is isomorphic to a maximal torus of $G_{E_v^T}$.
\finishproclaim

This follows from [Chi, ch. 3]: taking $\dbQ$ as the number field we get a $\dbQ$--torus such that the isomorphism part of 4.2.7 holds for $v$ provided the characteristic of $k(v)$ is big enough; by passing to a suitable finite field extension $E^T$ of $\dbQ$, we can control (cf. 4.2.4) the left aside primes as well as we can assume we have an $E^T$-torus which is embeddable in $H_{AE^T}$.

\Proclaim{4.2.7.1. Standard constructions.} \rm
For the constructions briefly reviewed here, we refer to [Chi] and [Pi]. Let $q$ be a rational prime relatively prime to $p$. We assume there is a prime $v_q$ of $E$ dividing $q$ and such $A$ has good reduction $A_{v_q}$ w.r.t. $v_q$. Let $\overline{v}_q$ be any prime of $\Ebar$ dividing $v_q$. Let $F(\overline{v}_q)\in\Gal(E)$ be a Frobenius automorphism of $\overline{v}_q$. Since the restriction of $\rho$ to the inertia group of $\overline{v}_q$ is trivial, 
$$\rho(F(\overline{v}_q))\in GL(H^1[{1\over p}])(\dbQ_p)$$ 
is well defined regardless of the choice of $F(\overline{v}_q)$; it is referred as a Frobenius element. Let $T_{v_q}$ be the $\dbQ$--model (obtained as in [Chi, 3.a]) of the connected component of the origin of the Zariski closure (in $GL(H^1[{1\over p}])$) of the set $\{\rho(F(\overline{v}_q)^n|n\in\dbZ\}$ of $\dbQ_p$-valued points of $GL(H^1[{1\over p}])$. So $T_{v_q}$ is a $\dbQ$--torus and $T_{v_q\dbQ_p}$ is a torus of $G_{\dbQ_p}$. 

$T_{v_q}$ depends only on (the Frobenius endomorphism $\pi(A_{v_q})$ of) $A_{v_q}$ (see [Chi, 3.a and 3.b]); moreover, in this statement we can replace $A_{v_q}$ by any abelian variety $A^\prime_{v_q}$ over a finite field extension $k_1(v_q)$ of $k(v_q)$ which is isogeneous to $A_{v_qk_1(v_q)}$. 

We consider now an abelian scheme $B^\prime$ over a complete, characteristic 0 DVR $R$ having $k_1(v_q)$ as its residue field, which lifts $A^\prime_{v_q}$ and has complex multiplication. Let $MT(B^\prime)$ be the Mumford--Tate group of $B^\prime$; as in 4.0, it is irrelevant which monomorphism $R\hookrightarrow\dbC$ we use to define it. Let $DC$ be the centralizer (in $GL(H_B^1(B^\prime_{\dbC},\dbQ))$) of the centralizer $C$ of $MT(B^\prime)$ in $GL(H_B^1(B^\prime_{\dbC},\dbQ))$. As any abelian variety over $k_1(v_q)$ has complex multiplication and as $\Lie(C)$ is naturally identified with the Lie algebra of (Betti realizations) of $\dbQ$--endomorphisms of $B^\prime_{\dbC}$, a sufficiently high, positive, integral power of any invertible element of the $\dbQ$--subalgebra of $\End(A^\prime_{v_q})$ generated by the Frobenius endomorphism of $A^\prime_{v_q}$, is naturally identified with a $\dbQ$--valued point of $DC$; so (cf. loc. cit.) $T_{v_q}$ can be identified with a subtorus of $DC$. If $B^\prime$ is equipped with a polarization $p_{B^\prime}$, then above we can replace $DC$ by the connected component of the origin of $DC\cap GSp(H_B^1(B^\prime_{\dbC},\dbQ),p_{B^\prime})$. So, if $\scrB$ is a semisimple $\dbQ$--algebra of $\End(H_B^1(B^\prime_{\dbC},\dbQ))$ formed by elements fixed by $MT(B^\prime)$, then $T_{v_q}$ is a torus of the intersection of $GSp(H_B^1(B^\prime_{\dbC},\dbQ),p_{B^\prime})$ with the subgroup of $GL(H_B^1(B^\prime_{\dbC},\dbQ))$ fixing all elements of $\scrB$.

It is well known (see [Chi, 3.b] and [Pi, 3.6-7]) that for an infinite set of primes $v_q$ of $E$ of good reduction for $A$ and relatively prime to $p$, $T_{v_q\dbQ_p}$ is a maximal torus of $G_{\dbQ_p}$. 
\finishproclaim

\Proclaim{4.2.7.2. Fact.}
There is a subtorus $T^{00}$ of $T^0$ which does not depend on $p$ and with the property that $G_{\dbQ_q}$ is generated by $G^{\der}_{\dbQ_p}$ and by $T^{00}_{\dbQ_p}$.
\finishproclaim

\proof
We assume $T_{v_q\dbQ_p}$ is a maximal torus of $G_{\dbQ_p}$. Let $C_A$ be the centralizer of $Z(H_A)$ in $GL(H^1_B(A_{\dbC},\dbQ))$. $Z(C_A)$ is a torus. Its Lie algebra, when viewed as a $\dbQ$--Lie subalgebra of $\End(H^1[{1\over p}])$, is generated by $\dbQ$--linear combinations of $\dbQ_p$-adic realizations of endomorphisms of $A_{v_q}$. Here we identify via $\overline{v}_q$, $H^1$ with $H^1_{\text{\'et}}(A_{v_q\overline{k(v)}},\dbZ_p)$. Let $C_{A_{v_q}}$ be the centralizer of $\Lie(Z(C_A))$ in the group of invertible elements of $\End(A_{v_q})\otimes_{\dbZ}\dbQ$. We identify $C_{A_{v_q}\dbQ_p}$ with a subgroup of $GL(H^1[{1\over p}])$.

We still denote by $\pi(A_{v_q})$ the $\dbQ_p$-adic realization $\rho(F(\overline{v}_q))$ of $\pi(A_{v_q})$. After replacing it by a positive, integral power of it dividing the order of $Z(C_A^{\der})$, as an element of $GL(H^1[{1\over p}])$ we can write it as a product
$$\pi(A_{v_q})=\pi_{C}\pi_{C^\bot},$$ 
with $\pi_{C}\in Z^0(C_A)(\dbQ_p)$ and with $\pi_{C^\bot}\in C_A^{\der}(\dbQ_p)$.
As $\pi(A_{v_q})\in H_A(\dbQ_p)$, after a possible similar replacement we get that $\pi_{C^\bot}\in H_A^{\der}(\dbQ_p)$ and that $\pi_{C}\in T^0(\dbQ_p)$. Using $C_{A_{v_q}}$ instead of $C_A$, we get that in fact $\pi_{C}\in T^0(\dbQ)$.

We take $T^{00}$ to be the connected component of the origin of the smallest subgroup of $T^0$ having $\pi_{C}$ as a $\dbQ$--valued point. The images of $G_{\dbQ_p}$ and $T_{v_q\dbQ_p}$ in $H_{A\dbQ_p}^{\ab}$ are the same torus $T^{01}$. As the images of $\pi(A_{v_q})$ and of $\pi_C$ in $H_{A\dbQ_p}^{\ab}(\dbQ_p)$ are the same, $T^{01}$ is the extension to $\dbQ_p$ of the image of $T^{00}$ in $T^{01}$. As $Z^0(G_{\dbQ_p})$ is contained in $T^0_{\dbQ_p}$ and as $G^{\der}_{\dbQ_p}$ is the connected component of the origin of $G_{\dbQ_p}\cap H_{A\dbQ_p}^{\der}$, $G_{\dbQ_p}$ is generated by $G^{\der}_{\dbQ_p}$ and $T^{00}_{\dbQ_p}$. The fact that $T^{00}$ does not depend on $v_q$ is checked in the standard way by considering another prime $v_l$ of $E$ dividing a rational prime different from $q$ and such that the Frobenius torus $T_{v_l}$ has the same rank as $T_{v_q}$. This ends the proof.

\Proclaim{4.2.7.2.1. Exercise.} \rm
Use the variant of 3.3.4 to give a second proof of 4.2.7.2. 
\finishproclaim

\Proclaim{4.2.8. Good primes.} \rm
The set of primes $p$ such that $G_{\dbQ_p}$ is a split, reductive group is of positive Dirichlet density and so is infinite. This is a direct consequence of 4.2.7 ($T$ splits over a finite field extension of $E^T$; the set of rational primes which split completely in this extension is --for instance, see [La, p. 168]-- of positive Dirichlet density).
\finishproclaim

\Proclaim{4.2.8.1. Variant.} \rm
We have a variant of 4.2.8: the set of primes $p$ such that $G_{\dbQ_p}$ and $H_{A\dbQ_p}$ are both split, reductive groups is of positive Dirichlet density and so is infinite (in the above argument we need to replace $E^T$ by its composite with a number field over which an arbitrarily chosen torus of $H_A$ splits).
\finishproclaim

\Proclaim{4.2.9. Corollary.} 
The conjecture 4.1 is implied by 1.2 and so by 4.1.2.
\finishproclaim

\proof
This is a direct consequence of 4.2.7-8, via 4.2.5.

\Proclaim{4.2.10. Products.} \rm
The conjecture 4.1 (or the expectation 1.2 or its variant 4.1.2) is true for $A$ iff it is true for $A^n\times (A^t)^m$, with $n$ and $m$ non-negative integers such that $n+m>0$.
\finishproclaim

\Proclaim{4.2.11. Isogenies.} \rm
The conjecture 4.1 is true for $A$ iff it is true for an isogeny $A^{\prime}$ of the extension of $A$ to a finite field extension $E^{\prime}$ of $E$.
\finishproclaim

The statements 4.2.10-11 are obvious.

\Proclaim{4.2.12. Characterizations of the CM type cases.} \rm
The following assertions are equivalent:

\medskip
{\bf a)} $A_{\Ebar}$ has complex multiplication.

{\bf b)} $H_A$ is a torus.

{\bf c)} For a (any) prime $p$, $G_{\dbQ_p}$ is a torus.

\medskip
The equivalence of a) (resp. c)) and b) follows from  definitions (resp. from 4.2.3).
\finishproclaim 

\Proclaim{4.2.13. On $\tilde G_{\dbQ_p}$.} \rm
We have $T^0_{\dbQ_p}=Z^0(\tilde G_{\dbQ_p})$, cf. 4.2.3. So if $G_{\dbQ_p}$ and $H_{A\dbQ_p}$ are split, then $\tilde G_{\dbQ_p}$ is also so. So we do not refer to the variant of 1.2 in which $G_{\dbQ_p}$ is replaced by $\tilde G_{\dbQ_p}$: we prefer to keep the original form of 1.2 (while simultaneously working with its simplest variant 4.1.2).   
\finishproclaim

\Proclaim{4.2.14. Pink's theorem.}
All non-trivial, simple factors of $G^{\ad}_{\overline{\dbQ_p}}$ are of $A_l$, $B_l$, $C_l$ or $D_l$ Lie type, $l\in\dbN$; moreover, the non-trivial, irreducible subrepresentations of the natural representation of the Lie algebra of any such simple factor on $H^1\otimes_{\dbZ_p} \overline{\dbQ_p}$, are associated to minimal weights. (cf. [Pi, 5.11]; see also [Ta2] for a slightly less general result).
\finishproclaim

This theorem is obtained as a corollary to (slight reformulation):

\Proclaim{4.2.14.1. Pink's theorem.} 
$G_{\overline{\dbQ_p}}$ is generated by cocharacters having the property that the resulting representations of $\dbG_m$ on $H^1\otimes_{\dbZ_p} \overline{\dbQ_p}$ are associated to the trivial and the inverse of the identity character of $\dbG_m$, the multiplicities being the same (cf. [Pi, 5.10]).
\finishproclaim

The proof of 4.2.14.1 (cf. [Pi, 3.15]) is obtained by combining the following three things.

\medskip
{\bf a)} Serre's theory (see [Pi, 3.5-8]; part of it was recalled in 4.2.4 and 4.2.7).

\smallskip
{\bf b)} Results (see [Pi, 2.3]) on the relation between Hodge cocharacters and Newton quasi-cocharacters (as defined and used in [Pi]). A variant of [Pi, 2.3] is contained in [RR, 4.2]; for a completely new proof of loc. cit. see [Va2, th. 1 of Appendix]. 

\smallskip
{\bf c)} Katz--Messing' result recalled in [Pi, 3.10].

\medskip
For two other approaches to b) (one particular and one general) in the context of abelian varieties, see 3.3.2.1 and 4.4.10 below.  
 
\Proclaim{4.2.15. Bogomolov's theorem.}
$G_{\dbQ_p}$ contains $Z(GL(H^1\fracwithdelims[]1p))$ (cf. [Bog]).
\finishproclaim

\Proclaim{4.2.16. References.} \rm
For cases the conjecture 4.1 was previously known to be true, see [Pi, the paragraph after the proof of 5.8] and [Pi, 5.13-15]. See also 5.4 6) below.
\finishproclaim

\Proclaim{4.2.17. Overview of what follows.} \rm
In 4.3 we start the proof of 1.2 in the general context. In 5.1.1-4 we accomplish its end in many situations. Some basic (new) techniques are presented in 4.3-4, while \S5 gathers many conclusions of \S3-4.

\smallskip
\Proclaim{4.3. Nice arrangements.} \rm
\finishproclaim

\Proclaim{4.3.0. On $N(A)$.} \rm
Let $p$ be a prime such that $G_{\dbQ_p}$ is a split, reductive group. 
The natural numbers $N_1(A)$, $N_2(A)$ and $N_3(A)$ are respectively defined below in 4.3.1, 4.3.4 and 5.1.0. As their notation points out, they depend only on $A$. We take
$$N(A):=\max\{N_1(A),N_2(A)\}.$$

The condition $p\Ge N_1(A)$ guarantees that $E$ is unramified above $p$ and $A$ has good reduction w.r.t. primes of $E$ dividing $p$. The condition $p\Ge N_2(A)$ says that we can obtain out of the faithful representation $G_{\dbQ_p}\hookrightarrow GL(H^1\fracwithdelims[]1p)$ ``whatever we desire"; i.e. it guarantees that the mentioned representation is ``strongly unramified", in some sense which will not be made explicit here. The condition $q\Ge N_3(A)$ says (cf. 3.3 F) and 5.1.0 below) that in some context we get a SHS enjoying nice properties. We are not at all interested here in a tight bound for $N(A)$ (cf. 4.2.5); so we treat $N_2(A)$ and $N_3(A)$ in a loose way, i.e. we do not stop here to get tight bounds for them. However, we do think that it is of some theoretical interest to have tight bounds for $N_2(A)$ (or similarly constructed numbers); accordingly we present the ideas in such a way, that in any concrete situation, the interested reader can ``trace out" such tight bounds.

We deal with the problem: if $p\Ge N(A)$ show that $G_{\dbQ_p}=H_{A\dbQ_p}$. An afterthought of the approach that follows towards the solution of this problem is included in 5.5.2.
Those who are interested just in this approach (to the proof 1.2), can look just at 4.3.1-3, at 4.3.4 but just working with a $\dbZ_p$-very well positioned family of tensors (so at 4.3.4.1, at ii), iv) and vi) of 4.3.4.2, at 4.3.4.6 with $\max\{d(A),4\}$ being replaced by 4, and at 4.3.4.8 a)), at 4.3.5-8, and at the whole of 4.4.

\Proclaim{4.3.1. $N_1(A)$.} \rm
There is a smallest integer $N_1(A)\Ge 3$ such that any prime $v$ of $E$ dividing a rational prime $l\Ge N_1(A)$, is unramified over $l$ and is a prime of good reduction for $A$. Obviously $N_1(A)$ depends only on $A$ (and on $E$). From now on we assume that $p\Ge N_1(A)$. 
\finishproclaim

\Proclaim{4.3.2. Principal polarizations.} \rm
We can assume $A$ has a principal polarization $p_A$. First argument: the Zarhin's trick (see [M-B, p. 205]) says $B:=A^4\times (A^t)^4$ has a principal polarization and so we can replace (cf. 4.2.10) if needed $A$ by $B$. Second argument: cf. [Mu2, cor. 1, p. 234] and 4.2.11.

This second approach is more convenient: we can assume that $A$ is a product of absolutely simple abelian varieties and that $p_A$ is a product of principal polarizations of these simple factors. We do not need this assumption in what follows: we include it just for the sake of future references. So let 
$$(A,p_A)=\prod_{i\in I_f} (A_i,p_{A_i}),$$
with $I_f$ a finite set, with $A_i$ an absolutely simple abelian variety, and with the pair $(A_i,p_{A_i})$ a principally polarized abelian variety over $E$, $\forall i\in I_f$. At the level of cohomology, we get a direct sum decomposition
$$H^1=\oplus_{i\in I_f} H^1(i),$$
with $H^1(i):=H^1_{\text{\'et}}(A_{i\Ebar},\dbZ_p)$.

For $i\in I_f$, let $H_A(i)$ be the Mumford--Tate group of $A_i$. We have:  

\Proclaim{Fact.} 
$H_A(i)^{\ad}$ is either a simple $\dbQ$--group or the $\dbQ$--rank of each simple factor of it is 0. 
\finishproclaim

In the last case, the Shimura pair $(H_A(i),X_A(i))$ attached to $A_i$ is of compact type, cf. [BHC]. The argument of this well known fact is simple and short and so we include it. The irreducible representations of a product of semisimple Lie algebras over $\dbC$ are obtained by taking tensor products of irreducible representations of its simple factors. But in the case of the representation of $\Lie(H_A(i)_{\dbC}^{\der})$ on $H^1_B\otimes_{\dbZ} \dbC$ we can not have such tensors products involving two non-compact factors of $\Lie(H_A(i)_{\dbR}^{\der})$ (we recall that the simple factors of $H_A(i)_{\dbR}^{\ad}$ are absolutely simple, cf. [De2, 1.2.1]). 

We can assume $(H_A(i)^{\ad},X_A(i)^{\ad})$ has a simple factor $(H(i_0),X(i_0))$ whose extension to $\dbR$ does not have compact factors; for instance, this happens if we have a simple factor of $H_A(i)^{\ad}$ whose $\dbQ$--rank is positive. The $\dbQ$--vector subspace $W(i)$ of $H^1_B(A_{i\dbC},\dbQ)$ generated by non-trivial, irreducible subrepresentations of $\Lie(H(i_0))$ is $H_A(i)$-invariant and fixed by the simple factors of $\Lie(H_A(i)^{\der})$ different from $\Lie(H(i_0))$; here we view $\Lie(H(i_0))$ naturally as a Lie subalgebra of $\Lie(H_A(i)^{\der})$ and so of $\grg\grl(W(i))$. 
So from Weyl's complete reducibility theorem, as $A_i$ is an absolutely simple abelian variety, we get $W(i)=H^1_B(A_{i\dbC},\dbQ)$. So $H_A(i)^{\ad}=H(i_0)$. This argues the Fact.

\Proclaim{4.3.2.1. Some adjoint Mumford--Tate groups.} \rm
We view $H_A(i)$ as a quotient of $H_A$. We get naturally a monomorphism
$$H_A^{\ad}\hookrightarrow\prod_{i\in I_f} H_A(i)^{\ad}.$$ 
Warning: $H_A^{\ad}$ is isomorphic to a product of simple factors of this product but it is not necessarily equal to the whole product; if it is not equal to the whole product, then two of adjoint Shimura pairs $(H_A(i)^{\ad},X_A(i)^{\ad})$, $i\in I_f$, have isomorphic simple factors. 
\finishproclaim

$p_A$ induces a perfect form (still denoted by $p_A$) 
$$p_A\colon H^1\otimes_{\dbZ_p} H^1\to\dbZ_p(-1)$$ 
(fixed by $\Gal(E)$). We still denote by $(H^1,p_A)$ the symplectic space over $\dbZ_p$ we get.

\Proclaim{4.3.3. Lemma.}
There is a $\dbZ_p$-lattice $H^1_{\prime}$ of $H^1\fracwithdelims[]1p$ such that we still get a perfect form $p_A\colon H^1_{\prime}\otimes_{\dbZ_p} H_{\prime}^1\to\dbZ_p$ and moreover the Zariski closure of $G_{\dbQ_p}$ in $GL(H_{\prime}^1)$ is a (split) reductive group $G_{\dbZ_p}$ over $\dbZ_p$. 
\finishproclaim

\proof
We get this by combining the following three facts: 

\medskip
{\bf a)} the inclusion $G_{\dbQ_p}\hookrightarrow G\Sp(H^1\fracwithdelims[]1p,p_A)$ of reductive groups over $\dbQ_p$ can be extended to an inclusion $G_{\dbZ_p}\hookrightarrow G\Sp_{\dbZ_p}$ of reductive groups over $\dbZ_p$;

\smallskip 
{\bf b)} for any $\dbZ_p$-reductive form $G\Sp_{\dbZ_p}$ of $G\Sp(H^1\fracwithdelims[]1p,p_A)$ there is a $\dbZ_p$-lattice $H^1_\prime$ of $H^1[{1\over p}]$ and there is $\beta\in\dbG_m(\dbQ_p)$ such that we get a perfect form ${\beta}p_A\colon H_{\prime}^1\otimes_{\dbZ_p} H_{\prime}^1\to\dbZ_p$ and $G\Sp_{\dbZ_p}=G\Sp(H_{\prime}^1,{\beta}p_A)$ (so $G\Sp_{\dbZ_p}(\dbZ_p)=G\Sp(H_{\prime}^1,{\beta}p_A)(H_{\prime}^1)$) , cf. [Ti, end of 2.5];

\smallskip 
{\bf c)} there is a cocharacter $\mu_p\colon\dbG_m\to G_{\dbQ_p}$ producing a direct sum decomposition $H^1\fracwithdelims[]1p=F^1_p\oplus F^0_p$, with $\mu_p$ acting on $F^0_p$ trivially and on $F^1_p$ as the inverse of the identity character of $\dbG_m$; it allows us to get ``rid of" $\beta$'s of b).  
\medskip

Such a cocharacter exists over $\overline{\dbQ_p}$ (cf. 4.2.14.1). As $G_{\dbQ_p}$ is split, such a cocharacter exists over $\dbQ_p$. This argues c); a) is obvious as the reductive groups are split. This ends the proof.

\smallskip
Warning: we do not require $H_{\prime}^1$ to be normalized by $\rho(\Gal(E))$.
\finishproclaim

\Proclaim{4.3.3.0. Remark.} \rm
We can assume that $H_{\prime}^1=\oplus_{i\in I_f} H_{\prime}^1\cap H^1(i)\fracwithdelims[]1p$. Argument: we can perform 4.3.3 for each simple factor $A_i$ of $A$, $i\in I_f$.
\finishproclaim

\Proclaim{4.3.3.1. Variants.} \rm
We have variants of 4.3.3 and 4.3.3.0: if $H_{A\dbQ_p}$ is split, then we can assume as well that its closure in $G\Sp(H_{\prime}^1,p_A)$ is a split, reductive group over $\dbZ_p$. To see this for 4.3.3, in 4.3.3 a) we have to start with inclusions of split, reductive groups $G_{\dbZ_p}\hookrightarrow H_{A\dbZ_p}\hookrightarrow G\Sp_{\dbZ_p}$; the rest is the same.  As $H_{A\dbQ_p}$ is split iff $H_A(i)_{\dbQ_p}$ is split $\forall i\in I_f$ (cf. 4.3.2.1), using entirely the same method we get the variant of 4.3.3.0.
\finishproclaim

\Proclaim{4.3.4. Proposition.} 
There is $N_2(A)\in\dbN$ such that if also $p\Ge N_2(A)$, then there is a family $(v_{\alpha})_{\alpha\in\scrJ}$ of tensors of $\scrT(H_{\prime}^1)$ which is $\dbZ_p$-very well positioned for the group $G_{\dbQ_p}$.
\finishproclaim

Most of 4.3.4.1-8 is dedicated to the proof of this Proposition: its proof ends with the proof of 4.3.4.8 a) below.

\Proclaim{4.3.4.1. Finiteness properties.} \rm
In what follows $\scrK$ is an arbitrary field of characteristic $0$. It is well known (for instance, see [Hu1]) that the finite numbers of a) and b) below do not depend on the choice of $\scrK$; we need $\scrK$ just to simplify the presentation. We have:

\medskip
\item{{\bf a)}}
For any $d\in\dbN$, there is only a finite number of isomorphism classes of split, semisimple groups over $\scrK$ of dimension not bigger than $d$.

\smallskip
\item{{\bf b)}} 
For any split, semisimple group over $\scrK$ there is only a finite number of isomorphism classes of representations of it of dimension smaller or equal to a given $r\in\dbN$.
\finishproclaim

For a), cf. the classification of split, simply connected semisimple groups over $\scrK$ and the finiteness of their centers. For b), cf. the classification of irreducible representations --all of them are absolutely irreducible-- of a split, semisimple Lie algebra over $\scrK$, the Weyl's complete reducibility theorem, and the dimension formula --in the case of a split, simple Lie algebra over $\scrK$-- of these representations. For all these facts see [Hu1]. [Hu1] is stated over $\overline{\scrK}$; as we are dealing with split, semisimple groups, the same remains true for any such field $\scrK$ (cf. also the uniqueness theorem of [SGA3, Vol. III, p. 313-314]). 

\Proclaim{4.3.4.2. Some invariants and $N_2(A)$.} \rm
Let  $g_A:=\dim_E(A)$. Let $d_A:=4g_A^2$ and $r_A:=2g_A$. In our case we can work most of the time with considerably smaller expressions of $d_A$ and $r_A$; for instance, the dimension of $G_{\dbQ_p}$ is smaller or equal to $\dim_{\dbQ}(H_A)$. But we recall (cf. 4.3.0) that we are not interested here for the best bound of $N_2(A)$. We need $N_2(A)$ to be greater than:

\medskip
\item{{\bf i)}}
all orders of centers of semisimple groups over $\scrK$ of dimension not bigger than $d_A$;

\smallskip
\item{{\bf ii)}}
all numbers $B(G)$, with $G$ running through all simple, adjoint groups over $\scrK$ of dimension not bigger than $d_A$; if $G$ is (resp. is not) of classical Lie type, then $B(G)$ is defined in [Va1, 5.7.2] (resp. we define for convenience $B(G):=30$; to be compared with [Hu2, p. 48-49] and [Va2, 2.2.23 B]);

\smallskip
\item{{\bf iii)}}
$\max\{2d_A,r_A\}$;

\smallskip
\item{{\bf iv)}}
the numbers $s(\grg,W)$ (defined in [Va1, 4.3.3]), for any split, semisimple Lie algebra $\grg$ over $\scrK$ of dimension not bigger than $d_A$, and for every representation $W$ of $\grg$ over $\scrK$ of dimension not bigger than $r_A$;

\smallskip
\item{{\bf v)}}
the numbers $k(\grg,W_0)$ (defined in 4.3.4.2.3 below), for any $\grg$ as in iv), and for every non-trivial, simple $\grg$-module $W_0$ of dimension not bigger than $r_A$;

\smallskip
\item{{\bf vi)}}
the natural numbers $2+d(A)$ and $\Pi(A)$ defined below in 4.3.4.2.1-2.
\finishproclaim

\Proclaim{4.3.4.2.1. The $d(A)$ invariant.} \rm
$d(A)\in\dbN$ is the smallest number such that for any prime $q$, $\tilde G_{\dbQ_q}$ is the subgroup of $GL(H^1_B\otimes_{\dbZ} \dbQ_q)$ fixing some homogeneous tensors of $\scrT(H^1_B\otimes_{\dbZ} \dbQ_q)$ of degree not bigger than $d(A)$.
\finishproclaim

\Proclaim{4.3.4.2.1.1. Argument.} \rm
To see why $d(A)$ exists, we first remark that $\tilde G_{\dbQ_q}$ is the subgroup of $GL(H^1_B\otimes_{\dbZ} \dbQ_q)$ fixing some homogeneous tensors of $\scrT(H^1\otimes_{\dbZ} \dbQ_q)$ of degree not bigger than a given $d(A,q)\in\dbN$, iff $\tilde G_{\overline{\dbQ_q}}$ is the subgroup of $GL(H^1_B\otimes_{\dbZ} \overline{\dbQ_q})$ fixing some homogeneous tensors of $\scrT(H^1_B\otimes_{\dbZ} \overline{\dbQ_q})$ of degree not bigger than $d(A,q)$.

Let $E^{T^0}$ be a number field over which $T^0$ splits. Let $\widetilde {GL}(q):=GL(H_B^1\otimes_{\dbZ} \dbQ_q)(\overline{\dbQ_q})$. The key point is: there is (cf. 4.3.4.1) a finite set $\{\tilde G_j|j\in S(A)\}$ of split, reductive subgroups of $GL(H^1_B(A_{\dbC},\dbQ)\otimes_{\dbQ} E^{T^0})$ having $T^0_{E^{T^0}}$ as the connected component of the origin of their centers and such that for any prime $q$, $\tilde G_{\overline{\dbQ_q}}$ is $\widetilde {GL}(q)$-conjugate to one of the subgroups $\tilde G_{j\overline{\dbQ_q}}$ (it does not matter which embedding $E^{T^0}\hookrightarrow\overline{\dbQ_q}$ we choose), $j\in S(A)$. The existence of $d(A)$ is implied by the first remark and by the fact that $S(A)$ is finite.   
\finishproclaim

\Proclaim{4.3.4.2.2. The $\Pi(A)$ invariant.} \rm
$\Pi(A)\in\dbN$ is defined as the smallest number such that for any representation $W_{\dbZ}$ over $\dbZ$ of dimension not bigger than $r_A$ of the Lie algebra $\grg_{\dbZ}$ of a split, semisimple group (over $\dbZ)$ of dimension not bigger than $d_A$, the symmetric bilinear form $TR$ on $\grg_{\dbZ}\fracwithdelims[]1{\Pi(A)}$ induced by the trace form on $\End(W_{\dbZ}\fracwithdelims[]1{\Pi(A)})$ is perfect. 4.3.4.1 guarantees $\Pi(A)\in\dbN$ is well defined ($TR$ depends only on the representation of $\grg_{\dbZ}\otimes_{\dbZ}\dbQ$ on $W_{\dbZ}\otimes_{\dbZ}\dbQ$).
\finishproclaim

\Proclaim{4.3.4.2.3. Notation.} \rm
Let $\grg$ be as in 4.3.4.2 iv). For a non-trivial, simple $\grg$-module $W_0$ of finite dimension, we denote by 
$$k(\grg,W_0)\in\dbN\cup\{0\}$$ 
the greatest coefficient of the decomposition of the highest weight of $W_0$ as a linear combination with non-negative coefficients of fundamental weights (here weights are w.r.t. some fixed Cartan Lie subalgebra of $\grg$; see [Hu1, ch. 20]).  
\finishproclaim

\Proclaim{4.3.4.3. Conclusion.} \rm
4.3.4.1 assures the existence of a smallest number $N_2(A)\in\dbN$ for which the needed properties i) to vi) of 4.3.4.2 are satisfied (cf. 4.3.4.2.1-2). Obviously $N_2(A)$ depends only on $g_A$. As $d_A\Ge 4$, we deduce from 4.3.4.1 iii) that $N_2(A)\Ge 5$. From now on we also assume that $p\Ge N_2(A)$. So $p\Ge N(A)$.
\finishproclaim

\Proclaim{4.3.4.4. Some decompositions.} \rm
Let $T_p^0:=Z^0(G_{\dbZ_p})$. It is a split torus (as $G_{\dbZ_p}$ is split). Let $H^1_{\prime}=\oplus_{\gamma\in\Gamma} H_{\prime}^{\gamma}$ be the maximal direct sum decomposition associated to the faithful representation $T_p^0\hookrightarrow GL(H_{\prime}^1)$; so $\Gamma$ is a set of characters of $T_p^0$ such that $\beta\in T_p^0(\dbZ_p)$ acts on $H_{\prime}^{\gamma}$ as the multiplication with $\gamma(\beta)\in\dbG_m(\dbZ_p)$.

For any $\gamma\in\Gamma$, let $H_{\prime}^{\gamma}\fracwithdelims[]1p=\oplus_{i\in I_{\gamma}} H^i\fracwithdelims[]1p$ be a direct sum decomposition in irreducible $G_{\dbQ_p}^{\der}$-representations. Let 
$$I:=\cup_{\gamma\in\Gamma} I_{\gamma}$$ 
(disjoint union). We get:
$$H_{\prime}^1\fracwithdelims[]1p=\oplus_{i\in I} H^i\fracwithdelims[]1p.$$
Let $H^i:=H^i\fracwithdelims[]1p\cap H_{\prime}^1$, $i\in I$.
 Let $I^{\text{tr}}$ be the subset of $I$ formed by those $i$ such that $G^{\der}_{\dbQ_p}$ acts trivially on $H^i\fracwithdelims[]1p$. Let $I^{\text{ntr}}:=I\setminus I^{\text{tr}}$. 
\finishproclaim

\Proclaim{4.3.4.5. Absolutely irreducible representations.} \rm
Let $i\in I^{\text{ntr}}$. The representation of $G^{\der}_{\dbQ_p}$ on $H^i[{1\over p}]$ is a finite tensor product $\otimes_{j\in I_i} \rho_j^i$ of irreducible representations factoring through quotients of $G^{\der}_{\dbQ_p}$ having split, simple adjoints. Each such $\rho_j^i$ extends to a representation $\tilde\rho_j^i$ of $G^{\der}_{\dbZ_p}/\Ker_i$, where $\Ker_i$ is the Zariski closure in $G^{\der}_{\dbZ_p}$ of $\Ker(\rho_j^i)$ (cf. [Ja, 10.4 of Part I]). $G^{\der}_{\dbZ_p}/\Ker_i$ is a semisimple group. The special fibre of $\tilde\rho_j^i$ is absolutely irreducible, cf. [Bo1, 6.4] (as $p\Ge N_2(A)$, 4.3.4.2 v) implies that loc. cit. applies). From this, [Ja, 10.4 and 10.9 of Part I] and [Bo1, 7.2-3] we get that the representation of $G^{\der}_{\dbF_p}$ on $H^i/pH^i$ is a tensor product of absolutely irreducible representations (indexed by $I_i$) and so it is absolutely irreducible. Obviously, if $i\in I^{\text{tr}}$, the representation of $G^{\der}_{\dbF_p}$ on $H^i/pH^i$ is trivial and 1 dimensional and so irreducible.
\finishproclaim

\Proclaim{4.3.4.5.1. Claim.} 
We can assume we have a direct sum decomposition
$$H_{\prime}^1=\oplus_{i\in I} H^i.$$
\finishproclaim

\proof: 
For each $\gamma\in\Gamma$ and for any $i\in I_{\gamma}$, let $J_i$ be the subset of $I_{\gamma}$ formed by elements $j$ such that the representations of $G^{\der}_{\dbQ_p}$ on $H^i[{1\over p}]$ and on $H^j[{1\over p}]$ are isomorphic.  Based on 4.3.4.5, it is enough to show that we can assume $\oplus_{j\in J_i} H^j$ is a direct summand of $H^1_{\prime}$. The case $i\in I^{\text{tr}}$ is trivial and so we assume $i\in I^{\text{ntr}}$. 

Let 
$$H(i):=H_{\prime}^1\cap \oplus_{j\in J_i} H^j[{1\over p}].$$ 
Let $B(i)$ be a Borel subgroup of the image $G(i)$ of $G^{\der}_{\dbZ_p}$ in $GL(H(i))$ (based on [Va1, 3.1.2.1 c)] we get easily that $G(i)$ is a reductive, closed subgroup of $GL(H(i))$. Let $T(i)$ be a split, maximal torus of $B(i)$. Let $w$ be the highest weight of the representation of $T(i)_{\dbQ_p}$ on $H^i[{1\over p}]$; here ``highest" is w.r.t. the basis of roots of the root system of the action of $T(i)_{\dbQ_p}$ on $\Lie(G^{\der}_{\dbQ_p})$ which corresponds naturally to $\Lie(B(i)_{\dbQ_p})$. Let $V(i)$ be the maximal direct summand of $H(i)$ on which $T(i)$ acts as scalar multiplication by its character corresponding to $w$. Its rank is the same as the number of elements of $J_i$ (see [Hu1, p. 108]). 

For $j\in J_i$, we choose $v_j\in V(i)$. We assume that the restrictions of $v_j$'s mod $p$ are forming an $\dbF_p$-basis $\scrB(i)$ of $V(i)/pV(i)$. We take $H^j[{1\over p}]$ to be the cyclic $\Lie(G(i)_{\dbQ_p})$-submodule of $H(i)[{1\over p}]$ generated by $v_j$ (cf. loc. cit.). As $\scrB(i)$ is an $\dbF_p$-basis of $V(i)/pV(i)$, from 4.3.4.5 we get that the intersection 
$$H^j/pH^j\cap \cap_{l\in J_i\setminus\{j\}} H^l/pH^l$$
is trivial, $\forall j\in J_i$. We conclude: $\oplus_{j\in J_i} H^i$ is $H(i)$. This ends the proof.

\Proclaim{4.3.4.6. Basic tensors.} \rm
Let $(v_{\alpha})_{\alpha\in\scrJ_0}$ be the family of tensors formed by all homogeneous tensors of $\scrT(H_{\prime}^1)$ fixed by $G_{\dbQ_p}$ and of degree at most $\max\{d(A),4\}$; so the set $\scrJ_0$ is infinite. In particular, it contains the following 4 types of tensors.

\medskip
\item{{\bf 1)}}
The projection $p_{\gamma}$ of $H_{\prime}^1$ on $H_{\prime}^{\gamma}$ associated to the direct sum decomposition $H_{\prime}^1=\oplus_{\gamma\in\Gamma} H^{\gamma}$ ($\gamma\in\Gamma$). So $p_{\gamma}\in\End(H^1_{\prime})=H^1_{\prime}\otimes_{\dbZ_p} (H^1_{\prime})^*$.

\smallskip
\item{{\bf 2)}}
The tensor $\pi(\Lie(G_{\dbQ_p}^{\der}),H_{\prime}^1\fracwithdelims[]1p)$ defined in [Va1, 4.2.1] (4.3.4.2.2 implies that it is enveloped by --i.e. it is integral w.r.t.-- $H_{\prime}^1$). 

\smallskip
\item{{\bf 3)}}
The tensors $B$ and $B^*$ defined in [Va1, 4.3.2] but for $\Lie(G_{\dbQ_p}^{\der})\subset\grg\grl(H_{\prime}^1\fracwithdelims[]1p)$ (cf. 2)).

\smallskip
\item{{\bf 4)}}
The projection $p^i$ of $H_{\prime}^1$ on $H^i$ associated to the direct sum decomposition of 4.3.4.5.1 ($i\in I$). So $p_i\in\End(H^1_{\prime})=H^1_{\prime}\otimes_{\dbZ_p} (H^1_{\prime})^*$.
\finishproclaim

\Proclaim{4.3.4.6.1. Remarks.} \rm
{\bf 1)} The tensors $\pi(\Lie(G_{\dbQ_p}^{\der}),H^1_{\prime}\fracwithdelims[]1p)$, $B$ and $B^*$ are fixed by any automorphism of $H_{\prime}^1\fracwithdelims[]1p$ normalizing $G^{\der}_{\dbQ_p}$.

{\bf 2)} We have $\deg(v_{\alpha})<p-1$, $\forall\alpha\in\scrJ_0$, cf. the inequalities 
$$p\Ge N_2(A)>\max\{2d_A,d(A)+1\}\Ge 8.$$ 
As $p\Ge 8$, there is $s(p)\in\dbN$ such that $p^{s(p)}$ times $\pi(G_{\dbQ_p}):=\pi(\Lie(G_{\dbQ_p}),H^1_{\prime}\fracwithdelims[]1p)$ is $v_{\alpha^0}$, for some $\alpha^0\in\scrJ_0$.
\finishproclaim

\Proclaim{4.3.4.7. The enlarged family of tensors.} \rm
The subgroup of $GL(H_{\prime}^1\fracwithdelims[]1p)$ fixing  $v_{\alpha}$, $\forall\alpha\in\scrJ_0$, is a subgroup of $\tilde G_{\dbQ_p}$ containing $G_{\dbQ_p}$, cf. 4.3.4.6 vi) and the definition of $(v_{\alpha})_{\alpha\in\scrJ_0}$. Let $(v_{\alpha})_{\alpha\in\scrJ}$, with $\scrJ$ a set containing $\scrJ_0$, be an enlarged family of tensors of $\scrT(H_{\prime}^1)$ such that $G_{\dbQ_p}$ is the subgroup of $GL(H_{\prime}^1\fracwithdelims[]1p)$ fixing $v_{\alpha}$, $\forall\alpha\in\scrJ$. We assume that all $p$-components of \'etale components of Hodge cycles of $A_{\dbC}$ which are elements of $\scrT(H^1_{\prime})$, are elements of this family of tensors indexed by $\scrJ$. This forces the set $\scrJ\setminus\scrJ_0$ to be infinite. 
\finishproclaim

\Proclaim{4.3.4.8. Claim.} 
{\bf a)} The family of tensors $(v_{\alpha})_{\alpha\in\scrJ_0}$ of $\scrT(H_{\prime}^1)$ is $\dbZ_p$-very well positioned for $G_{\dbQ_p}$. 

{\bf b)} If $\tilde G_{\dbQ_p}=G_{\dbQ_p}$, then this family is $H_{\prime}^1$-representation well positioned for $G_{\dbQ_p}$.
\finishproclaim

\proof
Let $R$ be an integral, faithfully flat $\dbZ_p$-algebra. Let $M\subset R\otimes_{\dbZ_p} H_{\prime}^1\fracwithdelims[]1p$ be a free $R$-module such that $M\fracwithdelims[]1p=R\otimes_{\dbZ_p} H_{\prime}^1\fracwithdelims[]1p$ and it envelopes the family of tensors $(v_{\alpha})_{\alpha\in\scrJ_0}$. Let $S:=R\fracwithdelims[]1p$.

The subfamily of tensors $(p_{\gamma})_{\gamma\in\Gamma}$ is $\dbZ_p$-well positioned for $T_p^0$ (see [Va1, 4.3.13]). So (cf. 4.3.4.6 1)) the Zariski closure $T^0_{pR}$ of $T_{pS}^0$ in $GL(M)$ is an $R$-torus.

The subfamily of tensors formed by $B$, $B^*$ and  $\pi(\Lie(G_{\dbQ_p}^{\der}),H_{\prime}^1\fracwithdelims[]1p)$ is $\dbZ_p$-very well positioned for $G_{\dbQ_p}^{\der}$ (cf. [Va1, 4.3.10 b) and 4.3.10.1 1)] and 4.3.4.2 iv)). So (cf. 4.3.4.6 2) and 3)) the Zariski closure of $G_S^{\der}$ in $GL(M)$ is a semisimple group over $R$.

From [Va1, 3.1.6] we deduce that the Zariski closure of $G_S$ in $GL(M)$ is a reductive group over $R$. So the family of tensors $(v_{\alpha})_{\alpha\in\scrJ_0}$ is $\dbZ_p$-very well positioned for $G_{\dbQ_p}$. This ends the proof of a) as well as of 4.3.4.

We now prove b). We assume that $R$ is a strictly henselian DVR $V$, that $\tilde G_{\dbQ_p}=G_{\dbQ_p}$, and that for any $\alpha\in\scrJ_0$, $v_{\alpha}$ has a non-zero image in $H_{\prime}^1/pH_{\prime}^1$ iff it has a non-zero image in $M/\pi_VM$, with $\pi_V$ a uniformizer of $V$ (cf. 2.3.1). We need to show: there is an isomorphism $M\arrowsim H_{\prime}^1\otimes_{\dbZ_p} V$ taking $v_{\alpha}$ into $v_{\alpha}$, $\forall\alpha\in\scrJ_0$. 4.3.4.7 allows us to express the existence of such an isomorphism in terms of a reduced $V$-scheme (which is sort of a torsor of $G_V$) not being a scheme over $\Spec(V/\pi_V\times V\fracwithdelims[]1p)$. So we can assume $V$ is complete and we are allowed to replace $V$ by its normalization $V_1$ in a finite field extension of $V\fracwithdelims[]1p$.

We have a direct sum decomposition of $G_V$-modules $M=\oplus_{i\in I} M^i$, with $M^i=p^i(M)$ (cf. 4.3.4.6 4)).

We have two hyperspecial subgroups of $G_{\dbQ_p}(V\fracwithdelims[]1p)$: $G_{V\fracwithdelims[]1p}(M)$ and $G_{V\fracwithdelims[]1p}(H_{\prime}^1\otimes_{\dbZ_p} V)$. They are $G_{\dbQ_p}^{\ad}(V\fracwithdelims[]1p)$-conjugate (cf. [Ti, p. 47]) and so by replacing $V$ by a $V_1$ as above and by replacing $M\otimes_V V_1$ with $g_1(M\otimes_V V_1)g_1^{-1}$, for some $g_1\in G_{\dbQ_p}(V_1\fracwithdelims[]1p)$, we can assume they are the same.

On the other hand, the representation of the special fibre of $G_V$ on $M^i/\pi_VM^i$ is absolutely irreducible, cf. 4.3.4.5. So for any $i\in I$, $M^i$ is $\delta_i H^i\otimes_{\dbZ_p} V$ for some $\delta_i\in\dbG_m(V\fracwithdelims[]1p)$. For $i\in I$, we write $\delta_i=\pi_V^{n_i}u_i$, with $n_i\in\dbZ$ and $u_i\in\dbG_m(V)$. Let 
$$\delta\in GL(H^1_{\prime}\otimes_{\dbZ_p} V\fracwithdelims[]1p)$$
be such that it takes $m_i\in H^i$ into $\pi_V^{n_i}m_i$, $\forall i\in I$. It centralizes $G_V$. We have 
$$\delta(H_{\prime}^1\otimes_{\dbZ_p} V)=M.$$ 
As the family of tensors $(v_{\alpha})_{\alpha\in\scrJ_0}$ is enveloped by $M$ and by $H_{\prime}^1\otimes_{\dbZ_p} V$, we deduce by induction that the family of tensors $(\delta^{-n}(v_{\alpha}))_{\alpha\in\scrJ_0}$ is enveloped by $H_{\prime}^1\otimes_{\dbZ_p} V$, $\forall n\in\dbN$; argument: $\delta^{-n}(v_{\alpha})$ is a $V$-linear combination of members of the family $(v_{\alpha})_{\alpha\in\scrJ_0}$, cf. 4.3.4.6.

Let $\alpha\in\scrJ_0$. We now check that $\delta(v_{\alpha})=v_{\alpha}$. We can assume that the image of $v_{\alpha}$ in $\scrT(H_{\prime}^1/pH_{\prime}^1)$ is non-zero. We start writing 
$$v_{\alpha}=\sum_{s\in\dbZ} v_{\alpha}^s,$$
with $v_{\alpha}^s$ a homogeneous tensor of $\scrT(H^1_\prime)$ of the same degree as $v_{\alpha}$ and such that 
$$\delta^{-1}(v_{\alpha}^s)=\pi_V^sv_{\alpha}^s.$$ 
We immediately get: $v_{\alpha}^s=0$ if $s<0$ and the image of $v_{\alpha}^0$ in $\scrT(H_{\prime}^1/pH_{\prime}^1)$ is non-zero. Using a suitable quotient of the difference $v_{\alpha}-\delta^{-1}(v_{\alpha})$ by a non-negative, integral power of $\pi_V$ and the fact this tensor is a $V$-linear combination of similar tensors $w_{\alpha}^s$ (of $\scrT((H^1_\prime)$) with $s\in\dbN$, we get (cf. the fact that $\forall\alpha^\prime\in\scrJ_0$, the image of $v_{\alpha^\prime}$ in $\scrT(H_{\prime}^1/pH_{\prime}^1)$ is non-zero iff its image in $\scrT(M/\pi_VM)$ is non-zero): $v_{\alpha}^s=0$, $\forall s\in\dbN$. So $v_{\alpha}=v_{\alpha}^0$ and so $\delta(v_{\alpha})=v_{\alpha}$.

So (in our present context) 2.3.1 b) is satisfied for such a ring $R$ while 2.3.1 c) follows from 4.3.4.7. These, together with a), end the proof of b) and so of the Claim.
 
\Proclaim{4.3.5. An adequate isogeny.} \rm
Let $A_{\dbC}^{\prime}$ be the abelian variety over $\dbC$ which is $\dbZ\fracwithdelims[]1p$-isogeneous to $A_{\dbC}$ and satisfies
$$H_{\text{\'et}}^1(A_{\dbC}^\prime,\dbZ_p)=H_{\prime}^1.$$ 
It is defined over a subfield $E^{\prime}$ of $\dbC$ which is a  finite field extension of $E$. We denote by $i_{E^{\prime}}\colon E^{\prime}\hookrightarrow\dbC$ the natural inclusion. Let $A^{\prime}$ be the (abelian variety) model of $A_{\dbC}^\prime$ over $E^{\prime}$. $A^{\prime}$ is obtained from $A_{E^{\prime}}$ by taking its quotient through a finite group subscheme $C$ of  $A_{E^{\prime}}$. If $n\in\dbN$ is such that $p^nH^1_{\prime}\subset H^1$, then the finite (abstract) group $C(\dbC)$ corresponds to the quotient group $H^1/p^nH_{\prime}^1$ (i.e. we have $H_{\text{\'et}}^1(A_{\dbC}/C_{\dbC},\dbZ_p)=p^nH_{\prime}^1$). The pull back of the polarization $p_A$ of $A$ under the multiplication endomorphism $p^n\colon A\to A$ descends to give birth to a principal polarization $p_{A^{\prime}}$ of $A^{\prime}$. The perfect form $p_{A^{\prime}}\colon H_{\prime}^1\otimes_{\dbZ_p} H_{\prime}^1\to\dbZ_p(-1)$ (of $\Gal(E^\prime)$-modules) induced by the polarization $p_{A^{\prime}}$ of $A^\prime$ is precisely $p_A$. 

From 4.3.3.0 we get that $A^{\prime}$ is a product of absolutely simple abelian varieties and that $p_{A^{\prime}}$ is a product of principal polarizations of these simple factors.
\finishproclaim

\Proclaim{4.3.6. Level structures.} \rm
Passing to a finite field extension of $E^{\prime}$ (unramified above primes not dividing $(N_1(A)-1)!$, cf. 4.3.1), we can assume that both $A^{\prime}$ and $A$ have a level-$4$ symplectic similitude structure. By this we mean similitude  isomorphisms 
$$k_4^{\prime}\colon (L/4L,\psi)_{E^\prime}\arrowsim (A^{\prime}[4],p_{A^{\prime}})$$
and $k_4\colon (L/4L,\psi)_E\arrowsim (A[4],p_A)$ of symplectic spaces over $\dbZ/4\dbZ$. Here $\psi$ is the symplectic form on 
$$W:=H_{1B}(A_{\dbC},\dbQ)=H_{1B}(A_{\dbC}^{\prime},\dbQ)$$ 
defined by the polarization $p_A$ of $A$, $L$ is a fixed $\dbZ$-lattice of $W$ such that $\psi\colon L\otimes_{\dbZ} L\to\dbZ$ is a perfect form, and we also identify $L/4L$ with a group scheme over $\Spec(\dbZ)$. For convenience, we take $L$ to be the first group of the Betti homology of $A_{\dbC}^\prime$ with coefficients in $\dbZ$; accordingly, we have a natural identification $L/4L_{E^\prime}=A^\prime[4]$ and we take $k_4^\prime$ to be defined by the identity isomorphism of $L/4L_{E^\prime}$. Let $K_p:=G\Sp(L,\psi)(L\otimes_{\dbZ} \dbZ_p)$. It is a hyperspecial subgroup of $G\Sp(W,\psi)(\dbQ_p)$. All below level-$4$ symplectic similitude structures are defined as above, i.e. are w.r.t. $L/4L$.
\finishproclaim

\Proclaim{4.3.6.1. Moduli setting.} \rm
We get a Siegel modular variety $\Sh(G\Sp(W,\psi),S)$ and an injective map 
$$f_A\colon (H_A,X_A)\hookrightarrow (G\Sp(W,\psi),S)$$ 
between Shimura pairs; here $X_A$ is the Hermitian symmetric domain defined by $H_A$ and by the Hodge structure of $W$ corresponding to $A_{\dbC}$ or to $A^{\prime}_{\dbC}$ (see also [Va1, 2.12 3)]). We also get a Shimura quadruple $(G\Sp(W,\psi),S,K_p,p)$. Let $\scrM$ be its integral canonical model (cf. [Va1, 3.2.9]).
So $GSp(L,\psi)(L\otimes_{\dbZ} \dbA_f^p)$ acts on $\scrM$.

We have $\dim_{\dbQ}(W)=2g_A$. Let $\scrA_{g_A,1,4}$ be the moduli scheme over $\Spec(\dbZ\fracwithdelims[]12)$ parameterizing isomorphism classes of principally polarized abelian schemes of dimension $g_A$ (over $\Spec(\dbZ\fracwithdelims[]12)$-schemes) having level-$4$ symplectic similitude structure. 
Its fibre over $\dbZ_{(p)}$ is the quotient of $\scrM$ by the compact, open subgroup of $G\Sp(W,\psi)(L\otimes_{\dbZ} \widehat{\dbZ}^p)$ of elements acting trivially on $L/4L$. Here $\widehat{\dbZ}^p$ is $\widehat{\dbZ}$ with the $p$-component omitted; so $\widehat{\dbZ}=\widehat{\dbZ}^p\times\dbZ_p$. 
Corresponding to the two triples $(A^{\prime},p_{A^{\prime}},k_4^{\prime})$ and $(A,p_A,k_4)$ (of elements defined above) we get two morphisms 
$$l^{\prime}\colon\Spec(E^{\prime})\to\scrA_{g_A,1,4}$$ 
and respectively $l\colon\Spec(E)\to\scrA_{g_A,1,4}$.
\finishproclaim

\Proclaim{4.3.6.1.1. Remark.} \rm
If $H_{A\dbQ_p}$ is a split group, then (cf. 4.3.3.1) the Zariski closure $H_{A\dbZ_{(p)}}$ of $H_A$ in $GL(L\otimes_{\dbZ} \dbZ_{(p)})$ is a reductive group. So we get another Shimura quadruple $(H_A,X_A,H_{A\dbZ_{(p)}}(H^1_{\prime}),v_A)$, with $v_A$ a fixed prime of the reflex field $E(H_A,X_A)$ dividing $p\Ge 8$. Let $\scrN$ be its integral canonical model, cf. [Va1, 6.4.1]. 
We get an injective map 
$$(H_A,X_A,H_A(H^1_{\prime}),v_A)\hookrightarrow (G\Sp(W,\psi),S,K_p,p),$$ 
and so, cf. [Va1, rm. 4) of 3.2.7], a natural morphism $\scrN\to\scrM$. In fact we get a SHS $(f_A,L\otimes_{\dbZ} \dbZ_{(p)},v_A)$ (as $p\Ge 2d_A$, [Va1, 5.8.6] applies).  
\finishproclaim
 
\Proclaim{4.3.6.2. Remark.} \rm
Passing, if needed, to a finite field extension of $E^{\prime}$, we can assume 
$$\rho(\Gal(E^{\prime}))\subset G_{\dbQ_p}(\dbQ_p),$$ 
i.e. we can assume that the tensors of the family $(v_{\alpha})_{\alpha\in\scrJ}$ are fixed by $\Gal(E^{\prime})$.
\finishproclaim

\Proclaim{4.3.7. $p$-adic setting.} \rm
Passing to a finite field extension of $E$ unramified above primes relatively prime to $(N_1(A)-1)!$, we can assume we have natural inclusions $E(H_A,X_A)\subset E\subset E^\prime\subset\dbC$. Let $v^\prime$ be a prime of $E^\prime$ dividing $v_A$ and let $v$ be the prime of $E$ divided by it. Warning: $v^{\prime}$ might be ramified over $p$. Let $e$ be the index of ramification of $v^{\prime}$ over $p$. Let $\dbF:=\overline{k(v^{\prime})}$. 
Let $V_0:=W(\dbF)$. Let $V$ be the completion of the maximal unramified cover of the ring of integers $O_{v^\prime}$ of $E_{v^{\prime}}^{\prime}$. $V$ is a finite, flat DVR extension of $V_0$ of degree $e$. Let $K_0:=V_0\fracwithdelims[]1p$ and $K:=V\fracwithdelims[]1p$. 

Let $A_{V_0}$ be the abelian scheme over $V_0$ having $A_{K_0}$ as its generic fibre and let $A_V^{\prime}$ be the abelian scheme over $V$ having $A_K^{\prime}$ as its generic fibre (cf. the inequality $p\Ge N_1(A)$, the definition of $N_1(A)$, and the fact that the good reduction property is preserved by passage to an isogeny). We still denote by $p_A$ and $p_{A^{\prime}}$ the principal polarizations of $A_{V_0}$ and respectively of $A_V^{\prime}$, whose restrictions to generic fibres are obtained from $p_A$ and respectively $p_{A^{\prime}}$ via pull backs.
\finishproclaim

\Proclaim{4.3.7.1. Notations.} \rm
We $l_{V_0}\colon\Spec(V_0)\to\scrA_{g_A,1,4}$ and $l_V\colon\Spec(V)\to\scrA_{g_A,1,4}$ be the morphisms corresponding respectively to the triples $(A_{V_0},p_A,k_4)$ and $(A_V^{\prime},p_{A^{\prime}},k_4^{\prime})$. Their generic fibres are the extension to $K_0$ and respectively to $K$ of $l$ and respectively of $l^\prime$. 

From the definition of $A^{\prime}$, we get a $\dbZ\fracwithdelims[]1p$-isogeny $$i_V\colon A_V^{\prime}\to A_V$$ 
inducing the isomorphism $H^1_{\text{\'et}}(A_{\Kbar},\dbQ_p)\arrowsim H^1_{\text{\'et}}(A_{\Kbar}^{\prime},\dbQ_p)$, defined by the inclusion $H^1\subset H^1_{\prime}\fracwithdelims[]1p$. At the level of generic fibres we get a $\dbZ\fracwithdelims[]1p$-isogeny $i_{\dbF}\colon A_{\dbF}^{\prime}\to A_{\dbF}$.
\finishproclaim

\Proclaim{4.3.8. An assumption.} \rm
From now on we assume (cf. 4.2.12) that $G_{\dbQ_p}$ is not a torus. So $G_{\dbQ_p}^{\der}$ is a non-trivial, semisimple group.
\finishproclaim

\Proclaim{4.4. Local deformation of $A_V^{\prime}$ and its de Rham tensors.} \rm
For what follows we rely heavily on [Va1, 5.2-3], [Va2, 3.6] and 4.3.
\finishproclaim

\Proclaim{4.4.1. Notations.} \rm
Let $\pi_V$ be a uniformizer of $V$. Let $Re$, $\Rtil e$, $B^+(V)$ and $\beta$ have the same significance as in [Va1, 5.2.1]. So (cf. 2.1 and loc. cit.) $\Rtil e$ (resp. $Re$) is the piano ring of $\dbF$ of resonance $e$ (resp. of convergent resonance $e$), while $\beta$ is a suitable element of the Fontaine's ring $B^+(V)$ of $V$.  

Let $(M^{\prime},F^1(M^{\prime}),\Phi_{M^{\prime}},\nabla^{\prime})$ be the filtered $F$-crystal over $Re/pRe$ obtained by taking the dual of the Lie algebra of the universal vector extension of the abelian variety $A_V^{\prime}$ (or of the $p$-divisible group associated to $A_V^{\prime}$) (cf. [Me], [BBM]; see also [Fa3] and [Va1, 5.2.2]). 

Let $M_V^{\prime}:=M^{\prime}\otimes_{Re} V$; here the natural $V_0$-epimorphism $Re\twoheadrightarrow V$ is defined by $\pi_V$, cf. 2.1 Fact and [Va1, 5.2.1]. $M_V^\prime$ is canonically identified with $H_{dR}^1(A_V^{\prime}/V)$, cf. [Va1, 5.2.2]. Let $$M_0:=H_{\text{crys}}^1(A_{\dbF}/W(\dbF))$$ and 
$$M_0^{\prime}:=H_{\text{crys}}^1(A_{\dbF}^{\prime}/W(\dbF)).$$ 
Let $F^1_0$ be the Hodge filtration of $M_0$ defined by $A_{V_0}$. Let $\phi_0$ and $\phi^{\prime}$ be the $\sigma$-linear endomorphisms of $M_0$ and respectively of $M_0^{\prime}$. The isocrystals $(M_0\fracwithdelims[]1p,\phi_0)$ and $(M_0^\prime\fracwithdelims[]1p,\phi^\prime)$ are naturally isomorphic (cf. the existence of $i_{\dbF}$).
\finishproclaim

\Proclaim{4.4.2. Fontaine's comparison theory.} \rm
Applying [Fa3, th. 7] to the $p$-divisible group of $A_V^{\prime}$, we get (cf. also [Va1, 5.2.5]) an inclusion of filtered $B^+(V)$-modules
$$\rho^{\prime}\colon M^{\prime}\otimes_{Re} B^+(V)\hookrightarrow H^1_{\prime}\otimes_{\dbZ_p} B^+(V),$$
which becomes an isomorphism $\rho_1^{\prime}$ by inverting $p\beta$.
Similarly, cf. [Va1, 5.2.2.1 and 5.2.4] (and [Va2, AE.0] which points out that we need to invert $p\beta$ and not just $p$), for $A_V$ we get an isomorphism of filtered $B^+(V)\fracwithdelims[]1{p\beta}$-modules
$$\rho_1\colon M_0\otimes_{V_0} B^+(V)\fracwithdelims[]1{p\beta}\arrowsim H^1\otimes_{\dbZ_p} B^+(V)\fracwithdelims[]1{p\beta}.$$
\noindent
$\rho_1$ is obtained (by inverting $p\beta$) from a monomorphism of filtered $B^+(V)$-modules
$$\rho\colon M_0\otimes_{V_0} B^+(V)\hookrightarrow H^1\otimes_{\dbZ_p} B^+(V),$$
as it can be easily checked starting from the fact that $A_V$ is definable over $V_0$.  
\finishproclaim

\Proclaim{4.4.2.1. Isomorphisms defined by isogenies.} \rm
The $\dbZ\fracwithdelims[]1p$-isogeny $i_V$ gives birth to two $B^+(V)\fracwithdelims[]1{p\beta}$-isomorphisms 
$$i_V^{dR}\colon M_0\otimes_{\dbZ_p} B^+(V)\fracwithdelims[]1{p\beta}\arrowsim M^{\prime}\otimes_{Re} B^+(V)\fracwithdelims[]1{p\beta}$$
and
$$i_V^{\text{\'et}}\colon H^1\otimes_{\dbZ_p} B^+(V)\fracwithdelims[]1{p\beta}\arrowsim H_{\prime}^1\otimes_{\dbZ_p} B^+(V)\fracwithdelims[]1{p\beta}$$
such that $\rho_1^{\prime}\circ i_V^{dR}=i_V^{\text{\'et}}\circ\rho_1$ (cf. the functorial part of [Fa3, th. 7]). 

$\rho^\prime$, $\rho$, $i_V^{dR}$ and $i_V^{\text{\'et}}$ commute with the natural (see [Fa3]) Frobenius endomorphisms.
\finishproclaim

\Proclaim{4.4.3. Tensors in the crystalline context.} \rm
Let $(w_{\alpha}^{\prime})_{\alpha\in\scrJ}$ be the family of tensors of $\scrT(M^{\prime}\fracwithdelims[]1p)$ obtained from $(v_{\alpha})_{\alpha\in\scrJ}$ through  $\rho^{\prime}$ (cf. Fontaine's comparison theory and 4.3.6.2). They are elements of the $F^0$-filtration of $\scrT(M^{\prime}\fracwithdelims[]1p)$, fixed by $\Phi_{M^{\prime}}$ and annihilated by $\nabla^{\prime}$. [Fa3, cor. 9] and 4.3.4.6.1 2) guarantee that $w_{\alpha}^{\prime}\in\scrT(M^{\prime})$, $\forall\alpha\in\scrJ_0$. All these are as in [Va1, 5.2.8-9].
\finishproclaim

\Proclaim{4.4.3.1. Key reductiveness property.} \rm
From the fact that $(v_{\alpha})_{\alpha\in\scrJ_0}$ is a $\dbZ_p$-well positioned family of tensors for the group $G_{\dbQ_p}$ and from 4.4.3, we deduce (as in [Va1, 5.2.12]) that the Zariski closure $G_{Re}$ in $GL(M^\prime)$ of the subgroup of $GL(M^{\prime}\fracwithdelims[]1p)$ fixing $w_{\alpha}^{\prime}$, $\forall\alpha\in\scrJ$, is a reductive group over $Re$. [Va1, 5.2.1.1] shows that the notation $G_{Re}$ is justified.
\finishproclaim

\Proclaim{4.4.3.2. Remark.} \rm
Under $(i_V^{dR})^{-1}$, $w_{\alpha}^{\prime}$ is mapped into an element $w_{\alpha}$ of the $F^0$-filtration of $\scrT(M_0\fracwithdelims[]1p)$ fixed by $\phi_0$ ($\alpha\in\scrJ$). For instance, this can be read out from the fact that $i_V^{dR}$ is the extension to $B^+(V)\fracwithdelims[]1{p\beta}$ of a similar isomorphism over $Re\fracwithdelims[]1p$.
\finishproclaim

\Proclaim{4.4.4. Cocharacters.} \rm
We consider the canonical split cocharacter 
$$\mu\colon\dbG_m\to GL(M_0)$$ 
of $(M_0,F^1_0,\phi_0)$. It factors through $G_{V_0}$ and produces a direct sum decomposition $M_0=F_0^1\oplus F_0^0$, with $\gamma\in\dbG_m(V_0)$
acting through $\mu$ on $F_0^i$ as the multiplication with $\gamma^{-i}$, $i=\overline{0,1}$. 
\finishproclaim

\Proclaim{4.4.4.1. Transfer via isogeny.} \rm
The $\dbZ\fracwithdelims[]1p$-isogeny $i_V$ transforms the cocharacter $\mu$ into a cocharacter $\mu^{\prime}\colon\dbG_m\to GL(M_V^{\prime}\fracwithdelims[]1p)$ factoring through the subgroup $G_K$ of $GL(M_V^{\prime}\fracwithdelims[]1p)$ fixing the image $t_{\alpha}\in\scrT(M_V^{\prime}\fracwithdelims[]1p)$ of $w_{\alpha}^{\prime}$ under the epimorphism $M^{\prime}\fracwithdelims[]1p\twoheadrightarrow M_V^{\prime}\fracwithdelims[]1p$ of $V_0$-modules, $\forall\alpha\in\scrJ$. It produces a direct sum decomposition $M_V^{\prime}\fracwithdelims[]1p=F^1_{\prime}\oplus F^0_{\prime}$, with $\gamma\in\dbG_m(K)$ acting through $\mu^{\prime}$ as the multiplication with $\gamma^{-i}$ on $F^i_{\prime}$,  $i=\overline{0,1}$; $F^1_{\prime}$ is the Hodge filtration of $M_V^{\prime}\fracwithdelims[]1p=H^1_{dR}(A_K^{\prime}/K)$ defined by $A_K^{\prime}$.
\finishproclaim

\Proclaim{4.4.4.2. Integral variant.} \rm
Repeating the arguments of [Va1, 5.3.1] we deduce the existence of a cocharacter $\mu_V^{\prime}$ of the subgroup $G_V$ of $GL(M_V^\prime)$, producing similarly a direct sum decomposition $M_V^{\prime}=F^1_{\prime V}\oplus F^0_{\prime V}$, with $F^1_{\prime V}:=F^1_{\prime}\cap M_V^{\prime}$. 
\finishproclaim

\Proclaim{4.4.5. Local deformations.} \rm
Repeating the arguments of [Va1, 5.3.2-12], we deduce the existence of a principally polarized abelian scheme $(A_{\Rtil e}^{\prime},p_{A_{\Rtil e}^{\prime}})$ over $\Spec(\Rtil e)$ such that the following properties hold.

\medskip
\item{{\bf 1)}}
Its principally quasi-polarized filtered $F$-crystal over $\Rtil e/p\Rtil e$ is 
$$(M^{\prime}\otimes_{Re}\Rtil e,F^1_{\Rtil e},\Phi_{M^{\prime}}\otimes 1,\nabla^{\prime},p_{A^{\prime}}),$$ 
with $F^1_{\Rtil e}$ a direct summand of $M^{\prime}\otimes_{Re} \Rtil e$ of dimension $g_A$ (we still denote by $\nabla^{\prime}$ its extension to $\Rtil e$).

\item{{\bf 2)}}
$w_{\alpha}^{\prime}\in\scrT(M^{\prime}\otimes_{Re} \Rtil e\fracwithdelims[]1p)$ belongs to the $F^0$-filtration of $\scrT(M^{\prime}\otimes_{Re} \Rtil e\fracwithdelims[]1p)$ defined by $F^1_{\Rtil e}$, is fixed by $\Phi_{M^{\prime}}\otimes 1$, and is annihilated by $\nabla^{\prime}$, $\forall\alpha\in\scrJ$. Moreover $w_{\alpha}^{\prime}\in\scrT(M^{\prime}\otimes_{Re} \Rtil e)$, $\forall\alpha\in\scrJ_0$.

\item{{\bf 3)}}
The pull back of the triple $(A_{\Rtil e}^{\prime},p_{A_{\Rtil e}^{\prime}},(w_{\alpha}^{\prime})_{\alpha\in\scrJ})$ through the $V_0$-morphism which at the level of rings is the $V_0$-epimorphism $\Rtil e\twoheadrightarrow V$ defined (see 2.1) naturally by $\pi_V$, is the triple $(A_V^{\prime},p_{A^{\prime}},(t_{\alpha})_{\alpha\in\scrJ})$.

\item{{\bf 4)}}
The pull back of $(A_{\Rtil e}^{\prime},p_{A_{\Rtil e}^{\prime}})$ through the $V_0$-valued Teichm\"uller lift of $\Spec(\Rtil e)$, is a principally polarized abelian variety $(A^0,p_{A^0})$ over $V_0$, whose special fibre $(A^0_{\dbF},p_{A^0_{\dbF}})$ is the special fibre of $(A_V^{\prime},p_{A^{\prime}})$; its associated principally quasi-polarized filtered $\sigma$-crystal over $\dbF$ is $(M_0^{\prime},F_0^{1\prime},\phi^{\prime},p_{A^0})$, where $F_0^{1\prime}:=F^1_{\Rtil e}\otimes_{\Rtil e} V_0\subset M_0^\prime:=M^\prime\otimes_{Re} V_0$.

\item{{\bf 5)}} 
There is a family $(t^0_{\alpha})_{\alpha\in\scrJ}$ of tensors of the $F^0$-filtration of $\scrT(M_0^{\prime}\fracwithdelims[]1p)$ (defined by $F_0^{1\prime}$), fixed by $\phi^{\prime}$, and such that we get a principally quasi-polarized Shimura filtered $\sigma$-crystal $\scrC_0:=(M_0^{\prime},F^{1\prime}_0,\phi^{\prime},G_{V_0},(t^0_{\alpha})_{\alpha\in\scrJ},p_{A^0})$. It is obtained naturally from the family $(w_{\alpha}^\prime)_{\alpha\in\scrJ}$, via the pull back of 4). We have $t_{\alpha}^0\in\scrT(M_0^{\prime})$, $\forall\alpha\in\scrJ_0$.
\finishproclaim

\smallskip
$A^{\prime}_{\Rtil e}$ is defined by a lift of the cocharacter $\mu^{\prime}_V$ of 4.4.4.2 to a cocharacter $\mu^{\prime}_{\Rtil e}\colon\dbG_m\to G_{\Rtil e}$. Let $\mu^{\prime}_{V_0}\colon\dbG_m\to G_{V_0}$ be obtained via the pull back of 4). We get a direct sum decomposition $M_0^{\prime}=F_0^{1\prime}\oplus F_0^{0\prime}$,  with $\gamma\in\dbG_m(V_0)$ acting through $\mu^{\prime}_{V_0}$ as the multiplication with $\gamma^{-i}$ on $F_0^{i\prime}$, $i=\overline{0,1}$.

We could identify $t_{\alpha}^0$ with $t_{\alpha}$, $\alpha\in\scrJ$; in what follows, not to create any unclarity, this will not be done.

\Proclaim{4.4.5.1. Back to an unramified context.} \rm
As $\Rtil e$ is a projective limit of artinian rings having $\dbF$ as their residue fields, the level-4 symplectic similitude structure of $A_V^\prime$ lifts to a level-4 symplectic similitude structure of $(A^\prime_{\Rtil e},p_{A^\prime_{\Rtil e}})$. So we get naturally a morphism 
$$l_{\Rtil e}:\Spec(\Rtil e)\to\scrA_{g_A,1,4}$$ 
through which (cf. 4.4.5 3)) the morphism $l_V$ (of 4.3.7.1) factors. Let 
$$l_{V_0}^0\colon \Spec(V_0)\to\scrA_{g_A,1,4}$$
be the composite of the pull back of 4.4.5 4) with it. 
\finishproclaim

\Proclaim{4.4.6. The $d$ number.} \rm
Let $d$ be the relative dimension of the quotient $G_{V_0}/P_{V_0}$, where $P_{V_0}$ is the parabolic subgroup of $G_{V_0}$ normalizing the direct summand $F^{1\prime}_0$ of $M_0^{\prime}$. Let $R:=V_0[[x_1,..,x_d]]$, with $x_i$ as independent variables. Let $\Phi_R$ be its Frobenius lift compatible which takes $x_i$ to $x_i^p$, $\forall i\in\{1,...,d\}$. Obviously we have $d\Le\dim_{\dbC}(X_A)$.

\Proclaim{4.4.6.1. Lemma.} 
We have $d\Ge 1$. 
\finishproclaim

\proof
If $v_{\alpha}$ is the $p$-component of the \'etale component of a Hodge cycle of $A_{\dbC}$ (i.e. of $A^{\prime}_K$, cf. 4.3.6.2 and 4.3.4.7) for some $\alpha\in\scrJ$, then $t_{\alpha}^0$ is the de Rham component of the same Hodge cycle (of $A^{\prime}_K$), cf. [Va1, 5.2.10]. So the group $H_{A\overline{K_0}}$ can be canonically viewed as a subgroup of $GL(M_0^{\prime}\otimes_{V_0} \overline{K_0})$. The monomorphisms we get $G_{\overline{K_0}}\hookrightarrow H_{A\overline{K_0}}\hookrightarrow GL(M_0^{\prime}\otimes_{V_0} \overline{K_0})$ are nothing else but the extension to $\overline{K_0}$ of the initial monomorphisms $G_{\dbQ_p}\hookrightarrow H_{A\dbQ_p}\hookrightarrow GL(H^1[{1\over p}])$. This is entirely similar to [Va1, 5.2.17.2] (being a property of Fontaine's comparison theory). If $d=0$, then $G_{\overline{K_0}}^{\der}$ is contained in a parabolic subgroup of $H_{A\overline{K_0}}^{\der}$, different from $H_{A\overline{K_0}}^{\der}$, and so it is included in a Levi subgroup of $H_{A\overline{K_0}}^{\der}$. From the structure of parabolic subgroups of simple factors of $H_{A\overline{K_0}}^{\der}$ (see [Bo2, 14.17-18]) we get easily that the centralizer of $G_{\overline{K_0}}^{\der}$ in $H_{A\overline{K_0}}^{\der}$ is non-trivial. This contradicts 4.2.3. So $d\Ge 1$.

\Proclaim{4.4.7. Theorem.} 
There is a morphism $l_R\colon\Spec(R)\to {\scrA_{g_A,1,4}}$ such that the following hold:

\medskip
{\bf a)} $l_{\Rtil e}$ and $l_V$ factor through it;

\smallskip
{\bf b)} $l_{V_0}^0$ is the composite of the Teichm\"uller lift $\Spec(V_0)\hookrightarrow\Spec(R)$ with $l_R$;

\smallskip
{\bf c)} The principally quasi-polarized $p$-divisible object of $\scrM\scrF_{[0,1]}^\nabla(R)$ defined by the principally polarized abelian scheme $(A_R,p_{A_R})$ over $R$ it gives birth to, is isomorphic to one of the form
$$\scrC_R:=(M_0^\prime\otimes_{V_0} R,F_0^{1\prime}\otimes_{V_0} R,g_{\text{univ}}(\phi^\prime\otimes 1),\nabla^0,p_{A^0}),$$
with $g_{\text{univ}}\in G_{V_0}(R)$ fixing the perfect form $p_{A^0}$ on $M_0^\prime\otimes_{V_0} R$ and such that $\scrC_R$ is a uni plus quasi-versal deformation of $\scrC_0$.
\finishproclaim

\smallskip
\proof
b) and c) are just a variant of [Va1, 5.4.4-8]: we can take $\scrC_R$ to be (see [Va2, 2.2.10]) the Shimura filtered $F$-crystal $(M_0^\prime,F_0^{1\prime},\phi^\prime,G_{V_0},N_{V_0},\tilde f)$, with $N_{V_0}$ a smooth subgroup of $G_{V_0}$ for which we have an open embedding $N_{V_0}\hookrightarrow G_{V_0}/P_{V_0}$. So a) follows from 2.4.1.

\Proclaim{4.4.7.1. Remark.} \rm
The connection $\nabla^0$ on $M_0^\prime\otimes_{V_0} R$ respects the $G_{V_0}$-action in the sense of [Fa3, rm. ii) after th. 10]. So, if $v_{\alpha}$ is the $p$-component of the \'etale component of a Hodge cycle of $A_{\dbC}$ (cf. 4.3.4.7), then $t^0_{\alpha}\in\scrT(M_0^\prime\otimes_{V_0} R\fracwithdelims[]1p)$ is the de Rham component of a Hodge cycle of (the generic fibre of) $A_R$. This is a consequence of [Va1, 5.2.10] (cf. the above proof of 4.4.6.1) and of [Va1, 5.4.4-8]. We conclude (as in [Va1, 5.4.5]):

\medskip\noindent
{\bf Corollary.} {\it The morphism $l_R$ factors through the normalization $\scrN(4)$ of ${\scrA_{g_A,1,4}}_{O_{(v)}}$ in the ring of fractions of $\Sh(H_A,X_A)_E$.}
\finishproclaim

\Proclaim{4.4.7.2. Remark.} \rm
As [Fa3] is stated over perfect fields, in 4.3.7 to 4.4.7 we could have worked as well with the rings of integers of $E_v$ and $E^\prime_{v^\prime}$ instead of $V_0$ and $V$. The only difference: we would get reductive groups over spectra of different integral rings involved, which are not necessarily split (however, their generic fibres are --cf. 4.3.6.2-- inner forms of logical extensions of $G_{\dbQ_p}$).
\finishproclaim

\Proclaim{4.4.8. Corollary.}
Let $T_v$ be the Frobenius torus of the reduction $A_v$ of $A$ w.r.t. $v$. If $A^0_{\dbF}$ (equivalently $A_v$) is an ordinary abelian variety, then we have:

\medskip
{\bf a)} its canonical lift $A^1$ is obtained from $A_R$ by pull back through a $V_0$-valued point $l^1_{V_0}$ of $\Spec(R)$;

\smallskip
{\bf b)} $T_v$ is isomorphic to a torus of $H_A$;

\smallskip
{\bf c)} $T_{v\dbQ_p}$ is isomorphic to a torus of $G_{\dbQ_p}$;

\smallskip
{\bf d)} for any rational prime $q$ we have $G_{\dbQ_q}=\tilde G_{\dbQ_q}$.
\finishproclaim

\proof
a) is a consequence of [Va2, 2.3.17.1-2 and 3.1.0 b)] and of [Va2, 3.11.1 c)], once we remark that the $G_{V_0}$-ordinary type produced by $\scrC_0$ is a usual ordinary type and the degree of definition $d_1$ of a Shimura-canonical lift $\scrC_1$ produced by $(M_0^\prime,\phi^\prime,G_{V_0})$ is $1$ (see [Va2, 3.1.1 and 3.11.3] for definitions). The first stated fact is a consequence of the ordinariness of $A^0_{\dbF}$. The second one is obtained in the same manner as [Va2, 4.6 P2]; it is as well a consequence of [Va2, 3.1.0 b) and c)] once we remark that the canonical split cocharacter of the canonical lift of any ordinary $\sigma$-crystal $(M_0^\prime,g_{V_0}\phi^\prime)$, with $g_{V_0}\in G_{V_0}(V_0)$, factors through $G_{V_0}$.

To check b), we first remark that $A^1$ is the extension to $W(\dbF)$ of an abelian variety $A^1_f$ over the Witt ring of a finite field $k_1$ containing $k(v^\prime)$ (cf. 4.4.7.2 or standard arguments on canonical lifts). From Corollary of 4.4.7.1 and the fact that each endomorphism $\pi_1$ of the special fibre $SF$ of $A^1_f$ lifts (after extension of scalars) to an endomorphism of $A^1$, the Mumford--Tate group $T_{A^1}$ of $A^1$ (or of $A^1_f$) is a torus of $H_A$ (uniquely determined up to $H_A(\dbQ)$-conjugation). So (cf. 4.2.7.1) a sufficiently high, integral power of the Frobenius endomorphism of $SF$ defines a $\dbQ$--valued point of $Z^0(C(T_{A_1}))$, where $C(T_{A_1})$ is the centralizer of $T_{A^1}$ in $GSp(W,\psi)$. So $T_v$ can be naturally identified with a subtorus of $Z^0(C(T_{A_1}))$. This identification is compatible with the natural identification (obtained via crystalline realizations) of $Z^0(C(T_{A_1}))_{K_0}$ with a torus of $GL(M_0^\prime[{1\over p}])$. So, as $A^1$ is a canonical lift, the generic fibre of the canonical split cocharacter of the filtered $\sigma$-crystal $(M_0^\prime,F_{\text{can}}^1,\phi^\prime)$ of $A^1$ is a cocharacter of ${T_{A^1}}_{K_0}$. It is as well the Newton cocharacter (see [Pi, 3.4]) of ${T_v}_{K_0}$ and the Hodge cocharacter of ${T_{A^1}}_{K_0}$. From Serre's result (see [Pi, 3.5]) and the fact that $T_{A^1}$ is a Mumford--Tate group, we get: as subtori of $Z^0(C(T_{A_1}))$, we have $T_v=T_{A^1}$.

Among many possible proofs of c) we chose the one which (we think) is the most instructive and useful.

\medskip\noindent
{\bf Claim.} {\it There is a canonical isomorphism 
$$ID:H^1_{\text{\'et}}(A^1_{\overline{B(\dbF)}},\dbZ_p)\tilde\to H^1_\prime$$
which takes the tensor $v_{\alpha}^1\in\scrT(H^1_{\text{\'et}}(A^1_{\overline{B(\dbF)}},\dbQ_p))$ obtained from $t^0_{\alpha}$ via Fontaine's comparison theory to $v_{\alpha}$, $\forall\alpha\in\scrJ_0$.} 

\medskip
To check this, we shift from the local context of $R$ to a global context as used in [Va2, 3.6 and 3.15.1]. Let $R^{\text{al}}:=V_0[x_1,...,x_d]$ be equipped with the Frobenius lift $Fr$ which takes each $x_i+1$ into $(x_i+1)^p$. The $p$-adic completion of each formally \'etale $R^{\text{al}}$-algebra introduced below is equipped with the Frobenius lift naturally induced by $Fr$. We consider the localization $R_1$ of $R^{\text{al}}$ in which all $x_i+1$'s are invertible and which has the property that the Frobenius lift $FL_{y^1}$ of the completion of $\Spec(R^{\text{al}})$ in an $\dbF$-valued point $y^1$ of it is of essentially multiplicative type (in the sense of [Va2, 3.6.18.0.1.1]; i.e. the rank, as used in [Va2, 3.6.18.1], of $FL_{y^1}$ is $d$), iff $y^1$ factors through $\Spec(R_1)$. Its existence is implied by [Va2, Claim of 3.6.18.9] and by the fact that $R^{\text{al}}$ is a UFD.  

Let $h_{\text{versal}}\in G(R_1)$ be such that it fixes the perfect form $p_{A^0}$ on $M_0^\prime\otimes_{V_0} R_1$, modulo the ideal $(x_1,...,x_d)$ of $R_1$ is the identity element of $G_{V_0}(V_0)$ and the morphism $\Spec(R_1)\to G_{V_0}/P_{V_0}$ it defines naturally is \'etale. We consider the $p$-divisible object with tensors
$$\scrC_{R_1}:=(M_0^\prime\otimes_{V_0} R_1^\wedge,F_0^{1\prime}\otimes_{V_0} R_1^\wedge,h_{\text{versal}}(\phi^\prime\otimes 1),(t_{\alpha}^0)_{\alpha\in\scrJ})$$
of $\scrM\scrF_{[0,1]}(R_1)$. Starting from it, we construct (as in [Va2, 3.6.1.3 and Theorem 2 of 3.15.1]) a uni plus quasi-versal, global deformation $(\scrC_{R_4},\nabla_{R_4})$ of $(M_0^\prime,F^1_{\text{can}},\phi^\prime,G_{V_0},(t_{\alpha}^0)_{\alpha\in\scrJ})$ over the $p$-adic completion $\Spec(R_4^\wedge)$ of an $\dbN$-pro-\'etale cover $\Spec(R_4)$ of the $p$-adic completion $\Spec(R_3^\wedge)$ of an \'etale scheme $\Spec(R_3)$ over a suitable affine, open, dense subscheme $\Spec(R_2)$ of $\Spec(R_1)$. We take $\Spec(R_2)$ to be an arbitrary affine, open, dense subscheme of $\Spec(R_1)$ for which the following two conditions hold:

\medskip
-- its special fibre is non-empty and by pulling back $\scrC_{R_1}$ through $\dbF$-valued points of $\Spec(R_1)$ factoring through it, we get Shimura filtered $\sigma$-crystals which, when viewed just as $\sigma$-crystals, are ordinary;

\smallskip
-- there is an $R_3$ as mentioned such that the pull back of $\scrC_{R_1}/p\scrC_{R_1}$ to $\Spec(R_3)$ has a connection $\nabla[p]$ on it which is uni plus quasi-versal, $\Spec(R_3/pR_3)$ is connected and non-empty, and the $V_0$-valued point of $\Spec(R_1)$ which at the level of rings takes each $x_i$ to $0$ factors through $\Spec(R_3)$. 

\medskip
Let $z_3:\Spec(V_0)\to\Spec(R_3)$ be this last factorization.
We apply [Va2, 3.6.18.6 a) and Corollary of 3.4.18.9] to the pull back $\scrC_{R_3}$ of $\scrC_{R_1}$ to $\Spec(R_3)$.
 We take $\Spf(R_4)$ to be any connected component of the moduli space of connections on (pull backs of) $\scrC_{R_3}$ (through formally \'etale morphisms) with the property that the connection $\nabla_{R_4}$ on the pull back $\scrC_{R_4}$ of $\scrC_{R_3}$ to $\Spec(R_4)$ we get is such that mod $p$ is the pull back of $\nabla[p]$. From [Fa2, 7.1] and [Va2, Fact of 3.6.19] we deduce that the pair $(\scrC_{R_4},\nabla_{R_4})$ is associated to a Shimura $p$-divisible group $(D_{R_4},(t_{\alpha}^0)_{\alpha\in\scrJ})$ over (warning) $R_4$. Let $z_4$ be a $V_0$-valued point $\Spec(R_4)$ lifting $z_3$. Let $R_5$ be the completion of $R_4$ in $z_4$. We still denote by $z_4$ its factorization through $\Spec(R_5)$.

From [Va2, 2.2.21 UP] we deduce the existence of a unique $V_0$-isomorphism
$$ISO:\Spec(R)\tilde\to\Spec(R_5)$$
such that $ISO\circ l^0_{V_0}=z_4$ and the pull back of $(\scrC_{R_4},\nabla)$ via it is isomorphic to $(\scrC_R,(t_{\alpha}^0)_{\alpha\in\scrJ})$ in such a way that in $l^0_{V_0}$ it is the natural identification of $(M_0^\prime,F_0^{1\prime},\phi,(t_{\alpha}^0)_{\alpha\in\scrJ})$. We also get that $ISO^*(D_{R_5})$ is the $p$-divisible group of $A_R$. Let $m_4:\Spec(V)\to\Spec(R_4)$ be defined naturally via $ISO\circ l_V$.

We apply [Fa2, 2.6] in the context of (truncations modulo positive, integral powers of $p$ of) $(\scrC_{R_4},\nabla_{R_4})$; we have the following 3 steps.

\medskip
{\bf i)} Loc. cit. remains valid in the context of the $\dbN$-pro-\'etale cover $\Spec(R_4)\to\Spec(R_3^\wedge)$. Argument: it involves anyway the passage from $R_3$ to the maximal finite, integral, normal extension $\bar R_3$ of $R_3^\wedge$ which is \'etale over $R_3^\wedge[{1\over p}]$; so $\bar R_3$ is naturally an $R_4$-algebra and so the Fontaine's rings (defined as in [Fa2, p. 36] or [Fa3, ch. 4]) $B^+(R_3)$ and $B^+(R_3^\prime)$ are canonically isomorphic, for any connected, \'etale cover $\Spec(R_3^\prime)$ of $\Spec(R_3^\prime)$ through which $\Spec(R_4)$ factors. As $Fr$ extends to a Frobenius of $V_0[x_1,...,x_d][{1\over {\prod_i=1}^d (x_i+1)}]$ which after inverting $p$ becomes \'etale, we get that $Fr$ extends naturally to a Frobenius lift of any such $R_3^\prime$ (and so to $R_4$ itself) which after inverting $p$ becomes \'etale. Working directly with $\dbZ_p$-coefficients, loc. cit. implies the existence of a filtered $B^+(R_3)$-monomorphism
$$M_0^\prime\otimes_{V_0} B^+(R_3)\hookrightarrow \scrL\otimes_{\dbZ_p} B^+(R_3),\leqno (3)$$
with $\scrL$ as a free $\dbZ_p$-module of rank $2g_A$ on which $\Gal(\bar R_3[{1\over p}]/R_4[{1\over p}])$-acts. $\scrL$ is the dual of the Tate module of $D_{R_4}$ (cf. proof of [Fa2, 7.1]). $B^+(R_3)$ can be viewed as an $R_4$-algebra in a way that respects Frobenius lifts, cf. [Fa2, p. 36]. So, as $M_0^\prime\otimes_{V_0} R_4^\wedge=H^1_{dR}(D_{R_4^\wedge})$, the left hand side of (3) can be identified with $H^1_{dR}(D_{B^+(R_3)})$. It has a canonical $\Gal(\bar R_3[{1\over p}]/R_4[{1\over p}])$-action (different from the natural on the second factor of $M_0^\prime\otimes_{V_0} B^+(R_3)$), cf. [Fa2, p. 36-37]. Also the right hand side of (3) is equipped with the tensor product $\Gal(\bar R_3[{1\over p}]/R_4[{1\over p}])$-action. (3) preserves these two $\Gal(\bar R_3[{1\over p}]/R_4[{1\over p}])$-actions and its cokernel is annihilated by an element $\beta_4\in B^+(R_3)$ having a similar meaning as $\beta$ of 4.4.1, cf. [Fa2, h) of p. 42-43]. 

\smallskip
{\bf ii)} The condition $p\Ge N_2(A)$ allows us to appeal to [Fa2, 2.6]  in the context of the mentioned (truncations of) crystalline tensors, cf. 4.3.4.6.1 2). So to $t_{\alpha}^0$ it corresponds naturally an element $v_{\alpha}^4\in\scrT(\scrL)$ fixed by $\Gal(\bar R_3[{1\over p}]/R_4[{1\over p}])$, $\forall\alpha\in\scrJ_0$. 

\smallskip
{\bf iii)} Loc. cit. is compatible with pull backs via $V$- and $V_0$-valued points. For instance, we refer to $m_4$ and we view $V$ as an $R_4$-algebra. Let $\bar V$ be the normalization of $V$ in $\overline{K}$; we also view it as an $R_4$-algebra. We choose arbitrarily an $R_4$-homomorphism $\bar m_4:\bar R_3\to\bar V$. To it, it is canonically associated a filtered $V_0$-homomorphism 
$$m_4^+:B^+(R_3)\to B^+(V)$$ 
compatible with the PD-structures and with the Galois action of the subgroup $SG$ of $\Gal(\bar R_3[{1\over p}]/R_4[{1\over p}])$ which stabilizes the prime ideal of $\bar R_3[{1\over p}]$ defining the $\overline{K}$-valued point of $\Spec(\bar R_3[{1\over p}])$ defined by $\bar m_4$; $SG$ is naturally identified with $\Gal(K)$. The scalar tensorization of (3) through $m_4^+$ is nothing else but the filtered $B^+(V)$-monomorphism of Fontaine's integral theory of [Fa3, th. 7] associated to $D_V:=m_4^*(D_{R_4})$; this is a consequence of the fact that (cf. i)) we can identify $M_0^\prime\otimes_{V_0} B^+(V)$ with $H^1_{dR}(D_{B^+(V)})=M^\prime\otimes_{Re} B^+(V)$ and that $m_4^+(\beta_4)=\beta$. This equality is a consequence of the fact that the definitions of $\beta_4$ and $\beta$ involve just Teichm\"uller and logarithmic maps. As $D_V$ is nothing else but the $p$-divisible group of $A_V^\prime$, this scalar tensorization is $\rho^\prime $ of 4.4.2. So the pair $(\scrL,(v_{\alpha}^4)_{\alpha\in\scrJ_0})$ is canonically identified with $(H^1_\prime,(v_{\alpha})_{\alpha\in\scrJ_0})$. 
Repeating the arguments for the $V_0$-valued of $\Spec(R_4)$ defined via $ISO\circ l^1_{V_0}$, the pair $(\scrL,(v_{\alpha}^4)_{\alpha\in\scrJ_0})$ is as well canonically identified with $(H^1_{\text{\'et}}(A^1_{\overline{B(\dbF)}},\dbZ_p),(v_{\alpha}^1)_{\alpha\in\scrJ_0})$. So the Claim follows.

\medskip
Under ID, the connected component of the origin of the algebraic envelope of the $p$-adic representation
$$\Gal(B(k_1))\to GL(H^1_\prime)(\dbQ_p)=GL(H^1_{\text{\'et}}(A^1_{\overline{B(\dbF)}},\dbZ_p))(\dbQ_p)$$ 
attached to $A_f^1$ is (cf. 4.3.4.7 and 4.3.4.6.1 2)) a torus $ST$ of $G_{\dbQ_p}$; warning: as above we worked with $\scrJ_0$, here we do use ([Bo2, 7.1] and) the fact that $\pi(G_{\dbQ_p})=p^{-s(p)}v_{\alpha^0}$. $ST$ is a subtorus of $T_{A^1\dbQ_p}$ (for instance, cf. 4.2.2). Moreover, $T_{v\dbQ_p}$ is a subtorus of $ST$. So by reasons of dimensions, we have $ST=T_{v\dbQ_p}=T_{A^1\dbQ_p}$. So c) follows. 

As $H_A$ is a Mumford--Tate group, the canonical homomorphism $T_{A^1}\to H_A^{\ab}$ is an epimorphism. So d) follows from the equality $T_v=T_{A^1}$ (for $q\neq p$, cf. also 4.2.7.2).
This ends the proof.

\Proclaim{4.4.8.1. Remark.} \rm
It is worth pointing out that the two facts of the first paragraph of the proof of 4.4.8 a) are a consequence of the splitness of $G_{\dbQ_p}$ and so hold without assuming that $A^0_{\dbF}$ is ordinary. To argue this, let us point out that $d_1$ can be computed directly in the context of the \'etale cohomology with $\dbQ_p$ coefficients, cf. [Va2, 4.1 and 4.2.1 a)]. Loc. cit. is worked in the context of a SHS but its essence holds in the abstract context of Shimura $\sigma$-crystals. First, global deformations as in the proof of 4.4.8 c) reduce the situation to Shimura-canonical lifts (the way we presented this proof, here we just need $p\Ge 5$ so that we can use $\pi(G_{\dbQ_p})$; however, see 4.4.8.3 below for $p=3$). Second, for them the arguments of [Va2, 4.2.3-7] apply entirely (see also [Va2, c) of 4.4.1 3)]). So if $G_{\dbQ_p}$ is split, $d_1$ is automatically $1$ (cf. [Va2, 4.6 P1-2, 4.6.1 5) and 4.6.3 A]). Moreover, loc. cit. implies as well that $\scrC_1$ is, after ignoring the extra Shimura structure, a usual canonical lift of an ordinary $\sigma$-crystal.
\finishproclaim

\Proclaim{4.4.8.2. Exercise.} \rm
{\bf a)} Use [Va2, 4.8.2 and 4.8.3 e)] to give a new proof of 4.4.8 c).

{\bf b)} Show directly 4.4.8 d) and then use 4.3.4.8 b) to give a new proof of the non-canonical form of the Claim of 4.4.8. Hint: use Lang's theorem and the form of 4.4.7 a) which keeps track of de Rham tensors. Warning: this new proof does require $p\Ge N_2(A)$.
\finishproclaim

\Proclaim{4.4.8.3. Remark.} \rm
[Va2, 2.2.4 H] points out that in fact the condition $p\Ge N_2(A)$ is not needed in 4.4.8 ii): in the proof of 4.4.8 c) we can work with $\scrJ$ instead of $\scrJ_0$ provided $p$ is odd.
\finishproclaim

\Proclaim{4.4.9. Corollary.} 
We do not assume anymore that $A^0_{\dbF}$ is ordinary.
We have: the restriction of $A_R$ to the generic point of $\Spec(R/pR)$ is an ordinary abelian variety. 
\finishproclaim

\proof
This is a consequence of [Va2, 3.12.1] and of 4.4.8.1. 

\Proclaim{4.4.10. Remark.} \rm
If $SL$ is the set of slopes of an isocrystal $(N,\phi_N)$ over a perfect field $k$, if $(N,\phi_N)=\oplus_{\alpha\in SL} (N_{\alpha},\phi_N)$ is the slope decomposition of it (see [Va2, Lemma of 2.2.3 3)]), and if $n\in\dbN$ is such that $n\alpha\in\dbN$, $\forall\alpha\in SL$, then the Newton quasi-cocharacter of $(N,\phi_N)$ is $1\over n$ times the cocharacter of $GL(N)$ through which $p\in\dbG_m(B(k))$ acts on $N(\alpha)$ as the multiplication with $p^{n\alpha}$. In [Va2, 2.2.24] the language of semisimple elements is used. 

We recall that $A^\prime_{V}$ is defined over the ring of integers of $E_{v^\prime}$. Based on this and on Katz--Messing result recalled in [Pi, 3.10] (see also [Mi2, 2.13]), we can identify $T_{vK_0}$ with a torus of the subgroup $G_{V_0}$ of $GL(M_0^\prime[{1\over p}])$. Under this identification, the Newton quasi-cocharacter of $(M_0^\prime,\phi^\prime)$ factors through $T_{vK_0}$ and coincides with the Newton quasi-cocharacter of $T_{vK_0}$ as defined in [Pi, 3.4].

Combining 4.4.9 with a result of Atiyah and Bott (see [Pi, 1.3] and [RR, 2.2] for variants of it) we regain (in a totally different manner) a form of 4.2.14.1 b), strong enough to be a substitute of it in getting 4.2.14.1 (and so 4.2.14).

Moreover, [Va2, 1.12.1 A) and B)] can be used as a refined substitute to the mentioned result of Atiyah and Bott. Loc. cit. implies directly (cf. also [Va2, 2.2.24.1, 3.9.7.2 and 4.12.12.6.3-4]): 

\medskip\noindent
{\bf NP.} {\it The Newton polygon quasi-cocharacter $\mu_{NP}^\prime$ of $(M_0^\prime,\phi^\prime)$ is $G_{V_0}(K_0)$-conjugate to the one $\mu_{NP}$ of the isocrystal of a cyclic diagonalizable Shimura $\sigma$-crystal of the form $(M_0^\prime,g_0\phi^\prime,T_{V_0})$, with $g_0\in G_{V_0}(V_0)$ and with $T_{V_0}$ a maximal torus of $G_{V_0}$.}

\medskip
Here cyclic diagonalizability is in the sense of [Va2, 2.2.22 1)]. Starting from the unique lift of $(M_0^\prime,g_0\phi^\prime,T_{V_0})$ to a Shimura filtered $\sigma$-crystal, we get (as in [Va2, 2.2.9 8)]) a $\dbZ_p$-structure $(M_{0\dbZ_p}^\prime,T_{\dbZ_p})$ of $(M_0^\prime,T_{V_0})$. Based on it, as in [Va2, Fact 2 of 4.1.1.1 and 4.1.2] we get:

\medskip\noindent
{\bf CONJ.} {\it A positive, integral multiple of $\mu_{NP}$ is a cocharacter $\tilde\mu_{NP}$ of $T_{K_0}$ which is a (very precise) sum of cocharacters of $T_{K_0}$ which under the action of suitable powers of $\sigma$ become $G_{V_0}(K_0)$-conjugate to $\mu_{K_0}^\prime$ (and so, in the language of [Pi], are Hodge cocharacters of $G_{K_0}$).}

\medskip
 For our present context of Shimura $\sigma$-crystals, CONJ is a much more precise form of [Pi, 2.3] which generalizes as well [Pi, 2.11].  
\finishproclaim

\Proclaim{4.4.10.1.} \rm
We assume the ranks of $T_v$ and $G_{\dbQ_p}$ are the same. As any two maximal tori of the centralizer $C_{\overline{K_0}}$ of ${\mu^\prime_{NP}}_{\overline{K_0}}$ in $G_{\overline{K_0}}$ are $C_{\overline{K_0}}(\overline{K_0})$-conjugate, from NP and CONJ we get:

\medskip\noindent
{\bf Corollary.}
{\it A positive, integral multiple of ${\mu^\prime_{NP}}_{\overline{K_0}}$ is a cocharacter of $T_{v\overline{K_0}}$ which is a sum of cocharacters of $T_{v\overline{K_0}}$ which under the action of suitable elements of $\Gal(K_0)$ become $G_{V_0}(\overline{K_0})$-conjugate to $\mu_{\overline{K_0}}^\prime$.}
\finishproclaim

\Proclaim{4.4.11. An approach.} \rm
We consider the set $S_{\text{ord}}$ of primes $w$ of $E$ relatively prime to $N(A)!$ and w.r.t. which $A$ has ordinary good reduction. We assume it has a positive Dirichlet density $\delta$. From 4.4.8 b) we get: $\forall w\in S_{\text{ord}}$, the Frobenius torus $T_w$ of the reduction $A_w$ of $A$ w.r.t. $w$ is a torus of $H_A$. Replacing $S_{\text{ord}}$ by a subset of it $S_{\text{ord}}^0$ which still has Dirichlet density $\delta$, we can assume (cf. [Chi, 3.8]) that the rank of $T_w$ is the same as the rank of $G_{\dbQ_p}$, $\forall w\in S_{\text{ord}}^0$. So if 4.1 does not hold for A, from 4.4.8 c) we get that for any prime $q$ there is a reductive, proper subgroup $G_{\dbQ_q}$ of $H_{A\dbQ_q}$ such that:

\medskip
-- $\forall v\in S_{\text{ord}}^0$, an $H_A(\dbQ_q)$-conjugate of ${T_v}_{\dbQ_q}$ is a maximal torus of $G_{\dbQ_q}$;

\smallskip
-- 4.2.3 (with $p$ replaced by $q$) holds.

\medskip
It should not be too difficult to show (at least in most of the cases) that there is a reductive subgroup $\bar G_{\dbQ}$ of $H_A$ such that:

\medskip
-- for any $v$ in a subset of $S_{\text{ord}}^0$ having Dirichlet density $\delta$, an $H_A(\dbQ)$-conjugate of $T_v$ is a maximal torus of $\bar G_{\dbQ}$;
 
\smallskip
-- the equivalent of 4.2.3 holds for monomorphisms $\bar G_{\dbQ}\hookrightarrow H_A\hookrightarrow GL(W)$; 

\smallskip
-- there are primes $p$ such that $\bar G_{\dbQ_p}$ is $H_A(\dbQ_p)$-conjugate to $G_{\dbQ_p}$. 

\medskip
For this an interpolation similar to the one of [Pi, \S6] is needed. To us, if $H_A^{\ad}$ is $\dbQ$--simple, it would be also very useful to show (if not that $\tilde G_{\dbQ}^{\ad}$ is $\dbQ$--simple) at least that all ranks of the Lie types of its simple factors are the same.

The existence of $\bar G$ ought to lead to a contradiction (perhaps using 5.2.1 below). This is one of the approaches to be taken in [Va7]. Based on [Va1, 1.7.1] in [Va7] we will see that an entirely similar approach is possible, without any assumption on $S_{\text{ord}}$; however, it looks much more difficult to reach a contradiction as above in the most general cases: for instance, there are limitations in the isotypical ideas of [LP1] or of [Pi, p. 212-213].
\finishproclaim

\Proclaim{4.4.12. Remark.} \rm
If $p>>0$, then [LP, 3.2] (resp. its proof) implies that $G_{\dbQ_p}$ is unramified (resp. that in 4.3.3 we can take $H^1_\prime=H^1$). So all the parts of 4.3-4 listed in the end of 4.3.0 remain true if we just assume that $G_{\dbQ_p}$ is unramified and that $p>>0$, without having to pass in 4.3.5 to an isogeny (and so to a potentially ramified context). Warning: this approach is ineffective.

We assume now that $p\Ge N(A)$ and that $G_{\dbQ_p}$ is just unramified. Then the mentioned places still hold provided we assume that 4.3.3 a) holds. Just one modification is required: in the proof of 4.3.4.8, we need to use all tensors of $\End(H^1_\prime)$ fixed by $G_{\dbZ_p}$ in order to conclude that $T^0_{pR}$ is a torus, cf. [Va1, 4.3.13]. In many situations it is easy to see that 4.3.3 a) holds. For instance, if no element of $H^1[{1\over p}]$ is fixed by $G_{\dbQ_p}^{\der}$, if $G_{\dbQ_p}$ splits over an unramified extension of odd degree and if all simple factors of $G_{\dbQ_p}^{\ad}$ are of some $A_n$ Lie type, with $n\in\dbN\setminus\{1\}$, then entirely similar to the proof of [Va1, 3.1.2.1 d)] we get that 4.3.3 a) holds. 

In case 4.3.3 a) does not hold, we can proceed as follows. We take $H^1_\prime$ such that the Zariski closure $G_{\dbZ_p}$ of $G_{\dbQ_p}$ in $GL(H^1_1)$ is reductive, cf. [Ja, 10.4 of Part I]. We apply [Va1, 6.5.1.1 v)] to the $\dbZ_{(p)}$-lattice $L^*\otimes_{\dbZ} \dbZ_{(p)}$  of $W$. [Va1, 6.5.1.1 v)]. We get, eventually after replacing $A$ by $A^2$, we can assume that $A^\prime$ has a principal polarization such that the perfect form $p_{A^\prime}$ on $H^1_\prime$ corresponding to it is fixed by $G_{\dbZ_p}$. Warning: as forms on $H^1[{1\over p}]$, $p_{A^\prime}$ can be different from $p_A$. So we still get the mentioned parts, but without knowing that isogenies of 4.3.5 and 4.3.7.1 respect the principal polarizations.
\finishproclaim

\bigskip
\noindent
{\boldsectionfont \S5. Some conclusions}
\bigskip

In this chapter we include some conclusions of what the techniques of \S3-4 can achieve.

\smallskip
\Proclaim{5.1. Applications to 4.1.} \rm
They are stated in 5.1.1-2 and 5.1.5. For the proof of 5.1.1 (resp. of 5.1.2) see 5.1.3 (resp. 5.1.4). 5.1.2.1-5 contains material leading to a deeper understanding of 5.1.1-2 and of what is still left to be proved. 5.1.6 complements [No, 2.2] while 5.1.7 generalizes 4.2.15.
\finishproclaim

\Proclaim{5.1.0. $N_3(A)$.} \rm
Let $N_3(A):=p(A)$, where the prime $p(A)$ is obtained in the same manner as in 3.3 E), starting from $(H_A^{\ad},X_A^{\ad})$ and a prime $v(A,E)$ of $E$ dividing it. Let $v(A)$ be a prime of $E$ dividing a rational prime $l(A)\Ge p(A)$, with $g.c.d.(p,v(A))=1$ and such that the ranks of $G_{\dbQ_p}$ and $T_{v(A)}$ are the same, cf. end of 4.2.7.1. Let $w(A)$ be the prime of $E(H_A,X_A)$ divided by $v(A)$, cf. 4.3.7.

Most common we do not mention the second injective map 
$$f_A^2:(X_A^2,X_A^2)\to (GSp(W_2,\psi_2))$$ 
(accompanying $F_A$ and) needed to truly define $p(A)$. Warning: we take $f_A^2$ to be the analogue of $\tilde f_2$ of 3.3.3.1 obtained starting from $(H_A^{\ad},X_A^{\ad})$; moreover, we often say that we can assume (via 3.3.3 or 3.3.3.1) that $f_A=f_A^2$, or that $f_A$ factors through $f_A^2$, or that $f_A$ or $H_A$ has such an such properties. Despite all these, strictly speaking we do define $p(A)$ via such a fixed $f_A^2$ and by these assumptions we just mean that we shift from $f_A$ to $f_A^2$ without introducing extra notations. Warning: $f_A^2$ is truly introduced just at the needed moment (i.e. in 5.1.3 S1 below). 

For the meaning of $N(A)$ see 5.3.0. We use the notations of 1.4 (1). Warning: $H_A(i)^{\ad}$'s of 4.3.2.1 are different in general from $H_j$'s of 1.4.
\finishproclaim

\Proclaim{5.1.1. Theorem.} 
If $p\Ge N(A)$, then the equality $\dim_{\dbC}(X_A)=d$ implies that $G_{\dbQ_p}=H_{A\dbQ_p}$.
\finishproclaim

\Proclaim{5.1.2. Theorem.} 
If $p\Ge N(A)$, then the split criterion is true (and so 4.1 is true for $A$, cf. 4.2.9) provided $\forall j\in\scrI$ one of the following conditions holds:

\medskip
{\bf a)} $(H_j,X_j)$ is of non-special $A_n$ type (see definition 5.1.2.1 2) below);

\smallskip
{\bf b)} $(H_j,X_j)$ is of $B_n$ type;

\smallskip
{\bf c)} $(H_j,X_j)$ is of $C_n$ type, with $n$ odd;

\smallskip
{\bf c')} $(H_j,X_j)$ is of $C_{2n}$ type, with $n$ odd and such that $4n$ is not of the form $C_{4m+2}^{2m+1}$, with $m\in\dbN$;

\smallskip
{\bf d)} $(H_j,X_j)$ is of $D_n^{\dbH}$ type, with $n$ odd and such that $2n$ is not of the form $C_{2^{m+1}}^{2^m}$, with $m\in\dbN$;

\smallskip
{\bf d')} $(H_j,X_j)$ is of non-inner $D_{2n}^{\dbH}$ type, with $n\in\{q,3q|q\,\text{a prime}\}\setminus\{9\}$;  

\smallskip
{\bf e)} $(H_j,X_j)$ is of $D_n^{\dbR}$ type, $n\in\dbN$, $n\Ge 4$.
\finishproclaim

\Proclaim{5.1.2.1. Definitions.} \rm
{\bf 1)} Let $f_1\colon (H_1,X_1)\hookrightarrow (GSp(W_1,\psi_1),S_1)$ be an injective map. A torus $T_1$ of $H_1$ is called $Sh$-good for $f_1$ if it is the smallest subgroup of $H_1$ through which a set of cocharacters $\dbG_m\to H_{1\dbC}$ which are $H_1(\dbC)$-conjugate to the cocharacters $\mu_{x_1}^*$ defined in [Va2, 2.3.1] factor. Here $x_1\in X_1$.

{\bf 2)} Let $(H_0,X_0)$ be a simple, adjoint Shimura pair of $A_n$ type, $n\in\dbN$. We say it is of non-special $A_n$ type if there is an injective map $f_1\colon (H_1,X_1)\hookrightarrow (GSp(W_1,\psi_1),S_1)$ such that:

\medskip
\item{{\bf a)}} $(H_1^{\ad},X_1^{\ad})=(H_0,X_0)$;

\smallskip
\item{{\bf b)}} for each $N\in\dbN$, there is a prime $q\Ge N$ such that for every torus $T_1$ of $H_1$ which is $Sh$-good for $f_1$, any reductive subgroup $G_1$ of $H_{1\dbQ_q}$ containing $T_{1\dbQ_q}$ as a maximal torus, whose adjoint has all simple factors of classical Lie type and for which we have $\End(W_1\otimes_{\dbQ} \dbQ_q)^{G_1(\dbQ_q)}=\End(W_1\otimes_{\dbQ} \dbQ_q)^{H_{1\dbQ_q}(\dbQ_q)}$, has the property that $G_1^{\ad}$ has simple factors of $A_n$ Lie type. 
\finishproclaim

\Proclaim{5.1.2.2. Numerical analysis.} \rm
Let $(H_0,X_0)$ be as in 5.1.2.1 2). We now list numerical conditions needed to be satisfied in order $(H_0,X_0)$ to be or or not of non-special $A_n$ type. There is nothing deep in obtaining them, so we just state the results. The results are based on the standard table (see [Pi, 4.2]) of minimal weights of split, simple Lie algebras of classical Lie type. When below numerical conditions are not satisfied, in order to do get examples of non-special $A_n$ types, we just need to take $f_1$ of 5.1.2.1 2) such that 3.3.3 1) and 2) hold for it.

We write
$$H_{0\dbR}=\prod_{j\in I_1} SU(a_j,b_j)_{\dbR}^{\ad}\times\prod_{j\in I_0} SU(n+1,0)_{\dbR}^{\ad},$$ 
with $a_j,b_j\in\dbN$, $a_j\Le b_j$, $a_j+b_j=n+1$, $\forall j\in I_1$. Here $I_0$ and $I_1$ are finite sets, keeping track of the non-compact and respectively compact factors of $H_{0\dbR}$. $I_1$ is not empty, cf. axiom [Va1, SV3]. For $j\in I_1$, let 
$$c_j:={a_j\over b_j}.$$
Let 
$$C:=\{c_j|j\in I_1\}\subset\{{1\over n},{2\over {n-1}},...,{[{{n+1}\over 2}]\over {n+1-[{{n+1}\over 2}]}}\}.$$
So $C$ has at least one element and at most $[{{n+1}\over 2}]$ elements.
Let $c$ (resp. $d$) be the minimum (resp. the maximum) of the numbers in $C$. If $1\in C$, then $n$ is odd. In what follows we consider pairs $(r,s)\in\dbN\times\dbN$ such that $2s-1\Le r$. For such a pair let 
$$c_{(r,s)}:={{C_r^{s-1}}\over {C_r^s}},$$
and 
$$e_{(r,s)}:=C_{r+1}^s.$$

We first point out one elementary thing. If $(H_0,X_0)$ is not of non-special $A_n$ type and if $d<1$, then the ``smallest property" of $T_1$ forces any subgroup $G_1$ of $H_{1\dbQ_q}$ (with $q$ a prime and with $T_1$ and $G_1$ as in 5.1.2.1 2)) to have the property:

\medskip\noindent
{\bf AT)} {\it All simple factors of $\Lie(G_1^{\ad})$ are of (potentially different) $A_m$ Lie types.}

\medskip
This can be read out from [Pi, Table 4.2] (see its list of multiplicities for the $B_n$, $C_n$ and $D_n$ Lie types) or from [Sa] (cf. also 3.3.3 1) and the assumptions of b) of 5.1.2.1 2)).

\medskip
{\bf Case 1: $c=d<1$.} If $(H_0,X_0)$ is not of non-special type, then (cf. AT)) for any $j\in I_1$ there is $t_j\in\dbN$ and pairs $(r_{j,m},s_{j,m})$, $m=\overline{1,t_j}$, such that the following hold:

\medskip

{\bf 1)} $\sum_{m=1}^{t_j} r_{j,m}$ does not depend on $j\in I_1$;

\smallskip
{\bf 2)} $c_{(r_{j,m},s_{j,m})}=c$, $\forall m\in\{1,...,t_j\}$;

\smallskip
{\bf 3)} $e_{(r_{j,m},s_{j,m})}a_j=C_{r_{j,m}}^{s_{j,m}-1}\times\prod_{l=1}^{t_j} e_{(r_{j,l},s_{j,l})}$, $\forall m\in\{1,...,t_j\}$;

\smallskip
{\bf 4)} $e_{(r_{j,m},s_{j,m})}b_j=C_{r_{j,m}}^{s_{j,m}}\times\prod_{l=1}^{t_j} e_{(r_{j,l},s_{j,l})}$, $\forall m\in\{1,...,t_j\}$;

\smallskip
{\bf 5)} If $t_j=1$, then $s_{j,1}\Ge 2$.

\medskip
{\bf Case 2: $c<d<1$.} 
If $(H_0,X_0)$ is not of non-special type, then (cf. AT)) for any $j\in I_1$ there are numbers $n_{j,m}\in\dbN$, with $m$ running as in Case 1, such that:

\medskip
{\bf 6)} the sum $\sum_{m=1}^{t_j} n_{j,m}$ does not depend on $j\in I_1$;

\smallskip
{\bf 7)} if $t_j=1$, then $n_{j,1}<n$;

\smallskip
{\bf 8)} there is a map $f(j)\colon\{1,...,t_j\}\to\{0,1\}$ having the following properties 

\medskip
-- if $f(j)(m)=0$, then we can write $n_{j,m}+1=a_{j,m}+b_{j,m}$, with $a_{j,m},b_{j,m}\in\dbN$ such that $c_{j,m}:={a_{j,m}\over b_{j,m}}\in C$;

\smallskip
-- if $f(j)(m)=1$, then there is $s_{j,m}\in\dbN$, $2s_{j,m}-1\Le n_{j,m}$, such that $c_{j,m}:=c_{(n_{j,m},s_{j,m})}\in C$;

\smallskip
-- the map $q_j\colon \{1,...,t_j\}\to C$ taking $m$ to $c_{j,m}$, is surjective;

\smallskip
-- $n+1=\prod_{m\in\{1,...,t_j\}, f(j)(m)=0} (n_{j,m}+1)\times\prod_{m\in\{1,...,t_j\}, f(j)(m)=1} e_{(n_{j,m},s_{j,m})}$.
\finishproclaim

\Proclaim{5.1.2.2.1. Remarks.} \rm
{\bf 1)} The conditions 5.1.2.2 1) and 6) just express 3.3.2 (via the fact that $T_1$ of 5.1.2.1 is defined over $\dbQ$). In both Cases of 5.1.2.2, $t_j$ counts the number of simply factors of the image of $G^{\der}_{1\overline{\dbQ_q}}$ in the simple factor $\scrF(j)$ of $H_{0\overline{\dbQ_q}}$ corresponding to $j$. In Case 2, the equation $f(j)(*)=0$ keeps track of which of these factors are embedded in $\scrF(j)$ via SD-representations. 

\smallskip
{\bf 2)} Of course, if $c<d=1$ some similar numerical conditions as in Case 2 of 5.1.2.2 can be obtained. We do not stop to list them here. Just one sample: we similarly get surjective maps $q_j$'s onto $C$. Warning: the property AT) does not a priori hold.
\finishproclaim  

\Proclaim{5.1.2.3. Examples.} \rm
We present some situations to which (the ideas of) 5.1.2.2 apply.

\smallskip
{\bf Example 1.} $(H_0,X_0)$ is of non-special $A_n$ type if there is $j\in I_1$ such that $g.c.d.(a_j,b_j)=1$ (so $c_j<1$ if $n>1$) and the pair $(a_j,b_j)$ is not of the form  $(C_r^{s-1},C_r^s)$, for $(r,s)$ a pair as above, with $s\Ge 2$. In particular, $(H_0,X_0)$ is of non-special type if there is $j\in I_1$ with $a_j=1$ (i.e. if ${1\over n}\in C$).

\smallskip
{\bf Example 2.} We assume that $c<d$ and that there is $j\in I_1$ such that $g.c.d.(a_j,b_j)=1$. Then $(H_0,X_0)$ is of non-special $A_n$ type. 

To see this, we first remark that $g.c.d.(a_j,b_j)=1$ implies that for any prime $q$, under the situation of 5.1.2.1 2), the image of $\Lie(G_{1\overline{\dbQ_q}})$ in the simple factor $\scrL_j$ of $\Lie(H_{1\overline{\dbQ_q}}^{\ad})$ corresponding to this $j$, is a simple Lie algebra $\scrL(0)$. Now, the surjectivity property mentioned in 5.1.2.2 and in 5.1.2.2.1 2) forces $\scrL(0)=\scrL_j$. In the same way we argue Example 1.
 
\smallskip
{\bf Example 3.} We assume that $c<d<1$ and that there is $j\in I_1$ such that $g.c.d.(a_j,b_j)=2$. Then $(H_0,X_0)$ is of non-special $A_n$ type.

\smallskip
{\bf Example 4.} We assume that $c<d<1$ and that there is $j\in I_1$ such that $g.c.d.(a_j,b_j)$ does not have factors which can be written as a product of numbers of the form $\tilde a_j+\tilde b_j$, with $\tilde a_j,\tilde b_j\in\dbN$ such that ${\tilde a_j\over\tilde b_j}\in C$, or of the form $e_{(\tilde r_j,\tilde s_j)}$, with the pair $(\tilde r_j,\tilde s_j)$ subject to above conditions and such that $c_{(\tilde r_j,\tilde s_j)}\in C$. Then $(H_0,X_0)$ is of non-special $A_n$ type, cf. 5.1.2.2 8). This example generalizes the previous one.

\smallskip
{\bf Example 5.} If $c=d=1$, then $(H_0,X_0)$ is of non-special type iff $n=3$ or $4$ does not divide $n+1$. The necessity can be seen by using a central quotient by $\mu_2$ of a product of two split, simply connected semisimple groups of $A_1$ and respectively of $A_m$ Lie type ($m\in 2\dbN+1$), embedded logically (via tensor products of $SD$-representations) in a split, simply connected semisimple group of $A_{2m+1}$ Lie type. 

To see the sufficiency part, we first remark that $c=d=1$ implies that $(H_0,X_0)$ is without involution and so we can use 3.3.3 2); we get: not all of $n_{j,m}$ are $1$. As ${{n+1}\over 4}\not\in\dbN\setminus\{1\}$, the image of $G_{1\overline{\dbQ_q}}$ in a simple factor of $H_{0\overline{\dbQ_q}}$ is a simple, adjoint group $SG$ of some $A_r$ Lie type. $c_{(r,s)}=1$ iff $r=2s-1$; if $r=2s-1$, then the fundamental representation of $\Lie(SG)$ associated to the weight $\overline{w}_s$ is self dual, cf. [Bou2, p. 188-189]. So using once more 3.3.3 2), in order not to get a contradiction with 4.2.3 we must have $r=n$.

\smallskip
{\bf Example 6.} If $c<d=1$ and if $(H_0,X_0)$ is not of non-special type, then for any $j\in I$ such that $c_j<1$, $2$ divides $a_j$ and $b_j$. This is a consequence of the surjectivity part of 5.1.2.2.1 2).

\smallskip
{\bf Example 7.} If $n+1$ is $4$ or a prime, then $(H_0,X_0)$ is of non-special type, regardless of how $C$ is. This is an elementary exercise.

\smallskip 
We do not stop here to restate 5.1.2.2-3 in terms of values of $n$, similar to statements of Examples 5 and 7. Referring to 1.4, i) (resp. iii)) of it is a particular case of Example 2 (resp. 6); moreover,  1.4 ii) (resp. iv)) is just Example 1 (resp. 7).     
\finishproclaim

\Proclaim{5.1.2.4. The simplest cases not settled by 5.1.2.} \rm
Here, for the sake of future references, we include the simplest (in the sense of involving adjoint groups of small rank and fixed parity) cases not settled by 5.1.2. For all of them we assume $H_A^{\ad}$ is a simple $\dbQ$--group but we do not assume any restriction on $p$ or on $G_{\dbQ_p}$. We denote by $IM$ the image of $G_{\overline{\dbQ_p}}$ in an arbitrarily chosen simple factor of $H_{A\overline{\dbQ_p}}^{\ad}$. It is an adjoint group (easy consequence of 4.2.3) but is not necessarily simple. We will use the sign $+$ to denote the Lie type of products of simple, adjoint ${\overline{\dbQ_p}}$-groups. 

\smallskip
{\bf A5.} $H_A^{\ad}$ is of $A_5$ Lie type and $C=\{1,{1\over 2}\}$. $IM$ could theoretically be of $A_1+A_2$ Lie type. So among the $7$ possibilities for the set $C$, only 1 is not settled.

\smallskip
{\bf A7.} $H_A^{\ad}$ is of $A_7$ Lie type, and $C=\{1\}$ or $C=\{1,{1\over 3}\}$. $IM$ could theoretically be of $A_1+A_3$ Lie type. So among the $15$ possibilities for the set $C$, only 2 are not settled.

\smallskip
{\bf A8.} $H_A^{\ad}$ is of $A_8$ Lie type and $C=\{{1\over 2}\}$. $IM$ could theoretically be of $A_2+A_2$ Lie type. So among the $15$ possibilities for the set $C$, only 1 is not settled.

\smallskip
{\bf A9.} $H_A^{\ad}$ is of $A_9$ Lie type. $IM$ could theoretically be of $A_4$ Lie type and then $C=\{{2\over 3}\}$, or $IM$ could theoretically be of $A_1+A_4$ and then either $C=\{1,{1\over 4}\}$ or $C=\{1,{2\over 3}\}$. So among the $31$ possibilities for the set $C$, only 3 are not settled.

\smallskip
{\bf A11.} $H_A^{\ad}$ is of $A_{11}$ Lie type. $IM$ could theoretically be of $A_1+A_5$, $A_1+A_1+A_2$, $A_2+A_3$ or $A_2+C_2$ Lie type. It is worth pointing out that $A_1+A_1+A_2$ and $A_2+C_2$ have the same rank and moreover, in both these situations, we have $C=\{1,{1\over 2}\}$. This is the simplest situation where a $C_m$ Lie type, with $m\Ge 2$, could theoretically show up among the factors of some $G_{\dbQ_p}^{\ad}$'s, for $H_A^{\ad}$ of some $A_n$ Lie type.  

\smallskip
{\bf C4, C8 and C10.} $H_A^{\ad}$ is of $C_4$ (resp. $C_8$ or $C_{10}$) Lie type. $IM$ could theoretically be of $A_1+A_1+A_1$ (resp. of $A_1+A_1+C_2$ or $A_5$) Lie type.

\smallskip
{\bf D4.} $(H_A^{\ad},X_A^{\ad})$ is of non-inner $D_4^{\dbH}$ type. $IM$ could theoretically be of $B_3$ or $B_1+B_2$ Lie type.  

\smallskip
{\bf D6.} $(H_A^{\ad},X_A^{\ad})$ is of non-inner $D_6^{\dbH}$ type. $IM$ could theoretically be of $C_1+C_3$ Lie type.

\smallskip
{\bf D8.} $(H_A^{\ad},X_A^{\ad})$ is of non-inner $D_8^{\dbH}$ type. $IM$ could theoretically be of $B_4$ or $A_1+A_1+A_1+A_1$ Lie type (the $D_5$ Lie type is excluded here by autodualities of [Pi, Table 4.2]).

\smallskip
{\bf D35.} $(H_A^{\ad},X_A^{\ad})$ is of $D_{35}^{\dbH}$ type. $IM$ could theoretically be of $A_7$ Lie type.

\smallskip
Also, it is worth remarking: the Shimura types $A_6$, $A_8$, $A_{10}$ and $A_{12}$ are settled by Example 7 of 5.1.2.3. 
\finishproclaim

\Proclaim{5.1.2.5. Remark.} \rm
Let $F_0:=F(H_0,X_0)$. If by chance we know that there is an infinite number of primes $p$ above which $F_0$ has only one prime, then one can reobtain easily most of the numerical conditions of 5.1.2-3, without being ``bothered" about mentioning $T_1$ in 5.1.2.1 2) or appealing to 3.3.3. But still this is not enough to get them all (we have in mind Example 5 of 5.1.2.3) or to get 5.1.2 a). So even if we have a good understanding of the (assumed to be very simple) arithmetics of $F_0$, still 5.1.3 S1 below (and so implicitly 3.3.3) looks to us unavoidable. However, in practice it is often easy to check that the numbers $t_j$ of the two Cases of 5.1.2.2 do not depend on $j\in I_1$, and accordingly, that the pairs (counted with multiplicities) $(s_{j,m},r_{j,m})$, $m\in\{1,...,t_j\}$, do not depend on $j\in I_1$. Of course, as in 5.1.2.4 A9, we might not have a unique possibility for such a $t_j$.
\finishproclaim

\Proclaim{5.1.3. The proof of 5.1.1.} \rm
It is achieved in three steps. We refer to the equality $d=\dim_{\dbC}(X_A)$ as the DD assumption (property).

Let $G_{V_0}^{\ad}=\prod_{i\in J} G_{iV_0}$ be the decomposition of the adjoint group of $G_{V_0}$ into simple $V_0$-factors. Let $I^{nc}$ be the subset of $J$ formed by those $i\in J$ such that the composition of $\mu^{\prime}_{V_0}$ of 4.4.5 with the natural quotient homomorphism $G_{V_0}\to G_{iV_0}$, is non-zero. Let $I^c:=J\setminus I^{nc}$. We can write (as $G_{\dbQ_p}$ is split)
$$G^{\ad}_{\dbZ_p}=G^{\text{adnc}}_{\dbZ_p}\times G^{\text{adc}}_{\dbZ_p}$$ 
as a product of two adjoint groups over $\dbZ_p$, where 
$$G^{\text{adnc}}_{\dbZ_p}=\prod_{i\in I^{nc}} G_i,$$ 
and $G^{\text{adc}}_{\dbZ_p}=\prod_{i\in I^c} G_i$
are product decompositions in split, adjoint, simple factors such that each $G_{iV_0}$ is the extension of $G_i$ to $V_0$ (so the notation $G_{iV_0}$ is justified). Let $\mu^{nc}_{V_0}$ be the composition of the cocharacter $\mu^{\prime}_{V_0}$ of 4.4.5 with the natural quotient homomorphism $G_{V_0}\twoheadrightarrow G_{V_0}^{\text{adnc}}$. We get (cf. 4.4.6.1) a Shimura group pair $(G^{\text{adnc}}_{V_0},[\mu^{nc}_{V_0}])$ over $V_0$ which is a product of simple, adjoint Shimura group pairs over $V_0$. 

Similarly we can write $H_{AK_0}^{\ad}=H_{AK_0}^{\text{adnc}}\times H_{AK_0}^{\text{adc}}$. Warning: we do not know a priori that the simple factors of $H_{AK_0}^{\text{adnc}}$ or $H_{AK_0}^{\text{adc}}$ are absolutely simple. 

\medskip
{\bf Step 1.} DD implies that $G_{K_0}^{\text{adnc}}=H^{\text{adnc}}_{AK_0}$. Fixing an embedding $i_K\colon K\hookrightarrow\dbC$ extending the embedding $E^\prime\hookrightarrow\dbC$ of 4.3.7, we check this equality over $\dbC$. DD says that the connected, compact Hermitian symmetric spaces defined (see [He]) by $G_{\dbC}^{\text{adnc}}$ and $H^{\text{adnc}}_{A\dbC}$ and their parabolic subgroups normalizing $F_{0}^{1\prime}\otimes_{V_0} \dbC$ have the same dimension; as one is a subspace of the other, they are the same. But these compact symmetric spaces determine these complex adjoint groups (this well known fact can be grasped from [He, ch. VIII, 6.1]: the situation gets easily reduced to the case when our both compact Hermitian spaces are simple). 
This forms the first step.

\medskip
{\bf Step 2.} $G^{\text{adnc}}_{\dbQ_p}$ is a direct factor of $H_{A\dbQ_p}^{\ad}$, cf. Step 1. Step 2 is to get (using this, 3.3 and the tools 4.2.1-7) that $G^{\der}_{\dbQ_p}=H^{\der}_{A\dbQ_p}$. We present three solutions to this statement: the first one, self contained, treating the general case, the second one treating the general case but based on the results stated in [Va1, 1.7], 
while the third one refers to the particular case $\End(A_{\Ebar})=\dbZ$. 
As the results stated in [Va1, 1.7] are proved only in [Va7], the Fact of the second solution below has a star, meaning its proof is not included here.

\smallskip
{\bf S1. Solution 1 (the general case).} As we have $l(A)\Ge N_3(A)$, 3.3.2 applies. Warning: the notations of 3.3.1-2 and of \S4-5 do not correspond: in 3.3.1-2 the group $G$ is defined over $\dbQ$ and is related to a Shimura pair $(G,X)$; moreover, the $p$ and $l$ of 3.3.1-2 are presently $l(A)$ and respectively $p$. From 3.3.2 and Step 1 we get: for any simple factor $\scrF_A^1$ of $H_{A\overline{\dbQ_p}}^{\ad}$, the image $\scrF_B^1$ of $\tilde G_{\overline{\dbQ_p}}$ in it, has the same rank as $\scrF_A^1$. From 4.2.3 we get directly (via 4.2.5.1) that $\scrF_A^1=\scrF_B^1$. 

We now check (based on similar reasons) that a simple factor of $\Lie(G^{\ad}_{\overline{\dbQ_p}})$ can not project epimorphically into two or more simple factors of $\Lie(H^{\der}_{A\overline{\dbQ_p}})$. We can assume that all simple factors of $H_{A\dbC}^{\ad}$ have the same Lie type $\text{LT}$ (cf. 3.3.3 and 3.3.5). 

As an exemplification, we first treat the following particular case. We assume that all factors of $H_{A\dbR}^{\ad}$ are non-compact and that $(H_A^{\ad},X_A^{\ad})$ has no simple factors of $A_n$ type with involution ($n\Ge 2$). So we can assume (cf. 3.3.3 and 3.3.3 3)) that $T^0=\dbG_m$. So if a simple factor of $\Lie(G^{\ad}_{\overline{\dbQ_p}})$ projects epimorphically into two or more simple factors of $\Lie(H^{\der}_{A\overline{\dbQ_p}})$, then the centralizers of $G_{\overline{\dbQ_p}}$ and of $H_{A\overline{\dbQ_p}}$ in $GL(H^1\otimes_{\dbZ_p}\overline{\dbQ_p})$ are distinct and so we reach a contradiction with 4.2.3.

We come back to the general case. We follow closely [De3, 2.3.7-10] and [Va1, 6.5-6]; warning: as in the general case of 3.3.3.1 we work with all simple factors of $(H_A^{\ad},X_A^{\ad})$ at once. In other words, we want to show that something similar to 3.3.3 3) allows us (via 4.2.3) to handle the present situation. We already know that all simple factors of $G^{\ad}_{\dbQ_p}$ are as well of $LT$ Lie type.

We can assume (cf. 3.3.3) that $f_A$ (of 4.3.6.1) factors through an injective map
$$f_A^2\colon (H_A^2,X_A^2)\hookrightarrow (GSp(W,\psi),S),$$
with $H_A^{2\ad}=H^{\ad}_A$ (so $(H_A^{2\ad},X_A^{2\ad})=(H_A^{\ad},X_A^{\ad})$) and with $f_A^2$ as a Hodge quasi product 
$${\times^{\scrH}}_{j\in\scrI} f_j^2$$ 
of injective maps 
$$f_j^2:(H_j^2,X_j^2)\hookrightarrow (GSp(W_j,\psi_j),S_j),$$
each one of them obtained as in 3.1 I) and I'); here $(H_j^2,X_j^2)$ is such that we have an identification $(H_j^{2\ad},X_j^{2\ad})=(H_j,X_j)$, $j\in\scrI$.  
We can assume $f^2_A$ is uniformized in the sense of 3.3.3.1. So there is a totally real number field $F$ such that each $f_j^2$ factors through the composite of the following two injective maps
$$(H_j^3,X_j^3)\hookrightarrow (H_j^1,X_j^1)\operatornamewithlimits{\hookrightarrow}\limits^{f_j^1} (GSp(W_j,\psi_j),S_j)$$
and we have the following properties:

\medskip
{\bf 1)} $H_j^{3\ad}$ is the Weil restriction from $F$ to $\dbQ$ of an absolutely simple, adjoint $F$-group of $LT$ Lie type.  

\smallskip
{\bf 2)} $H_j^{1\ad}$ is the Weil restriction from $F$ to $\dbQ$ of an absolutely simple, adjoint $F$-group of $LT_1$ Lie type. Here $LT_1$ is $LT$ if $LT=A_n$ and is the Lie type of $G(A)$ of 3.1 if $(H_j,X_j)$ has the same type as $(G,X)$ of 3.1. 

\smallskip
{\bf 3)} $f_j^1$ is a PEL type embedding constructed as in 3.1 I') (we recall that, as stated in its beginning paragraph, 3.1 I') handles all $A_n$ types).

\smallskip
{\bf 4)} Two simple factors of $\Lie(H_{A\dbC}^{\der})$ or of $\Lie(H_{A\dbC}^{s\der})$, with $s\in\{2,3\}$, can not map injectively into the same simple factor of $\Lie(H_{A\dbC}^{1\der})$. 

\medskip
If $LT=D_n$, with $n\in\dbN$, $n\Ge 5$, then (due to 4.2.14.1 and table [Pi, 4.2]) we can assume all simple factors of $(H_A^{\ad},X_A^{\ad})$ are either of $D_n^{\dbR}$ or of $D_n^{\dbH}$ type. So if $LT\neq D_4$, then the Lie type of $H^{1\ad}_j$ does not depend on $j$. $f_A^2$ factors through
$$f_A^1:={\times^{\scrH}}_{j\in\scrI} f_j^1:(H_A^1,X_A^1)\hookrightarrow (GSp(W,\psi),S).$$
\indent
Let (as in 3.1 I')) $\scrS$ be the set of extremal points of the Dynkin diagram $\scrD$ of $\Lie(H_{A\dbC}^{1\der})$. $\Gal(\dbQ)$ acts on $\scrS$ and so we identify $\scrS$ with the $\Gal(\dbQ)$-set of $\dbQ$--homomorphisms from $K_\scrS$ to $\dbC$, for $K_{\scrS}$ a suitable product indexed by $j\in\scrJ$ of totally imaginary quadratic extensions of $F$. By enlarging $F$, we can assume $K_{\scrS}$ is a product of $s$ copies of the same totally imaginary quadratic extension $K$ of $F$; here $s$ is the number of simple factors of $H_{A\dbR}^{1\ad}$. $K_{\scrS}$ acts faithfully on $W$ and so we identify $\Res_{K_{\scrS}/\dbQ} \dbG_m$ with a torus of $GL(W)$. We can assume that all homomorphisms $\Res_{\dbC/\dbR} \dbG_m\to H_{A\dbR}^1$ defining elements of $X_A^1$ factor through the extension to $\dbR$ of the subgroup of $GL(W)$ generated by $H_A^{1\der}$, by $Z(GL(W))$ and by the maximal subtorus $T_c$ of $\Res_{K_{\scrS}/\dbQ} \dbG_m$ which over $\dbR$ is compact. This is nothing else but the logical variant of iii) of 3.1 I') for injective maps which are Hodge quasi products. We recall (see 4.0) that $T^0=Z^0(H_A)$. 
We have:

\medskip\noindent
{\bf Fact.} {\it  There are irreducible $T^0_{\overline{\dbQ_p}}$-submodules of $W\otimes_{\dbQ} \overline{\dbQ_p}$ which contain faithful representations of any normal, semisimple subgroup of $H_{A\overline{\dbQ_p}}^{2\der}$ whose adjoint is product of simple factors of $H_{A\overline{\dbQ_p}}^{\ad}$ onto which the same simple factor of $G_{\overline{\dbQ_p}}^{\ad}$ projects naturally.}

\medskip
As in the particular case, we get that this Fact implies $G^{\der}_{\dbQ_p}=H_{A\dbQ_p}^{\der}$. We start its proof by pointing out that $f_A^1$ is a PEL type embedding (cf. 3) and 2.6). We can assume (cf. 3.3.1 Fact') that the Frobenius torus $T_{v(A)}$ is (isomorphic to) a torus of a particular inner form of $H_A^1$. In what follows, as we are interested in ``minimizing" $T^0$, based on 3.3.1 AB, not to complicate notations we assume that this inner form is trivial. Let $T^1_{v(A)}$ be the centralizer of $T_{v(A)}$ in $H_A^1$. We need two Subfacts.

\medskip\noindent
{\bf Subfact 1.} {\it $T^1_{v(A)}$ is a maximal torus of $H_A^1$.}

\medskip
\proof
Based on 1), 2) and 3.1 D), Subfact 1 can be reformulated as: with the notations of 3.1 D), the centralizer in $G(A)_{\dbC}$ of the image through $h_{\dbC}$ of a maximal torus of $\tilde G^{\sc}_{\dbC}$, is a maximal torus of $G(A)_{\dbC}$. If $(G,X)$ is of $A_n$, $C_n$ or $D_n^{\dbH}$ this statement is obvious. For the remaining cases, we can assume we are dealing with the spin representation of a complex, simple Lie algebra of $B_n$ or $D_n$ Lie type. For such a representation, the statement can be easily read out from its description in [Bou2, p. 195-197 and 209]. This proves Subfact 1.

\medskip\noindent
{\bf Subfact 2.} {\it Using the embedding $i_K$ of Step 1, the cocharacter $\mu$ of 4.4.4 and the cocharacter $\mu_{V_0}^\prime$ of 4.4.5, over $\dbC$, become $H_A(\dbC)$-conjugate to the cocharacter $\mu(\dbC)$ of 4.0.}

\medskip
\proof
Let $P_{A\dbC}$ be the parabolic subgroup of $H_{A\dbC}$ normalizing $F^{1,0}$ of 4.0. It is enough to show: if $\mu_1$ is a cocharacter of $P_{A\dbC}$ producing as in 4.0 a direct sum decomposition $F^{1,0}\oplus F^{0,1}_1$, then $\mu_1$ and $\mu(\dbC)$ are $P_{A\dbC}(\dbC)$-conjugate. As any two maximal tori of $P_{A\dbC}$ are $P_{A\dbC}(\dbC)$-conjugate, we can assume that $\mu_1$ and $\mu(\dbC)$ are commuting. This implies $F^{0,1}=F^{0,1}_1$. So $\mu_1=\mu(\dbC)$. This proves Subfact 2.

\medskip
Using the fact that $T_{v(A)}$ is a $\dbQ$--torus of $H_A^1$, we define a $\Gal(\dbQ)$-invariant equivalence relation on $\scrS$ as follows. We say that two elements of $\scrS$ are in relation if the following two conditions hold:

\medskip
{\bf a)} {\it they are nodes of Dynkin diagrams of two simple factors of $\Lie(H_{A\overline{\dbQ_p}}^{2\der})$ on which the same simple factor $\Lie(\scrF)$ of $\Lie(G^{\der}_{\overline{\dbQ_p}})$ maps injectively;}

\smallskip
{\bf b)} {\it w.r.t. the image of $T_{v(A)\overline{\dbQ_p}}^1$ in $H_{A\overline{\dbQ_p}}^{2\der}$ they are either both original nodes or are both ending nodes of Dynkin diagrams (based on a) and Subfact 1 this makes sense; it is 4.2.7.1 which allows us to view $T_{v(A)\dbQ_p}$ as a torus of $G_{\dbQ_p}$).}

\medskip
Let $\scrS^1$ be the set of equivalence classes of this equivalence relation. As $\Gal(\dbQ)$ acts on $\scrS^1$, as above we define $K_{\scrS^1}$ and a subtorus $T_c^1$ of $\Res_{K_{\scrS^1}/\dbQ} \dbG_m$; we identify $K_{\scrS^1}$ with a $\dbQ$--subalgebra of $K_{\scrS}$. Let $F_1:=F\cap K_{\scrS^1}$. $K_{\scrS^1}$ is a a product of at most $s$ field extensions of $F_1$. 

Due to a), b) and the fact that each map $f_j^1$ is obtained following the pattern of $f_1$ of 3.1 I) and I') (mainly we need ii) and iii) of 3.1 I')), the extension to $\dbC$ of $\mu_{V_0}^\prime$ factors through the extension to $\dbC$ of the subgroup of $H_A^1$ generated by $H_A^{1\der}$ and $\Res_{K_{\scrS^1}/\dbQ} \dbG_m$. So (cf. Subfact 2 and the fact that $H_A$ is a Mumford--Tate group) $T^0$ is a subtorus of the torus of $GL(W)$ generated by $Z(GL(W))$ and by $T_c^1$. Recalling (see end of ii) of 3.1 I')) how $\Res_{K_{\scrS}/\dbQ} \dbG_m$ (and so how $\Res_{K_{\scrS^1}/\dbQ} \dbG_m$) acts on $W$, we get the above Fact.   

\smallskip
We conclude $G^{\der}_{\dbQ_p}=H_{A\dbQ_p}^{\der}$.
So $\tilde G_{\dbQ_p}=H_{A\dbQ_p}$. This forms the first solution.

\smallskip
{\bf S2. Solution 2 (the general case).} What follows is included just to emphasize that the machinery of 3.2 (needed for 3.3.2) can be avoided, provided we are eager to assume another (in fact much more technical) machinery of Langlands--Rapoport conjecture.

We apply our results (on the Langlands--Rapoport conjecture) stated in [Va1, 1.7] in the context of the SHS $(f_A,H_{1B}(A,\dbZ_{(l(A))}),w(A))$ (cf. iii) of 3.3 E) and 5.1.0); the set of primes of $E(H_A,X_A)$ referred to in loc. cit. is of positive Dirichlet density (cf. also [Va1, end of 1.7.1]) and so we can assume $w(A)$ is one of such primes.

So, as in the part of 3.3.1 referring to [Zi, 4.4], we can assume (cf. [Mi1, p. 186-191]):

\medskip\noindent
{\bf Fact*.} 
{\it There is a maximal torus of a conjugate of $G_{\dbQ_p}$ by an element of $H_A(\dbQ_p)$ which is the extension to $\dbQ_p$ of a torus $T_A^p$ of $H_A$.} 

\medskip
Step 1 implies: $T_A^p$ projects into a maximal torus of each simple, non-compact factor of $H_{A\dbQ_p}$. As there is no non-trivial, compact factor of $H_{A\dbR}^{\ad}$ defined over $\dbQ$, $T_{A\dbC}^p$ projects into a maximal torus of each simple factor of $H_{A\dbC}^{\ad}$. Applying again 4.2.3 and 4.2.5.1 we get $G^{\der}_{\dbQ_p}=H_{A\dbQ_p}^{\der}$. This forms the second solution.

\smallskip
{\bf S3. Solution 3 (the case $\End(A_{\Ebar})=\dbZ$).} 
In many cases we know a priori that the Fact of S2 is true and so we do not need to refer at all to Langlands--Rapoport conjecture.  
For instance, if $\End(A_{\Ebar})=\dbZ$ and if $H_{A\dbR}^{\ad}$ has no compact factors of $C_r$ Lie type, with $r\in\dbN\setminus\{1,2\}$, then we can assume as well that $A$ has ordinary reduction w.r.t. $w(A)$. So the Fact of S2 follows as well from [Pi, 7.2] and [Va2, 4.6 P2] via (the mixed characteristic $(0,l(A))$ analogue of) 4.4.8 b).

In general, if $\End(A_{\Ebar})=\dbZ$ we can avoid relying on Fact of S2 as follows. [Pi, 5.13 (c) and 5.12] together with Step 1 imply: we can assume all simple factors of $H^{\ad}_{A\overline{\dbQ_p}}$ and of $G^{\ad}_{\dbQ_p}$ have the same Lie type. From this and 4.2.3 we get (by simple reasons of counting dimensions) that the number of simple factors of $H^{\ad}_{A\overline{\dbQ_p}}$ is the same as of $G^{\ad}_{\dbQ_p}$.  We get: $\tilde G_{\dbQ_p}^{\der}=H_{A\dbQ_p}^{\der}$. This forms the third solution. 

\medskip
{\bf Step 3.} The third step is to get $G_{\dbQ_p}=H_{A\dbQ_p}$ knowing that their derived groups are equal. This is achieved using 4.2.7 and [Wi]. We present the details. From 4.2.7.2 we deduce the existence of a subgroup $H_A^\prime$ of $H_A$ such that $G_{\dbQ_p}=H_{A\dbQ_p}^{\prime}$. So $\mu(\dbC)$ factors through $H_{A\dbC}^\prime$, cf. Subfact 2 of S2. So $H_A^\prime=H_A$ and so $G_{\dbQ_p}=H_{A\dbQ_p}$. This forms the third step; it ends the proof of 5.1.1.   
\finishproclaim

\Proclaim{5.1.3.1. Remark.} \rm
If $H_{A\dbR}^{\ad}$ has no compact factors, then we have nothing to prove in Step 2. Similarly, if $\End(A_{\Ebar})=\dbZ$, then Step 3 is nothing else but 4.2.15. 
\finishproclaim

\Proclaim{5.1.3.2. Exercise.} \rm
Using [An, 1.5.1-2] reobtain (in a different way) the result of [An, 1.6.1 3)]) pertaining to Mumford--Tate groups. Hint: use 5.1.2.

\Proclaim{5.1.4. The proof of 5.1.2.} \rm
We follow closely the proof of 5.1.1. It is enough to concentrate on the adjoint groups of $G_{\dbQ_p}$ and of $H_{A\dbQ_p}$, cf. Step 3 of 5.1.3. From 3.3.3, 3.3.5 and 5.1.3 S1 we deduce that we can assume $H_A^{\ad}$ is a simple $\dbQ$--group. Even more, it is enough to show that there is a simple factor $SF$ of $H_{A\overline{\dbQ_p}}^{\ad}$ onto which $G_{\overline{\dbQ_p}}$ projects naturally, cf. 5.1.3 S1.  We show this by treating the cases a) to e) one by one.

\smallskip
{\bf a)} This case is a consequence of 3.3.3 and of [Va2, 2.3.5.1]. In other words, we can assume that the embedding $f_A^2$ is a PEL type embedding and that $f_A$ factors through $f_A^2$. Let $f_1$ be as in 5.1.2.1 2) but for $(H_0,X_0):=(H_A^{\ad},X_A^{\ad})$. We can assume $H_1$ is a Mumford--Tate group. As in 3.3.1 we can assume $G_{\dbQ_p}^{\ad}$ contains a maximal torus which is $H_A^{\ad}(\dbQ_p)$-conjugate to the extension to $\dbQ_p$ of the image $T_A$ of $t_{v(A)}$ in $H_A^{\ad}=H_A^{2\ad}$. Strictly speaking, in 3.3.1 we did not bother to explain in full generality why in fact $G_1=G_1^\prime$; accordingly, we ought to use an inner form of $H_A^{\ad}$ which is still a simple $\dbQ$--group and we ought to twist $H_A^2$ by this inner form. But from the point of view of the Fact below these are irrelevant: in 3.3.1 we have $G_{1\dbQ_p}=G^\prime_{1\dbQ_p}$ and their extensions to $\dbC$ are equipped with the same inner conjugacy classes of cocharacters (cf. 3.3.1 HW and the paragraphs after; the cocharacters are obtained as in Subfact 2 of 5.1.3 S1); so we will not complicate notations by passing to an inner form of $H_A^{\ad}$. Let $T_1$ be the $\dbQ$--torus of $H_1$ generated by $Z^0(H_1)$ and by the $\dbQ$--torus of $H_1^{\der}$ whose image in $H_1^{\ad}=H_A^{\ad}$ is $T_A$. So we can assume that:

\medskip\noindent
{\bf Fact.} {\it $T_1$ is the smallest torus of $H_1$ such that a set of cocharacters of $H_{1\dbC}$ which are $H_1(\dbC)$-conjugate to the cocharacters of $\mu_{x_1}^*$ of [Va2, 2.3.1] (here $x_1\in X_1$) factor through $T_{1\dbC}$}.

\medskip
This is a consequence of Serre's result of [Pi, 3.5] applied to $T_{v(A)}$ and of the analogues of 4.4.10.1 and of Subfact 2 of 5.1.3 S2 in mixed characteristic $(0,l(A))$; warning: 4.2.14.1 or [Pi, 3.16 and 5.10] the way they are stated, do not imply this Fact (however, we could get this Fact tracing through the proofs of loc. cit.). So case a) follows from def. 5.1.2.1 2) (strictly speaking we might have to first replace $p$ by a $q$ as in b) of 5.1.2.1 2) and then to quote 4.2.5).

\smallskip
For cases b) to e) we rely freely on 4.2.14.1. 

\smallskip
{\bf b)} This case is nothing else but [Pi, 4.3].

\smallskip
{\bf c)} This case is the well known case of Serre, cf. [Chi, p. 330] or [Pi, 4.5-6 and Table 4.6].

\smallskip
{\bf c')} This case is a consequence of [Pi, 4.5-6 and Table 4.6].

\smallskip
{\bf d)} This case is just the orthogonal (versus symplectic) variant of c). It can be easily read out from the table [Pi, 4.2].

\smallskip
{\bf d')} To treat this case, we ``move out" from the split criterion. In other words, for this case we rely on 4.2.5 and we work with a prime $v_F$ of $F(H_A^{\ad},X_A^{\ad})$ whose residue field has a prime number $l$ of elements and for which there is a prime $v_I$ of $I(H_A^{\ad},X_A^{\ad})$ dividing it and whose residue field has more than $l$ elements. We assume that $I(H_A^{\ad},X_A^{\ad})$ is unramified above $l$, that $H_A^{\ad}$ is unramfied over $\dbQ_l$, and that $A$ has good reduction w.r.t. some prime of $E$ unramified above $l$. Let $F_l^1:=F(H_A^{\ad},X_A^{\ad})_{v_F}$. We consider a simple factor $G_1$ of $H_{AF_l^1}^{\ad}$ which is not an inner form of a split $F_l^1$-group, cf. the definition of non-inner types (see 2.5) and the choice of $v_F$. It is an unramified, non-split, adjoint group over $F_l^1$ which splits over the unramified quadratic extension $F_l^2$ of $F_l^1$.

There is a cocharacter $\mu_1$ of $G_{1F_l^2}$ such that $(G_1,[\mu_1])$ is a Shimura group pair and its centralizer in $G_{1F_l^2}$ has a derived subgroup of $A_{2n}$ Lie type, cf. 4.2.14.1. So $[\mu_1]$ is not fixed by $\Gal(F_l^2/F_l^1)$. Let $IM$ be the image of $G_{F_l^1}$ in $G_1$. We can assume that a $G_1(\overline{\dbQ_l})$-conjugate of $\mu_{1\overline{\dbQ_l}}$ factors through $IM_{\overline{\dbQ_l}}$.

We assume $G_1$ is different from $IM$. If $n=2q$, with $q$ a prime, then $IM$ is of $C_1+C_q$ Lie type; if $n=6q$, with $q$ a prime different from $3$, then $IM$ is of $C_1+C_{3q}$ or of $C_3+C_q$ Lie type. This follows from 4.2.3 and the table [Pi, 4.2], via 4.2.14.1. So $IM$ is absolutely simple. But any adjoint Shimura group pairs over $F_l^1$ involving absolutely simple groups of $C_n$ Lie types is such the conjugacy class of cocharacters defining it is fixed by $\Gal(F_l^1)$. So $[\mu_1]$ is fixed by $\Gal(F_l^2/F_l^1)$. Contradiction. So $IM=G_1$ and so $G_{\dbQ_l}$ projects onto a simple factor of $H_{A\dbQ_l}^{\ad}$.

\smallskip
{\bf e)} We can assume that 3.3.3 1) and 2) hold (cf. 3.3.3). 4.2.14.1 implies (see [Pi, 4.3]) that the image $\scrF_1$ of $G_{\overline{\dbQ_p}}$ in a simple, non-compact factor $\scrF_0$ of $H_{A\overline{\dbQ_p}}^{\ad}$ is either $\scrF_0$ itself or a product $PB$ of two simple, adjoint groups of $B_{n_1}$ and respectively $B_{n_2}$ Lie type, with $n_1+n_2=n-1$; here the so called $B_0$ Lie type corresponding to the trivial group is allowed. But the representation of $\Lie(\scrF_0)$ on $H^1\otimes_{\dbZ_p} \overline{\dbQ_p}$ involves both the half spin representations and any irreducible subrepresentation of $\Lie(H_{A\overline{\dbQ_p}}^{\der})$ on $H^1\otimes_{\dbZ_p} \overline{\dbQ_p}$ factors through a simple factor of it (i.e. no tensor products are involved). Moreover, 

\medskip\noindent
{\bf SAME.} {\it The representation of $\Lie(PB)$ on a simple, non-trivial $\Lie(\scrF_0)$-submodule $H^1_{\text{half}}$ of $H^1\otimes_{\dbZ_p} \overline{\dbQ_p}$ is the same (it is the tensor product of the spin representations of the two simple factors of $\Lie(PB)$, cf. loc. cit.), regardless of which half spin representation of $\Lie(\scrF_0)$ is involved in defining $H^1_{\text{half}}$.}

\medskip
So if $\scrF_1$ is not $\scrF_0$, we reach a contradiction with 4.2.3 (cf. 3.3.3 1) and 2) and SAME). So $\scrF_0=\scrF_1$. 

This ends the proof of the case e) and so of 5.1.2.
\finishproclaim

\smallskip
In what follows the abelian motives over a field of characteristic 0 are understood in the sense of absolute Hodge cycles, see [DM]. Using realizations in the Betti and \'etale cohomologies, as for abelian varieties we speak about their Mumford--Tate groups and about the Mumford--Tate conjecture for them.

\Proclaim{5.1.5. Corollary (the independence property).}
If the Mumford--Tate conjecture is true for the members of a finite set of abelian motives over a number field $E$, then it is true for their direct sum $DS$.
\finishproclaim

\proof
Let $A_{DS}$ be an abelian variety over $E$ such that $DS$ is an object of the Tannakian category of abelian motives over $E$ generated by its $H^1$-motive $H^1(A_{DS})$ and by Tate-twists (see [DM]). We can assume $A_{DS}$ has a principal polarization $p_{A_{DS}}$ (see 4.3.2). 
We work with a prime $p$ big enough so that we get a SHS $(f_{A_{DS}},L_{(p)},v_{A_{DS}})$, with $f_{A_{DS}}:(H_{A_{DS}},X_{A_{DS}})\hookrightarrow (GSp(W,\psi),S)$ as the injective map defined by $(A_{DS},p_{A_{DS}})$ (see [Va1, 2.12 3)]) and with $L_{(p)}$ as the first group of the Betti homology of $A_{DS}$ with coefficients in $\dbZ_{(p)}$ (the choice of a prime $v_{A_{DS}}$ of $E(H_{A_{DS}},X_{A_{DS}})$ dividing $p$ is irrelevant here). 
[Va2, Fact 1 of 2.3.5.2] speaks about some factors of $(f_{A_{DS}},L_{(p)},v_{A_{DS}})$; they are, in some sense, the ``simple" factors. We have a variant of loc. cit. which pays attention to arbitrary factors of $(H_{A_{DS}}^{\ad},X_{A_{DS}}^{\ad})$. Using such a variant (strictly speaking it can be performed only over $\dbQ$ but for future references it is convenient to point out that we have such a variant over $\dbZ_{(p)}$), we can assume that each simple factor of the adjoint group of the Mumford--Tate group $H_{A_{DS}}$ of $A_{DS}$ is naturally isomorphic to a simple factor of the adjoint group of the Mumford--Tate group of $DS$. Using 5.1.3 S2 and our hypothesis we deduce that the Mumford--Tate conjecture holds for the derived part of $H_{A_{DS}}$. So the Mumford--Tate conjecture holds for $A_{DS}$, cf. Step 3 of 5.1.3. The Corollary follows from this and the fact that the Mumford--Tate groups of $H^1(A_{DS})$ and of its direct sum the motive $T(1)$ of the Tate-twist are the same.

\Proclaim{5.1.5.1. Sample.} \rm
The following sample is just meant to justify the terminology ``independence property". We assume $\End(A_{\dbC})=\dbZ$. So $(H_A^{\ad},X_A^{\ad})$ is simple; we assume it is of $D_{2n+3}^{\dbR}$ type, $n\in\dbN$. We refer to an $f_2$ of 3.3 B) but obtained starting from $(H_A^{\ad},X_A^{\ad})$. We consider an abelian variety $B$ over $E$ which is a pull back (as in 3.3 D)) of an $\scrA_2$ as in 3.3 B); for this we might have to replace $E$ by a finite field extension of it. We assume that the adjoint group of the Mumford--Tate group $MT$ of $A\times_E B$ is $H_A^{\ad}\times H_A^{\ad}$. The connected component of the origin of the algebraic envelope of the image of the $p$-adic Galois representation of $A\times_E B$ could theoretically be a reductive group having $H^{\ad}_{A\dbQ_p}$ as its adjoint. But 5.1.5 says this is not the case (as 4.1 holds for $A\times_E B$, cf. 5.1.2 e)). 
\finishproclaim

\Proclaim{5.1.6. The case $d=1$.} \rm
We assume $d=1$. So $G^{\text{adnc}}_{V_0}$ is a $PSL_2$-group. Let $\bar G_{V_0}$ be the subgroup of $G_{V_0}$ generated by $\mu_{V_0}^\prime$ and by the normal, semisimple subgroup of $G^{\der}_{V_0}$ having $G^{\text{adnc}}_{V_0}$ as its adjoint. It is a $GL_2$-group. We also assume that its representation on $M_0^\prime$ is a direct sum of standard $2$-dimensional representations. As $G_{\dbQ_p}$ is split, the triple $(M_0^\prime,\phi,\bar G_{V_0})$ is a Shimura $\sigma$-crystal. So the abelian variety $A_v$ over $k(v)$ is either ordinary or supersingular; of course, if the Frobenius torus $T_v$ has rank at least $2$, then $A_v$ must be ordinary and so 4.4.8 applies. This complements [No, 2.2].
\finishproclaim 

We have the following generalization of 4.2.15.

\Proclaim{5.1.7. Exercise.} \rm
Show that in 4.2.7.2 we always have $T^{00}=T^0$. Hint: just copy Step 3 of 5.1.3.
\finishproclaim 

\smallskip
\Proclaim{5.2. A criterion.} \rm 
It is desirable to have a concrete approach for checking DD. In [Va7] we approach the proof of the equality $d=\dim_{\dbC}(X_A)$ via [Va2, 3.6] and the following Lemma and its variants of 5.2.1.
\finishproclaim

\Proclaim{5.2.0. Lemma.} We assume that 4.4.7 remains true if we replace $R$ by an integral, \'etale $V_0[y_1,...,y_d][{1\over {\prod_{i=1}^d y_i}}]$-algebra $R_0$, with $y_i$'s as independent variables. In other words, we assume that $l_R$ factors through a morphism $l_{R_0}\colon\Spec(R_0)\to\scrA_{g_A,1,4}$, with $R_0$ as mentioned, for which the following two things hold:

\medskip
-- the Frobenius $\Phi_R$ is obtained from a Frobenius $\Phi_{R_0}$ of the $p$-adic completion $R_0^\wedge$ of $R_0$ (via the factorization $\Spec(R)\to\Spec(R_0)$; due to the versal property in 4.4.7, it is automatically formally \'etale);

\smallskip
-- the principally quasi-polarized $p$-divisible object of $\scrM\scrF_{[0,1]}^\nabla(R_0)$ defined by the principally polarized abelian scheme $(A_{R_0},p_{A_{R_0}})$ over $R_0$ we get through $l_{R_0}$, is isomorphic to one of the form
$$\scrC_{R_0}:=(M_0^\prime\otimes_{V_0} R_0,F_0^\prime\otimes_{V_0} R_0,g_{\text{univ}}^0(\phi^{\prime}\otimes 1),\nabla^0,p_{A^0}),$$
with $g_{\text{univ}}^0\in G_{V_0}(R_0^\wedge)$ fixing $p_{A_0}$ and such that $\scrC_{R_0}$ is a versal deformation of $\scrC_0$. 

\medskip
We also assume that the $\dbQ_p$-\'etale Tate tensor of $A_{R_0^\wedge\fracwithdelims[]1p}$ corresponding via [Fa1, 2.6] to $t^0_{\alpha}\in\scrT(M_0^\prime\otimes_{V_0} R_0^\wedge\fracwithdelims[]1p)$ is in fact defined over $R_0\fracwithdelims[]1p$, $\forall\in\scrJ_0$; so we get a family (still denoted by) $(v_{\alpha})_{\alpha\in\scrJ_0}$ of $\dbQ_p$-\'etale Tate cycles of $A_{R_0\fracwithdelims[]1p}$. Then $d=\dim_{\dbC}(X_A)$.
\finishproclaim

\proof
$l_{R_0}$ factors through $\scrN(4)$, cf. Corollary of 4.4.7.1. As in the proof of 4.4.8 a), we get (via 4.4.9) the existence of a Zariski dense set $ZDS$ of $\dbF$-valued points of $\Spec(R_0/pR_0)$ such that:

\medskip\noindent
{\bf ORD.} {\it $y_0^*(A_{R_0})$ is an ordinary abelian variety and its canonical lift is the pull back of $A_{R_0}$ through a $V_0$-valued point of $\Spec(R_0)$ lifting $y_0$, $\forall y_0\in ZDS$.}

\medskip
So (cf. [Mo, 3.11 and 5.2]) $l_{R_0}$ dominates a Shimura subvariety of ${\scrA_{g_A,1,4}}_{O_{(v)}}$. This subvariety is a connected component $Y$ of $\scrN(4)$, cf. 4.4.7 a) and Corollary of 4.4.7.1.

Dropping the assumption that $R_0$ is smooth over $V_0$ of relative dimension $d$, we can assume that the affine $Y_{V_0}$-scheme $Y^0_{V_0}:=\Spec(R_0)$ is obtained from an integral, finite type $Y_{O_{(w)}}$-scheme $Y^0$, with $O_{(w)}$ a local, \'etale $\dbZ_{(p)}$-algebra containing $O_{(v)}$, by pull back. This is a general fact: it can always be achieved by passing (if needed) to a finite field extension $F$ of $E$ included in $K_0$ and by replacing $Y^0_{V_0}$ by a scheme of finite type over it (so $O_{(w)}:=F\cap V_0$). 

Using standard Galois descent, from 4.4.7 a) and from 4.3.6.2 we deduce that the family of $\dbQ_p$-\'etale Tate cycles $(v_{\alpha})_{\alpha\in\scrJ_0}$ is defined over the generic fibre $Y^0_F$ of $Y^0$. There is a quasi-finite, dominant morphism $Y^1_F\to Y_F$ factoring through the morphism $Y_F^0\to Y_F$ (argument: localizing, we can assume that $Y_F^0$ is smooth over $Y_F$; so we can use [EGA, 17.16.2]). 

We get a family, still denoted by $(v_{\alpha})_{\alpha\in\scrJ_0}$, of $\dbQ_p$-\'etale Tate cycles defined over $Y_F^1$.  
With this we are essentially done. 
There are many directions to be persuaded to end the proof.
We choose one.

We can assume $Y^1_F$ is geometrically connected: this is achieved by replacing $F$ with a finite field extension of it (we are not anymore interested in having $F$ as a subfield of $K_0$). Let $H_A^1$ be the maximal connected subgroup of $H_A$ with the property that $H^1_{A\dbQ_p}$ is included in $\tilde G_{\dbQ_p}$. We consider the extension $Y^1_{\dbC}\to Y_{\dbC}$ of the morphism $Y^1_F\to Y_F$ to $\dbC$ (via the inclusion $F\subset\dbC$).
The key point is: 

\medskip\noindent
{\bf DENSE.} {\it $Y_{\dbC}^1$ has an analytically dense set $ADS$ of complex points, over which we get (by pulling back $A_{R_0}$) abelian varieties having complex multiplication.}

\medskip
Argument: for $Y_{\dbC}$ this is a property of Shimura varieties and so we just need to point out that the morphism $Y_{\dbC}^1\to Y_{\dbC}$ is quasi-finite and dominant. $\forall\alpha\in\scrJ_0$, in all these points $v_{\alpha}$ becomes a $\dbQ_p$-linear combination of $p$-components of \'etale components of Hodge cycles, cf. 4.2.6. So all the Hodge structures on $W\otimes_{\dbQ} \dbC$ defined by these abelian varieties over $\dbC$, correspond to cocharacters of $G\Sp(W,\psi)_{\dbC}$ factoring through $H^1_{A\dbC}$. As the set $ADS$ is analytically dense in $Y^1_{\dbC}$, this remains true for all complex points of $Y^1_{\dbC}$. This implies (cf. the very definition of $H_{A}$) that $H^1_A=H_A$. So $H_{A\dbQ_p}^1=\tilde G_{\dbQ_p}=H_{A\dbQ_p}$ and so also $d=\dim_{\dbC}(X_A)$. This ends the proof.   
\finishproclaim

\Proclaim{5.2.1. Remark.} \rm
Let $G^{\text{dernc}}_{\dbQ_p}$ be the normal, semisimple subgroup of $G^{\text{der}}_{\dbQ_p}$ having $G^{\text{adnc}}_{\dbQ_p}$ as its adjoint. If by chance we only know that the subfamily of $(v_{\alpha})_{\alpha\in\scrJ_0}$ formed by tensors fixed by the centralizer of $G^{\text{dernc}}_{\dbQ_p}$ in $H_{A\dbQ_p}$ is defined over a scheme $Y^0_{V_0}$ as in 5.2, then the same arguments show that $G^{\text{adnc}}_{\dbQ_p}$ is a direct factor of $H_{A\dbQ_p}^{\ad}$. So $G_{K_0}^{\text{adnc}}=H_{AK_0}^{\text{adnc}}$. So $d=\dim_{\dbC}(X_A)$. Moreover, 4.4.8.3 points out that we have a variant of 5.2 in which $\scrJ_0$ is replaced by $\scrJ$ and $\tilde G_{\dbQ_p}$ is replaced by $G_{\dbQ_p}$.
\finishproclaim

\Proclaim{5.3. Serre's volumes.} \rm
We come back to the general situation described in 4.0: $A$ is an arbitrary abelian variety over a number field $E$ and $p$ is an arbitrary prime; in order to emphasize $p$, $\rho$ of 4.0 is now denoted by $\rho_p$. Let
$$\prod_{p\,\text{a prime}} \rho_p:\Gal(E)\to GL(H^1_{\text{\'et}}(A_{\Ebar},\dbA_f))(\dbA_f)$$ 
be the product $\dbA_f$-representation of all these $\rho_p$'s. 
Let $E^{\text{conn}}$ be the smallest finite field extension of $E$ such that all Hodge cycles of $A_{\Ebar}$ are defined over $E^{\text{conn}}$. It is a Galois extension of $E$: for any prime $p$, we have $\Gal(E^{\text{conn}})=\rho_p^{-1}(H_{A\dbQ_p}(\dbQ_p))$.

We now assume that the Mumford--Tate conjecture is true for $A$. 
This implies that 
$E^{\text{conn}}$ is the same field extension of $E$ as the one defined (see [Pi, 3.6]) by Serre. Two of the most important global invariants of the isogeny class of $A$ are the Shimura pair $(H_A,X_A)$ (see 4.3.6.1) and $E^{\text{conn}}$ (or just $[E^{\text{conn}}:E]$). 

We consider the Haar measure $\scrM(p)$ on $H_A(\dbQ_p)$ normalized by the fact that a maximal compact, open subgroup $H(p)$ of $H_A(\dbQ_p)$ of maximum volume, has volume 1. So if $H_A$ is unramified over $\dbQ_p$, then $H(p)$ is a hyperspecial subgroup of $H_A(\dbQ_p)$, cf. [Ti, p. 55]. We define Serre's (local) $p$-volume $\scrV_p(A)$ of $A$ by:
$$\scrV_p(A):=\scrM(p)(\rho(\Gal(E^{\text{conn}}))).$$
Similarly we define the total Serre's volume $\scrV(A)$ of $A$ by:
$$\scrV(A):=\prod_{p\,\text{a prime}}\scrV_p(A).$$
Obviously $\scrV_p(A)$ is a rational number of the interval (0,1]. So $\scrV(A)$ is a real number of the interval [0,1]. It is expected that $\scrV(A)$ is a rational number of (0,1] (see [Se, 11.4?]). Working with a Haar measure $\scrM(\dbA_f)$ of $H_A(\dbA_f)$, normalized in the same manner, we define the global Serre's volume of $A$ by:
$$\scrV_g(A):=\scrM(\dbA_f)(\prod_{p\,\text{a prime}} \rho_p(\Gal(E^{\text{conn}}))).$$
We have 
$$0\Le\scrV_g(A)\Le\scrV(A).$$
It is expected (see loc. cit.) that $0<\scrV_g(A)$ and that $\scrV(A)$ and $\scrV_g(A)$ are rational numbers. 
$\scrV_g(A)$ and, with $p$ a prime, $\scrV_p(A)$ are as well, in our opinion, very important global invariants of the isogeny class of $A$.
\finishproclaim

\Proclaim{5.3.1. Remarks.} \rm 
{\bf 1)} The study of Serre's volumes is closely related to the study of finite quotients of $\Gal(E)$. We expect the possibility of proving that any finite group (resp. any simple, finite group) $C$ is a quotient of $\Gal({\dbQ})$ by just using computations of local Serre's volumes of particularly chosen abelian varieties (resp. absolutely simple abelian varieties) over $\dbQ$. At least in the case when $C$ is a simple finite group associated to a Chevalley group (over a finite field) of classical Lie type, it should be possible to construct an epimorphism $\Gal(\dbQ)\twoheadrightarrow C$, starting from 5.1.2 and 3.2.1. We will come back to all these in [Va6]. 

{\bf 2)} We do believe that everything one would like to know about finite quotients of $\Gal(\dbQ)$ is codified in SSO. This represents our second approach towards the understanding of such quotients (see [Va2, 3.6.18.4.6 H] for the first one).

{\bf 3)} If the Zariski closure $H_{A\dbZ_p}$ of $H_{A\dbQ_p}$ in $GL(H^1)$ is a reductive group over $\dbZ_p$, then the study of the fact that $\scrV_p(A)$ is or is not 1, is closely related to the study of the fact that the natural (adjoint) representation of $H^{\der}_{A\dbZ_p}(\dbF_p)$ on $\Lie(H^{\der}_{A\dbZ_p})\otimes_{\dbZ_p} \dbF_p$ is or is not irreducible. In connection to the second study, see [Va3] for a quite elementary approach.

{\bf 4)} We think the terminology Serre's volumes is justified, as J.-P. Serre was the first one to study the numbers (representations) involved; for instance, see the report [Se] and its references.
\finishproclaim

\Proclaim{5.4. Extra remarks. 1)} \rm
 5.1.2 and 5.1.5 imply (this is standard) similar results for finitely generated fields of characteristic 0. Moreover (cf. the proof of 5.1.5), we can as well replace abelian varieties by abelian motives over such fields.

{\bf 2)} Many previous approaches (for instance see [Pi]) to the proof of 4.1 were intermingled with approaches to the proof of the ordinary reduction conjecture for abelian varieties over number fields. So a natural question arises: what the proofs of 5.1.1-2 and the new techniques of \S3-4 bring new in this direction? Our approach to the proofs of 5.1.1-2 was mainly based on Shimura varieties and much less on the ``inner" study of Frobenius tori (we needed mainly just their construction and 4.2.4-5). So here we do not bring anything new, besides more information on Frobenius tori (see 3.3.2, 4.4.8 b) and c) and 5.1.7). However, we invite the reader to deal with the following two Exercises (their complete proofs will be presented in [Va7]; warning: the first one is hard and so it is starred). 

\smallskip
{\bf Exercise 1*.} We consider the situation of 4.0. We assume $(H_A,X_A)$ is such that all simple symmetric domain factors of a connected component of $X_A$ are (see [He, p. 518]) either of type  $A III(q=1)$ or of type $BD I(q=2)$ (i.e. we assume that each simple, non-compact factor of $H_{A\dbR}^{\ad}$ is either $SU(n,1)^{\ad}_{\dbR}$ or $SO(n,2)_{\dbR}^{\ad}$, $n\in\dbN$). Using 5.1.2, [Pi] (the part referring to ordinariness), and the results on the Langlands--Rapoport conjecture stated in [Va1, 1.7.1], show that the ordinary reduction conjecture is true for $A$. Hint: use 5.1.3 S2, [Va2, 4.9.23] and 3.3.3. 

\smallskip
{\bf Exercise 2.} Using Exercise 1 and 5.1.3.2, show that the ordinary reduction conjecture is true for all motives over number fields mentioned in [An, 1.5]. Hint: use [An, 1.5.2].

\smallskip
{\bf 3)} The family of tensors of 4.3.4.8 a) is strongly $\dbZ_p$-very well position for $G_{\dbQ_p}$. This is a consequence of its proof and of [Va1, 4.3.6 2)] and the reference to it in [Va2, AE.0].

{\bf 4)} In [Va7] we will see that instead of $N_3(A)$ we can work as well with $N_3^\prime(A)\in\dbN\setminus\{1,2\}$ defined by: for any prime $q\Ge N_3^\prime(A)$, $H_A$ is unramified over $\dbQ_q$.

{\bf 5)} We now assume that $\End(A_{\Ebar})=\dbZ$ and that the simple, adjoint Shimura pair $(H_A^{\ad},X_A^{\ad})$ is of $D_n^{\dbR}$ type. Here $n\Ge 4$ has to be an odd number (see autodualities in [Pi, 4.6]). 5.1.2 allows another proof of [Pi, p. 230-236]. In other words, 3.2 f) represents another (direct) way (based on 5.1.2) to get directly the lifting property of the $l$-adic representations mentioned in the first paragraph of [Pi, p. 230] and due to J. -P. Wintenberger. 

{\bf 6)} From the point of view of possibilities of types of simple Shimura varieties of adjoint, abelian type, 5.1.2 covers about 80 percents (cf. 1.7.1). To our knowledge 5.1.2 does encompass all previous results pertaining to 4.1 except [Pi, 5.14]. Warning: one can not obtain a general form of [Pi, 5.14] (in the pattern of 5.1.2 C)) by just using the ideas of this paper; however, restricted to the case $\End(A_{\dbC})=\dbZ$, loc. cit. does have such a general form (i.e. to [Pi, 5.15] one can add the $C_n$ Lie types mentioned in 5.1.2 c) and c'); warning: the proof of loc. cit. can be entirely adapted to get this as well). See [Va5] for the general principle behind [Pi, 5.14] which encompasses loc. cit. and extends significantly this paper in the context when the simple, adjoint Mumford--Tate groups have over $\dbR$ only one simple, non-compact factor; [Va5] deals as well with some cases in which we have at least two such factor.

{\bf 7)} Shimura varieties of $A_n$ (resp. of either $B_n$ or $D_n^{\dbR}$) type are also called unitary (resp. orthogonal or spin Shimura varieties). 5.1.2 supports the point of view of [Va1, 1.2.6.2] that the study of unitary and spin Shimura varieties is easier than the one of Shimura varieties of $C_n$ or $D_n^{\dbH}$ type. Based on this point of view, on [Va1, 1.7.1], on 3.2.4 and on 5.3.1 2), we do expect that the mathematical role played by Siegel modular forms in the twentieth century will be substituted in this century by the one played by spin modular forms. 

{\bf 8)} For the remaining cases of the Mumford--Tate conjecture, we can assume $H_A^{\ad}$ is $\dbQ$--simple (and if it is helpful that 3.3.3 1), 2) or 3) holds), cf. 5.1.3-4.

{\bf 9)} For any $A_n$, $B_n$, $C_{2n+1}$ ($n\Ge 1$) or $D_n$ ($n\Ge 5$) Lie type $\tau$, there is a simple, adjoint Shimura pair $(H_0,X_0)$, with $H_0$ an absolutely simple group of $\tau$ Lie type, such that (with the notations of 4.0) the Mumford--Tate conjecture is true for $A$, provided $(H_A^{\ad},X_A^{\ad})$ is simple and $H_{A\dbR}^{\ad}$ has a factor isomorphic to $H_{0\dbR}$. Argument: for instance, we can take $H_0$ such that $H_0$ is $SU(1,n)^{\ad}$ or $SO(2,m)^{\ad}$, with $m\in\{2n-2,2n-1\}$, or a split, adjoint group of $C_{2n+1}$ Lie type; so 5.1.2 a), b) and e) apply. We also have a slightly weaker variant for the $D_4$ Lie type: we need to require $(H_A,X_A)$ to be of $D_4^{\dbR}$ type. But this is not yet known to be true for arbitrary $C_{2n}$ Lie types (cf. 6) and 5.1.2.4).    

{\bf 10)} From Subfact 2 of 5.1.3 S1, as $G_{\dbQ_p}$ is split, we get $k(v_A)=\dbF_p$ (cf. also [Mi3, 4.6-7]; for $v_A$ see 4.3.7).
This conforms as well with 4.4.9 and [Va2, 4.6 P1].

{\bf 11)} Many ideas used in this paper can be used as well for other classes of motives, provided (at least some of) the equivalent results of [De2, 2.3.9-13], 4.2.3, [Sa] and [Zi, 4.4] are known to be true (in some form). Also we would like to mention: in connection to 4.1, the results [De2, 2.3.9-13] and [Zi, 4.4] were neither used nor quoted previously.

{\bf 12)} Many of the ideas of 4.3-4 and of 5.2 can be extended to other contexts (for instance, those in which we have an abelian variety over $W(k)$, with $k$ a perfect field).
\finishproclaim

\Proclaim{5.5. Final complements.} \rm
\finishproclaim

\Proclaim{5.5.1. A method of extending the main result of [Wi] to some ramified contexts.} \rm
We use the notations of 2.4. We do not assume a priori the existence of the cocharacter $\mu_{K}$. It is desirable to have easily checkable conditions on the family of tensors $(w_{\alpha})_{\alpha\in\scrJ}$ which guarantee the existence of such a cocharacter. Some time ago we proposed one such condition:

\medskip
--  the family of tensors $(w_{\alpha})_{\alpha\in\scrJ}$ is moreover $H^1_{\text{\'et}}(H_V,\dbZ_p)$-representation well positioned for $\tilde G_{\dbQ_p}$.

\medskip
However, it turned out that this extra condition is not needed, as we have:

\medskip\noindent
{\bf Lemma.} {\it Under the assumptions on the family $(w_{\alpha})_{\alpha\in\scrJ}$ of $\dbQ_p$-\'etale Tate cycles of 2.4, the cocharacter $\mu_K$ exists automatically.}

\medskip
\proof
Using Tate's classical result [Tat, cor. 2 of p. 180] and standard arguments of algebraic geometry, we deduce that such a cocharacter exists over a finite field extension $K_1$ of $K$. Denoting by $V_1$ the ring of integers of $K_1$, as in 2.4 we get a cocharacter $\mu_{V_1}$ of $G_{V_1}$. We can assume that the residue field of $V_1$ is a Galois extension $k_1$ of $k$. As the special fibre of $\mu_{V_1}$ lifts to a cocharacter of $G_{W(k_1)}$ (cf. [SGA3, Vol. II, 3.6 of p. 48]), from the very definition [Va2, 2.2.8 1)] we get: the quadruple $(M_0\otimes_{W(k)} W(k_1),\phi\otimes 1,G_{W(k_1)},(t_{\alpha})_{\alpha\in\scrJ})$ is a Shimura $\sigma$-crystal over $k_1$; so, with the terminology of [Va2, 2.2.9 1')], $(M_0,\phi,G_{W(k)},(t_{\alpha})_{\alpha\in\scrJ})$ is a quasi Shimura $\sigma$-crystal over $k$. 

But the Galois group $\Gal(k_1/k)$ acts on the $F^1$-filtrations $F^1(M_0\otimes_{W(k)} W(k_1))$ of $M_0\otimes_{W(k)} W(k_1)$ such that the quadruple $(M_0\otimes_{W(k)} W(k_1),F^1(M_0\otimes_{W(k)} W(k_1)),\phi,G_{W(k_1)})$ is a Shimura filtered $\sigma$-crystal over $k_1$. As $F^1(M_0\otimes_{W(k)} W(k_1))$ mod $p$ is defined over $k$, by induction on $n$ we can assume that $F^1(M_0\otimes_{W(k)} W(k_1))$ mod $p^n$ is defined over $W_n(k)$. For instance, the space of such filtrations mod $p^2$ is a $k_1$-vector subspace of the extension to $k_1$ of the $k$-vector space parameterizing lifts of the kernel of $\phi$ mod $p$ to a direct summand of $M_0/p^2M_0$; so [Bo2, 14.2 of p. 30] applies. The induction follows.

So we can assume that $F^1(M_0\otimes_{W(k)} W(k_1))$ is the extension to $W(k_1)$ of a direct summand $F^1$ of $M$. So we can construct a Shimura $p$-divisible group $(H,(t_{\alpha}^0)_{\alpha\in\scrJ})$ over $R$ as in 2.4; we assume the Kodaira--Spencer map of $H$ is injective. So, as $(H_{V_1},(t_{\alpha}^V)_{\alpha\in\scrJ})$ is obtained from $(H_{W(k_1)},(t_{\alpha}^0)_{\alpha\in\scrJ})$ by pull back, $(H_{V},(t_{\alpha}^V)_{\alpha\in\scrJ})$ is as well obtained from $(H,(t_{\alpha}^0)_{\alpha\in\scrJ})$ by pull back. From this and from the fact that the canonical split cocharacter of $(M_0,F^1,\phi)$ factors through $G_{W(k)}$ (as it can be seen by passage from $k$ to $k_1$), the Lemma follows.
\finishproclaim

\Proclaim{5.5.2. A representation well positioned result for suitable Dedekind rings.} \rm
4.2.14 and [Bou2, ch. VIII, \S7.3] imply that the condition 4.3.4.2 iv) can be discarded in the context of 4.3: the numbers of 4.3.4.2 iv) for the representation $G_{\dbQ_p}\hookrightarrow GL(H^1\fracwithdelims[]1p)$ are all 2 (to be compared with [Va1, 5.7.1]). However, the conditions i) to vi) of 4.3.4.2 are such that they handle any arbitrary situation. 

Let $D$ be an integral Dedekind ring of characteristic 0 having an infinite number of maximal ideals. Let $K(D)$ be its field of fractions. We assume that for any $N\in\dbN$, there is only a finite number of prime ideals of $D$ of characteristic smaller than $N$. Let $M_D$ be a projective $D$-module of finite rank. Let $G_{K(D)}$ be a reductive subgroup of $GL(M_D\otimes_D K(D))$. We have:
\finishproclaim

\Proclaim{Theorem.}
There is a family of tensors $\scrF=(v_{\alpha})_{\alpha\in\scrJ_D}$ of $\scrT(M_D\otimes_D K(D))$ fixed by $G_{K(D)}$ and a non-zero ideal $I_D$ of $D$, such that for any maximal ideal $m$ of $D$ not containing $I_D$ and such that the Zariski closure of $G_{K(D)}$ in $GL(M_D\otimes_D D_{(m)})$ is reductive, $\scrF$ is enveloped by $M_D\otimes_D D_{(m)}$ and is strongly $M_D\otimes_D D_{(m)}$-representation well positioned for the group $G_{K(D)}$. Moreover, we can take $I_D$ to depend only on the representation of $Z^0(G_{K(D)})$ on $M_D\otimes_D K(D)$ (and not on $G^{\der}_{K(D)}$).
\finishproclaim

\proof
Let $n\in\dbN$ (resp. $n_{M_D}\in\dbN$) be obtained in the same way as $\max\{d(A),4\}$ (resp. as $N_2(A)$) in 4.3.4.2-6 but with $2g_A$ replaced by $\dim_D(M_D)$. Let $n_1\in\dbN$, $n_1\Ge n$, be such that $G_{K(D)}$ is the subgroup of $GL(M_D\otimes_D K(D))$ fixing a family of homogeneous tensors of $\scrT(M_D\otimes_D K(D))$ of degree at most $n_1$. We take $\scrF$ to be the family of tensors formed by all homogeneous tensors of $\scrT(M_D)$ fixed by $G_{K(D)}$ and of degree at most $n_1$. Let $I_D$ be the intersection of all maximal ideals of $D$ of whose residue fields have characteristic at most $n_{M_D}$. Let $I_D(G_{K(D)})$ be the largest subset of $D$ such that the Zariski closure of $G_{K(D)}$ in $M_D\otimes_D D_{(\tilde m)}$ is a reductive group $G_{D_{(\tilde m)}}$, for any maximal ideal $\tilde m$ of $D$ not containing $I_D(G_{K(D)})$. 

For the definition of $N_2(A)$ in 4.3.4.2 it was irrelevant that we are dealing with a split group over $\dbQ_p$. On the other hand, the only place of the proof of 4.3.4.8 b) where it was relevant that we are dealing with a split group over $\dbQ_p$ is the one in which we appealed to 4.3.4.6 1) and 4). But $G_{D_{(\tilde m)}^{\sh}}$ is split and so by using $D_{(\tilde m)}^{\sh}$-linear combinations of members of $\scrF$ we can still appeal to 4.3.4.6 1) and 4).

Moreover, in the proof of 3.4.8.4 b), the role of $\dbZ_p$ and $V$ was just one of a faithfully flat monomorphism $O_1\hookrightarrow O_2$ between two DVR's of mixed characteristic, with $O_2=O_2^{\sh}$. So, we can appeal to it, in our present context of $D_{(\tilde m)}\hookrightarrow D_{(\tilde m)}^{\sh}$, provided $\tilde m$ does not contain $I_D$. So the Theorem follows from 4.3.4.8 b) and 5.4 3).  

\bigskip
\noindent
{\boldsectionfont Erratum to [Va1]}
\bigskip

What follows next will be transferred to the final version of [Va2]. It is included here just for the time being (convenience of the reader).

The second paragraph $SP$ of [Va1, Case 2 of 6.6.5.1] is not correct as it stands. To explain this we use the notations of [Va1, 6.5.1 and 6.6.5]. Due to the fact that we took $W_{\dbZ_{(p)}}=K_{(p)}\otimes_{F_{(p)}} W_{(p)}$, the subgroup $G_{4\dbZ_{(p)}}$ of $GSp(W_{\dbZ_{(p)}},\psi)$ fixing all endomorphisms of $W_{\dbZ_{(p)}}$ fixed by $\tilde G_{3\dbZ_{(p)}}$, is reductive and has a derived subgroup which over $V_0$ is isomorphic to two copies of $G_{3V_0}^{\der}$ in such a way that the embedding of $G_{3V_0}^{\der}$ in $G^{\der}_{4\dbZ_{(p)}}$ is the diagonal embedding. Even if we replace $\tilde G_3$ by the smallest subgroup of it through which all homomorphisms ${\Res}_{\dbC/\dbR} \dbG_m\to G_{3\dbR}$ defining elements of $X_3$ factor, in general we can not ``get rid" of the second copy of $G_{3V_0}^{\der}$. 
\smallskip
There is a very simple way to adjust $SP$ so that we do get (with perhaps different notations) an injective map $\tilde f_3:(\tilde G_3,\tilde X_3)\hookrightarrow (GSp(W,\psi),S^0)$ such that $\tilde G_{3\dbZ_{(p)}}$ is the subgroup of $GSp(W_{\dbZ_{(p)}},\psi)$ fixing all endomorphisms of $W_{\dbZ_{(p)}}$ fixed by $\tilde G_{3\dbZ_{(p)}}$. It goes by: we entirely ``skip" the use of $K_{(p)}$ as follows.

We work with $W_{(p)}$ instead of $W_{\dbZ_{(p)}}$. $K_S$ is a totally imaginary quadratic extension of $F$ and so it makes sense to speak about ${K_S}_{(p)}$. The torus $GT:={\Res}_{{K_S}_{(p)}/\dbZ_{(p)}} \dbG_m$ acts naturally in a faithful way on $W_{(p)}$. Let $G_{4\dbZ_{(p)}}$ be the subgroup of $GL(W_{(p)})$ generated by $\tilde G_{\dbZ_{(p)}}$ and $GT$. Let $\tilde G_{4\dbZ_{(p)}}$ be the subgroup of $GL(W_{(p)})$ generated by $\tilde G_{\dbZ_{(p)}}^{\der}$, by $Z(GL(W_{(p)}))$ and by the maximal subtorus of $GT$ which over $\dbR$ is compact.

Let $RE$ be the set of real embeddings of $F$. For each $e_F\in RE$, let $V(e_F)$ be the maximal $\dbR$-vector subspace of $W_{(p)}\otimes_{\dbZ_{(p)}} \dbR$ on which the factor of $G_{4\dbR}^{\der}$ corresponding to $e_F$ acts non-trivially. So $GT_{\dbR}$ acts on it via its factor $F(e_F)$ which is a copy of ${\Res}_{\dbC/\dbR} \dbG_m$ and which is defined naturally by $e_F$. We have a direct sum decomposition
$$W_{(p)}\otimes_{\dbZ_{(p)}} \dbR=\oplus_{e_F\in RE} V(e_F)$$
left invariant by $G_{4\dbR}$.
Let $x\in X_2$. We consider a monomorphism $h_x:{\Res}_{\dbC/\dbR} \dbG_m\hookrightarrow GL(W_{(p)}\otimes_{\dbZ_{(p)}} \dbR)$ having the properties:

\medskip
\item{{\bf 1)}} if $e_F$ is such that $G^s\times_{\Spec(F)} {}_{e_F}\Spec(\dbR)$ is non-compact, then the resulting homomorphism ${\Res}_{\dbC/\dbR} \dbG_m\to GL(V(e_F))$ is the one obtained by composing the homomorphism ${\Res}_{\dbC/\dbR} \dbG_m\to G_{2\dbR}$ defining $x$ with the natural homomorphism $G_{2\dbR}\to GL(V(e_F))$;

\smallskip
\item{{\bf 2)}} if $e_F$ is such that $G^s\times_{\Spec(F)} {}_{e_F}\Spec(\dbR)$ is compact, then the resulting homomorphism ${\Res}_{\dbC/\dbR} \dbG_m\to GL(V(e_F))$ is a monomorphism whose image is naturally identified with the image of $F(e_F)$ in $GL(V(e_F))$.

\medskip
$h_x$ factors through $G_{4\dbR}$. The Hodge structure of $W_{(p)}\otimes_{\dbZ_{(p)}} \dbR$ it defines has type $\{(-1,0),(0,-1)\}$. So we can define two Shimura pairs $(G_4,X_4)$ and $(\tilde G_4,\tilde X_4)$ similar to $(G_3,X_3)$ and $(\tilde G_3,\tilde X_3)$ of [Va1, p. 507]. Moreover, taking perfect forms $\psi$ and $\tilde\psi$ on $W_{(p)}$ as mentioned in [Va1, v) of 6.5.1.1] and in SP, the subgroup $G_{5\dbZ_{(p)}}$ of $GSp(W_{(p)},\psi)$ fixing all endomorphisms of $W_{(p)}$ fixed by $\tilde G_{4\dbZ_{(p)}}$, is reductive and has $G^{\der}_{4\dbZ_{(p)}}=\tilde G^{\der}_{4\dbZ_{(p)}}$ as its derived subgroup. Warning: we do not have to replace $W_{(p)}$ by a direct sum of two copies of itself (as this is implicitly done by [Va1, 6.6.5 d1)]). If $X_5$ is the $G_{5\dbR}$-conjugacy class of $h_x$ viewed as a homomorphism of $G_{5\dbR}$, then the pair $(G_5,X_5)$ is a Shimura pair whose adjoint is $(G,X)$. Moreover, we get a PEL type embedding
$$f_5:(G_5,X_5)\hookrightarrow (GSp(W_{(p)}[{1\over p}],\psi),S_0)$$
which has the desired property (i.e. we can take $\tilde f_3:=f_5$).

\medskip\noindent
{\bf Remark.} In [Va1, 6.5.1.1 and Case 1 of 6.6.5.1], as above we can ``get rid" of $K$ for the $D_{2l+1}^{\dbR}$ type. But this is not true in general for the $B_l$, $C_l$, $D_{2l}^{\dbR}$ and $D_{l+2}^{\dbH}$ types. Here $l\in\dbN\setminus\{1\}$. However, we can always choose $W_{(p)}$ such that we can still ``get rid" of $K$, cf. 3.1 I') and 3.2.1.

\bigskip
\references{37}
{\nspace{

\Ref[An]
Y. Andr\'e,
\sl On the Shafarevich and Tate conjectures for hyperk\"ahler varieties,
\rm Math. Ann. {\bf 305} (1996), p. 205--248.

\Ref[BBM]
P. Berthelot, L. Breen, and W. Messing,
\sl Th\'eorie de Dieudonn\'e crystalline II,
\rm LNM {\bf 930} (1982), Springer--Verlag.

\Ref[Bl]
D. Blasius,
\sl A $p$-adic property of Hodge cycles on abelian
varieties, 
\rm Proc. Sympos. Pure Math. {\bf 55} (1994), Part 2, p. 293--308.

\Ref[BHC]
A. Borel and Harish-Chandra,
\sl Arithmetic subgroups of algebraic subgroups,
\rm Ann. of Math. (3) {\bf 75} (1962), 485--535.

\Ref[BM]  P. Berthelot and W. Messing, 
\sl Th\'eorie de Dieudonn\'e crystalline III, 
\rm Grothendieck's Festschrift I, Birkh\"auser, Progr. in Math. {\bf 86} (1990), p. 173--247.

\Ref[Bog]
F. A. Bogomolov,
\sl Sur l'algebricit\'e des repr\'esentations $l$-adiques,
\rm C. R. Acad. Sc. Paris Ser I Math. {\bf 290} (1980), p. 701--703.

\Ref[Bo1]
A. Borel,
\sl Properties and linear representations of Chevalley groups,
\rm LNM {\bf 131} (1970), p. 1--55, Springer--Verlag.

\Ref[Bo2]
A. Borel,
\sl Linear algebraic groups,
\rm Grad. Texts Math. {\bf 126} (1991), Springer--Verlag.

\Ref[Bou1] N. Bourbaki, 
\sl Groupes et alg\`ebres de Lie, 
\rm Chapitre {\bf 4--6} (1968), Hermann.

\Ref[Bou2]
N. Bourbaki,
\sl Groupes et alg\`ebres de Lie,
\rm Chapitre {\bf 7--8} (1975), Hermann.

\Ref[Chi]
W. Chi,
\sl $l$-adic and $\lambda$-adic representations associated to abelian varieties defined over number fields,
\rm Amer. J. Math. {\bf 114} (1992), p. 315--353.

\Ref[De1]
P. Deligne,
\sl Travaux de Shimura,
\rm S\'eminaire  Bourbaki 389, LNM {\bf 244} (1971), p. 123--163, Springer--Verlag.

\Ref[De2]
P. Deligne,
\sl Vari\'et\'es de Shimura: Interpr\'etation modulaire, et
techniques de construction de mod\`eles canoniques,
\rm Proc. Sympos. Pure Math. {\bf 33} (1979), Part 2, p. 247--290.

\Ref[De3]
P. Deligne,
\sl Hodge cycles on abelian varieties,
\rm Hodge cycles, motives, and Shimura varieties, LNM {\bf 900} (1982), p. 9--100, Springer--Verlag.

\Ref[DM]
P. Deligne and J. Milne, 
\sl Tannakian categories,
\rm Hodge cycles, motives, and Shimura varieties, LNM {\bf 900} (1982), p. 101--128, Springer--Verlag.

\Ref[EGA]
A. Grothendieck et al.,
\sl Etude locale des sch\'emas et des morphismes de sch\'ema,
\rm Publ. Math. IHES {\bf 20} (1964), {\bf 24} (1965), {\bf 28} (1966), {\bf 32} (1967).

\Ref[Fa1]
G. Faltings,
\sl Endlichkeitss\"atze f\"ur abelsche Variet\"aten \"uber Zahlk\'orpern,
\rm Inv. Math. {\bf 73} (1983), p. 349--366.

\Ref[Fa2]
G. Faltings,
\sl Crystalline cohomology and $p$-adic  Galois
representations, 
\rm Alg. Analysis, Geometry, and Nr.
Theory (1990), p. 25--79, Johns Hopkins Univ. Press.

\Ref[Fa3]
G. Faltings,
\sl Integral crystalline cohomology over very ramified
valuation rings,
\rm J. of Am. Math. Soc., Vol. {\bf 12} (1999), no. 1, p. 117--144.

\Ref[Ha]
G. Harder,
\sl \"Uber die Galoiskohomologie halbeinfacher Matrizengruppen II,
\rm  Math. Z. {\bf 92} (1966), p. 396--415.

\Ref[Hu1]
J. E. Humphreys,
\sl Linear algebraic groups,
\rm Grad. Texts Math. {\bf 21} (1975), Springer--Verlag.

\Ref[Hu2] 
J. E. Humphreys, 
\sl Conjugacy classes in semisimple algebraic groups, 
\rm Math. Surv. and Monog., Vol. {\bf 43} (1995), Am. Math. Soc. 

\Ref[He]
S. Helgason,
\sl Differential geometry, Lie groups, and symmetric spaces,
\rm Academic Press (1978).

\Ref[Ja] 
J.C. Jantzen, 
\sl Representations of algebraic groups, 
\rm Academic Press (1987).

\Ref[Ko]
R. E. Kottwitz,
\sl Points on some Shimura Varieties over finite fields,
\rm J. of Am. Math. Soc., Vol. {\bf 5} (1992), no. 2, p. 373--444.

\Ref[La]
S. Lang,
\sl Algebraic number theory,
\rm second edition, Grad. Texts Math. {\bf 110} (1994), Springer--Verlag.

\Ref[LP1]
M. Larsen and R. Pink,
\sl Determining representations from invariant dimensions,
\rm Inv. Math. {\bf 102} (1990), p. 377--398.

\Ref[LP2]
M. Larsen and R. Pink,
\sl Abelian varieties, $l$-adic representations and $l$-independence, 
\rm Math. Ann. {\bf 302} (1995), p. 561--579.

\Ref[Me]
W. Messing,
\sl The crystals associated to Barsotti--Tate groups
with applicactions to abelian schemes,
\rm LNM {\bf 264} (1972), Springer--Verlag.

\Ref[Mi1]
J. S. Milne, 
\sl The points on a Shimura variety modulo a prime of good
reduction,
\rm The Zeta function of Picard modular surfaces, Les Publ. CRM (1992), p. 153--255.

\Ref[Mi2]
J. S. Milne,
\sl Motives over finite fields,
\rm Proc. Symp. Pure Math. {\bf 55} (1994), Part 1, p. 401--459.

\Ref[Mi3]
J. S. Milne,
\sl Shimura varieties and motives,
\rm Proc. Symp. Pure Math. {\bf 55} (1994), Part 2, p. 447--523.

\Ref[Mi4]
J. S. Milne,
\sl On the conjecture of Langlands and Rapoport,
\rm manuscript 9/1995.

\Ref[Mo]
B. Moonen,
\sl Special points and linearity properties of Shimura varieties,
\rm thesis (1995), Utrecht University. 

\Ref[M-B] 
L. Moret-Bailly, 
\sl Pinceaux de vari\'et\'es ab\'eliennes, 
\rm J. Ast\'erisque {\bf 129} (1985).

\Ref[MFK]
D. Mumford, J. Fogarty and F. Kirwan,
\sl Geometric invariant theory,
\rm third enlarged edition, A Series of Modern Surv. in Math. {\bf 34} (1994), Springer--Verlag

\Ref[Mu1]
D. Mumford, 
\sl Families of abelian varieties,
\rm Algebraic Groups and Discontinuous Subgroups, Proc. Sympos. Pure Math., Am. Math. Soc., Vol. {\bf 9} (1966), p. 347--352, Providence, R. I.

\Ref[Mu2]
D. Mumford,
\sl Abelian varieties,
\rm Tata Inst. of Fund. Res. (1988), Oxford Press.

\Ref[No]
R. Noot,
\sl Abelian varieties, Galois representations and properties of ordinary reduction,
\rm Comp. Math. {\bf 97} (1995), p. 161--171.

\Ref[Oo]
F. Oort,
\sl Some questions in algebraic geometry,
\rm preprint Utrecht University, 6/95.

\Ref[Pi]
R. Pink,
\sl $l$-adic algebraic monodromy groups, cocharacters, and the Mumford--Tate conjecture,
\rm J. reine angew. Math. {\bf 495} (1998), p. 187--237.

\Ref[RR]
M. Rapoport and M. Richartz,
\sl On the classification and the specialization of $F$-crystals with additional structure,
\rm Comp. Math. {\bf 103} (1996), p. 153--181.

\Ref[Sa]
I. Satake,
\sl Holomorphic imbeddings of symmetric domains into a Siegel space,
\rm  Am. J. Math. {\bf 87} (1965), p. 425--461.

\Ref[Se]
J.-P. Serre, 
\sl Propri\'et\'es conjecturales des groupes de Galois motiviques et des repr\'esentations $l$-adiques,
\rm Proc. Sympos. Pure Math., Vol. {\bf 55} (1994), Part 1, p. 377--400.

\Ref[SGA3]
M. Demazure, A. Grothendieck, et al.,
\sl Schemes en groupes, 
\rm Vol. I--III, LNM {\bf 151--153} (1970), Springer--Verlag.

\Ref[SGA4]
M. Artin, A. Grothendieck and J. Verdier,
\sl Th\'eorie des topos et cohomologie \'etale des sch\'emas,
\rm Vol. I--III, LNM {\bf 269-70, 305} (1972-3), Springer--Verlag.

\Ref[Sh]
G. Shimura,
\sl Moduli of abelian varieties and number theory, Proc. Sympos. Pure Math., Vol. {\bf 9} (1966), Am. Math. Soc., p. 312--322, Providence, R. I.

\Ref[Ta1]
S. G. Tankeev,
\sl Cycles on abelian varieties and exceptional numbers,
\rm Izv. Ross. Akad. Nauk. Ser. Mat. {\bf 60} (1996), no. 2, p. 391--424.

\Ref[Ta2] 
S. G. Tankeev,
\sl On the weights of an $l$-adic representation and the arithmetic of Frobenius eigenvalues,
\rm (Russian) Izv. Russ. Akad. Nauk
Ser. Mat. {\bf 63} (1999), no. 1, p. 185--224, translation in Izv. Math. {\bf 63} (1999), no. 1, p. 181--218.

\Ref[Tat]
J. Tate,
\sl p-divisible groups,
\rm Proceedings of a conference on local fields, Driesbergen 1966 (1967), p. 158--183, Springer--Verlag.

\Ref[Ti]
J. Tits,
\sl Reductive groups over local fields, 
\rm Proc. Sympos. Pure Math. {\bf 33} (1979), Part 1, p. 29--69.

\Ref[Va1]
A. Vasiu,
\sl Integral canonical models for Shimura varieties of preabelian type,
\rm Asian J. Math., Vol. {\bf 3} (1999), no. 2, p. 401--518.

\Ref[Va2]
A. Vasiu,
\sl Points of integral canonical models of Shimura varieties of preabelian  type, p-divisible groups, and applications, 
\rm first part, submitted for publication in J. Ast\`erisque; up dated version 8/16/01, p. 1--595.

\Ref[Va3]
A. Vasiu,
\sl Letter to Prof. H. W. Lenstra, 11/16/98,
\rm www.math.arizona.edu/$\tilde{}$ adrian (to be submit. for publ. in an enlarged form).

\Ref[Va4]
A. Vasiu,
\sl Moduli schemes and the Shafarevich conjecture (the arithmetic case) for pseudo-polarized $K3$ surfaces,
\rm manuscript, 9/99.

\Ref[Va5]
A. Vasiu,
\sl Shimura varieties and the Mumford--Tate conjecture, part II,
\rm first announcement, 2/2000, www.math.arizona.edu/$\tilde{}$ adrian, to be submit. for publ.

\Ref[Va6]
A. Vasiu,
\sl Serre's volumes and finite quotients of the $\Gal(\dbQ)$,
\rm to be submit. for publ.

\Ref[Va7]
A. Vasiu,
\sl Points of the integral canonical models of Shimura varieties of  preabelian  type, p-divisible groups, and applications, 
\rm the second and the third parts, to be submit. for publ.

\Ref[Za]
Y. G. Zarhin,
\sl Abelian varieties having a reduction of K3 type,
\rm  Duke Math. J. {\bf 5} (1992), p. 511--527.

\Ref[Zi]
T. Zink,
\sl Isogenieklassen von Punkten von Shimuramannigfaltigkeiten mit Werten in einem endlichen K\"orper,
\rm Math. Nachr. {\bf 112} (1983), p. 103--124.

\Ref[Wi]
J.-P. Wintenberger,
\sl Un scindage de la filtration de Hodge pour certaines variet\'es algebriques sur les corps locaux,
\rm Ann. of Math. (2) {\bf 119} (1984), p. 511--548.

}}

\bigskip
\line{\hfill\vbox{\hbox{\it address:}
\medskip
\hbox{Adrian Vasiu}
\hbox{University of Arizona}
\hbox{Department of Mathematics}
\hbox{617 N. Santa Rita}
\hbox{P.O. Box 210089}
\hbox{Tucson, AZ-85721-0089, U.S.A.}
\medskip
\hbox{e-mail: adrian at math.arizona.edu\ \ \ \ \ }}}

\enddocument